\documentclass[preprint, review, 12pt]{elsarticle}

\usepackage{amssymb}
\usepackage{fullpage}
\usepackage{amsthm}
\usepackage{amsmath}
\usepackage{caption}
\usepackage{subcaption}
\usepackage{float}
\usepackage{hyperref}
\usepackage{bm}
\usepackage{color}

\journal{Journal Name}

\begin{document}

\begin{frontmatter}



\title{3rd-order Spectral Representation Method: Part I -- Multi-dimensional random fields with fast Fourier transform implementation}


\author{Lohit Vandanapu}
\author{Michael D. Shields}

\address{Dept.\ of Civil Engineering, Johns Hopkins University}

\begin{abstract}
This paper introduces a generalised 3rd-order Spectral Representation Method for the simulation of multi-dimensional stochastic fields with asymmetric non-linearities. The simulated random fields satisfy a prescribed Power Spectrum and Bispectrum. The general $d$-dimensional simulation equations are presented and the method is applied to simulate 2D and 3D random fields. The differences between samples generated by the proposed methodology and the existing classical Spectral Representation Method are analysed. An important feature of this methodology is that the formula can be implemented efficiently with the Fast Fourier Transform, details of which are presented. Computational savings are shown to grow exponentially with dimensionality as a testament of the scalability of the simulation methodology.
\end{abstract}

\begin{keyword}
Stochastic Fields \sep Random Fields \sep Spectral Representation \sep Fast Fourier Transform \sep Simulation


\end{keyword}

\end{frontmatter}


\section{Introduction}
\label{S:1}

Stochastic processes and random fields are used extensively in the field of engineering, from studying the dynamics of wind \cite{gurley1997analysis, Shields2013_1}, ocean waves \cite{elgar1985observations, Gurley1999}, and seismic loads \cite{deodatis1996non} on structures to simulation of material microstructures \cite{feng2014statistical,shields2016discussion}. Because of their importance, numerous methods have been developed for the simulation of stochastic processes and random fields. Simulation is particularly useful in the context of Monte Carlo simulations of large, complex non-linear systems where analytical analysis of the uncertainty in the system is not possible. Moreover, simulation of stochastic processes and random fields finds applications beyond simple Monte Carlo simulations and is important for essentially any simulation-based uncertainty quantification framework. 

Until recently, simulation methods for stochastic processes and random fields have derived only from second-order properties of the process or field. Consider a standard probability space $(\Omega, \mathcal{F},\mathcal{P})$ where $\Omega$ is the sample space, $\mathcal{F}$ the sigma algebra of events, and $\mathcal{P}$ a probability measure. In these simulation methods, the process/field indexed on $\bm{x}\in \mathcal{D}$ is represented in terms of a stochastic expansion of the general form
\begin{equation}
\label{eq:stochastic_expansion}
A(\bm{x}, \omega) \approx \hat{A}(\bm{x}, \omega) = \sum_{i=1}^{n}\theta_{i}(\omega)\psi_{i}(\bm{x})
\end{equation}
where $\theta_{i}(\omega), \omega\in\Omega$ are independent random variables and $\psi_{i}(\boldsymbol{x}), \bm{x}\in\mathcal{D}$ are deterministic basis functions. Many such stochastic expansions have developed. Among these methods the most popular ones are the Spectral Representation method(SRM)\cite{Shinozuka1972_1,Shinozuka1972_2,Shinozuka1991} and the Karhunen-Loeve Expansion(KLE)\cite{Huang2001,Ghanem1991}. Each of these methods operates by finding a set of random variables $\theta_{i}(\omega)$ along with a set of compatible basis functions $\psi(\bm{x})$ which satisfy $C(\boldsymbol{x_{1}}, \boldsymbol{x_{2}}) = \mathbb{E}[A(\boldsymbol{x_1})A(\boldsymbol{x_2})] = \mathbb{E}[\hat{A}(\boldsymbol{x_1})\hat{A}(\boldsymbol{x_2})]$.

For the SRM method, $\psi_{i}(\boldsymbol{x})$ are the harmonic functions (Fourier basis) and $\theta_{i}(\omega)$ are random variables whose amplitude is derived from the power spectrum (Fourier transform of the covariance function $C(\boldsymbol{x_{1}}, \boldsymbol{x_{1}})$). Likewise for the K-L expansion, $\psi_{i}(\boldsymbol{x})$ are the eigen-functions of the covariance function and $\theta_{i}(\omega)$ are standard normal random variables scaled by the square root of the appropriate eigenvalues.

While each of these methods has its advantages, all such methods have a common disadvantage in that they are only second-order representative, i.e they can only match the process up to its covariance function. Unless acted upon by a nonlinear operator, these fields are asymptotically Gaussian as the number of terms $n$ increases \cite{Grigoriu2010}. In signal processing terms, the simulated stochastic processes and random fields by the above methods are equivalent to the output of a linear system acted upon by Gaussian random noise. This simplification breaks down in case of real world scenarios such as seismic waves propagating through different strata of soil, non-linear wind loads on structures, ocean waves acting on an off-shore structural system, or turbulent flow of a fluid governed by the Navier-Stokes equation. Thus, the second-order representation inherently limits these methods as they fail to match the higher order properties of the stochastic fields, which dominate the tail behaviour and in turn plays a crucial role in uncertainty quantification, reliability etc. The stochastic fields generated from these non-linear systems possess asymmetric non-linear wave interactions which need to be included in the stochastic expansion, details of which were first introduced in \cite{Shields2017} and are reviewed in the subsequent sections. 

Methods for the simulation of non-Gaussian stochastic fields work primarily by non-linear transformation of the stochastic expansion Eq.\ \eqref{eq:stochastic_expansion}]. One class of such nonlinear transformations works by introducing correlated random variables with deterministic basis functions such as Hermite and Legendre polynomials \cite{Puig2002,Liu2017}. These stochastic processes match the marginal statistical moments to a certain order along with the covariance function. Perhaps the most commonly used method is the explicit Cumulative Distribution Function (CDF) based transformation \cite{Grigoriu2002} given by
\begin{equation}
    Y(x) = F^{-1}(\Phi(A(x)))
\end{equation}
where $A(x)$ is a standard Gaussian random process, $\Phi(.)$ is the standard normal CDF and $F^{-1}(.)$ is the inverse CDF of the prescribed non-Gaussian distribution. This method is generally referred to as the `translation process'. Efficient algorithms for the translation of scalar, vector, stationary, and non-stationary stochastic processes simulated by either SRM or KLE method have been developed in recent years\cite{Shields2011, Shields2013_1,Kim2015,Shields2013_2}. Another class of methods for simulation of non-Gaussian stochastic processes are the polynomial-chaos expansion methods \cite{Sakamoto2002}. Wavelet-based simulation methodologies have been developed and applied extensively in the case of non-stationary stochastic processes \cite{zeldin1996random,Phoon2004}.

This work is concerned with the efficient simulation of multi-dimensional non-Gaussian random fields. We specifically consider third-order, asymmetrically non-linear random fields (i.e.\ fields that possess quadratic phase interactions leading to an asymmetrically skewed distribution) prescribed by a known power spectrum and bispectrum. This extends the generalized third-order spectral representation method proposed in \cite{Shields2017} to multiple spatial dimensions and introduces a fast Fourier transform (FFT) implementation of the simulation algorithm that greatly improves the computational efficiency.

As a brief note, stochastic processes and random fields here are considered probabilistically equivalent with the only difference being that stochastic processes are indexed on time (one-dimensional, denoted by $t$ or $\tau$, with $\omega$ representing frequency under FFT) while random fields are indexed in space (up to three-dimensional, denoted by $x$ or $\xi$, with $\kappa$ representing wave-number under FFT). Given that this work focuses on simulating multi-dimensional quantities, we will generally present concepts in the context of random fields (using the $x$, $\xi$, $\kappa$ notation).


\section{Properties of Random Fields}
\label{S:2}
Prior to introducing any concepts in simulation, it is important first to review several important properties of random fields. First, we discuss the notion of stationarity (of various orders). We then review generalized moments, cumulants, and spectra for stationary random fields.

\subsection{Stationary Random Fields}
\label{S:2.2}

In numerous fields of engineering, we encounter random fluctuations that are probabilistically invariant under a translation in space, time etc. Random fields that are invariant across the indexing variable are referred to as `stationary'. For the development of the proposed methodology we present 3 different orders of stationary random fields.

\subsubsection{Strictly or Strongly Stationary Random Fields}

A random field $A(x)$ is said to be strictly stationary, or strongly stationary, if the full joint probability measure is invariant to a shift in index. That is, consider that $A(x)$ has complete $n$-dimensional joint cumulative distribution function given by $F_A(a(x_1),a(x_2),\dots,a(x_n))$, the random field is strongly stationary if
\begin{equation}
    F_A(a(x_1),a(x_2),\dots,a(x_n))=F_A(a(x_1+\xi),a(x_2+\xi),\dots,a(x_n+\xi)),\quad \forall \xi
    \label{eqn:strong_stationary}
\end{equation}
It follows from Eq.\ \eqref{eqn:strong_stationary} that all lower-dimensional distributions are similarly invariant to a shift in index and that all characteristics of the joint distribution (i.e.\ moments, cumulants, etc.) are independent of $\xi$.





\subsubsection{$k^{th}$-order Stationary Random Fields}

Strict stationarity is, as the name implies, is a very strong condition that can be assured only in very special cases (e.g.\ white noise, stationary Gaussian random fields). It is useful therefore, to introduce less strict conditions of stationarity.

A random field $A(x)$ is considered $k^{th}$-order stationary if the probability measure up to $k^{th}$-order is invariant to a shift in index. Let $F_A^{(k)}(a(x_1),a(x_2),\dots,a(x_k)), k<n$ denote the $k^{th}$-order joint cumulative distribution function of $A(x)$. The random field is $k^{th}$-order stationary if
\begin{equation}
    F_A^{(k)}(a(x_1),a(x_2),\dots,a(x_k))=F_A^{(k)}(a(x_1+\xi),a(x_2+\xi),\dots,a(x_k+\xi)),\quad \forall \xi
    \label{eqn:k_stationary}
\end{equation}
Again, it follows that all measures of order lower than $k$ are similarly invariant to a shift in index and that all characteristics of the $k^{th}$-order joint distribution are independent of $\xi$.

Of particular interest here is the notion of the $3^{rd}$-order stationarity, which implies that the bispectrum (defined below in Section \ref{S:3.1}) is invariant (i.e.\ does not change in time or space). As we will see, random fields generated according to the proposed method are $3^{rd}$-order stationary.




\subsubsection{Weak or Wide-Sense Stationary Random Fields}

A random field is considered to be weakly, or wide-sense stationary if the joint probability distribution up to $2^{nd}$-order is invariant to a shift in index. In other words, a weakly stationary random field is a $k^{th}$-order random field with $k=2$.

Weakly stationary random fields form a particularly important class of random fields because, for practical reasons, it is often only possible to ensure $2^{nd}$-order stationarity. More importantly here, existing simulation methods, rely heavily on the weak stationarity of the simulated random field. That is, because existing expansions are derived from $2^{nd}$-order properties of the random field (i.e.\ power spectrum or two-point correlation), the simulated fields are, by construction, $2^{nd}$-order stationary. Even methods that simulate non-stationary random fields (e.g. spectral representation using evolutionary spectra \cite{Shields2013_2} and Karhunen-Loeve expansion \cite{Kim2015}), often rely on the notion of instantaneous second-order stationarity. That is, at any point in time/space, the random field is $2^{nd}$-order stationary. The second-order properties, however, may vary in space/time. See Priestley's definition of the evolutionary spectrum for further details \cite{Priestley1965}.






\subsection{Moments and Cumulants} 
\label{S:2.1}
Only when a random field is Gaussian, and in a few other special cases, is the full joint probability density of the random field known. For practical purposes, many random fields are therefore defined through some subset of properties of the field -- typically its moments, cumulants, or spectra. These properties are reviewed in the following. We note however that, given the classical moment problem, the full probability density of the random field is identifiable from the moments only when Carleman's Condition is satisfied -- that is only when the infinite moment series has positive radius of convergence \cite{akhiezer1965classical}. Consequently, moments (cumulants/spectra) provide only a limited view of the random field since we realistically cannot expect to know infinite moments, nor can we be assured that Carleman's Condition will hold.


\subsubsection{Random Variables and Random Vectors}

Given a real random variable $X$, the $r^{th}$-order moments ($m_r$) and cumulants ($c_r$) are given by
\begin{equation}
m_{r} = \mathbb{E}[X^{r}] = \frac{\partial^{r} \phi(\lambda)}{{\partial \lambda}^{r}} \bigg\rvert_{\lambda=0}
\end{equation}
\begin{equation}
c_{r} = \frac{\partial^{r} \ln \phi(\lambda)}{{\partial \lambda}^{r}}\bigg\rvert_{\lambda=0}
\end{equation}
where $\phi(\lambda)$ is the characteristic function. This definition naturally extends to a random vector $\mathbf{X}=\{X_{1}, X{2}, \dots X_{n}\}$, where the $r=k_1+k_2+\dots+k_n$-th order joint moments and cumulants are given by \cite{Brillinger1965}
\begin{equation}
m_{k_{1}, k_{2} \dots k_{n}} = \frac{\partial^{r} \phi(\lambda_{1}, \lambda_{2} \dots \lambda_{n})}{{\partial \lambda_{1}}^{k_{1}} \lambda_{2}^{k_{2}} \dots \lambda_{n}^{k_{n}}} \bigg\rvert_{\lambda_{1}=\lambda_{2}=\dots=\lambda_{n}=0}
\end{equation}
\begin{equation}
c_{k_{1}, k_{2} \dots k_{n}} = \frac{\partial^{r} \ln \phi(\lambda_{1}, \lambda_{2} \dots \lambda_{n})}{{\partial \lambda_{1}}^{k_{1}} \lambda_{2}^{k_{2}} \dots \lambda_{n}^{k_{n}}}\bigg\rvert_{\lambda_{1}=\lambda_{2}=\dots=\lambda_{n}=0}
\end{equation}
where $\phi(\lambda_{1}, \lambda_{2} \dots \lambda_{n})$ is the joint characteristic function of $\mathbf{X}$.

In general, the $r=n^{th}$-order cumulants can be expressed in terms of the moments through the following relationships \cite{nikias1993higher}:
\begin{align}
\label{eq:1_5}
\begin{split}
c_{k_1,k_2,\dots,k_n} = c[X_1^{k_1},X_2^{k_2},\dots,X_n^{k_n}] = \sum (-1)^{p-1} (p-1)! E\big[\prod_{i\in s_1} X_i\big] E\big[\prod_{i\in s_2} X_i\big] \dots E\big[\prod_{i\in s_p} X_i\big] \\
\end{split}
\end{align}
where the summation extends over all groups $\{s_1,s_2,\dots,s_p\}, \ p=1,2,\dots,n$ of the integers $1, 2,\dots, n$. For example, some third-order cumulants are given by
\begin{align}
\label{eq:1_6}
\begin{split}
c_{111} = c[X_1,X_2,X_3] & = E[X_1 X_2 X_3] - E[X_1 ] E[X_2 X_3] -E[X_2] E[X_1 X_3] - E[X_3] E[X_1 X_2]\\ 
& + 2 E[X_1] E[X_2] E[X_3]\\
c_{120} = c[X_1,X_2^2] & = c[X_1,X_2,X_2]=E[X_1 X_2^2] - E[X_1 ] E[X_2^2] -2 E[X_2] E[X_1 X_2] \\ 
& + 2 E[X_1] E[X_2]^2\\
c_{300} = c[X_1^3] & = c[X_1, X_1, X_1] =  E[X_1^3] - 3E[X_1 ] E[X_1^2] + 2 E[X_1]^3 \\
\end{split}
\end{align}
Note that when $\mathbf{X}$ is jointly Gaussian, all cumulants $c_r$ of order $r>2$ are zero. Thus, non-zero higher order cumulants imply non-Gaussianity.

\subsubsection{Random Fields}
\label{S:2.3}

We can similarly denote the $n^{th}$-order moments for any stationary random field $A(x)$ by
\begin{equation}
m_{n}^{A}(\xi_{1}, \dots, \xi_{n-1}) = \mathbb{E}[A(x)A(x+\xi_{1}) \dots A(x+\xi_{n-1})].
\end{equation}
The cumulants of a stationary random field, meanwhile, can be expressed in terms of the moments by applying Eq.\ \eqref{eq:1_5}, yielding
\begin{align}
\label {eq:1_8}
\begin{split}
c_1^A &= m_1^A \\
c_2^A (\xi) &= m_2^A(\xi) - m_1^A \\
c_3^A (\xi_1, \xi_2) &= m_3^A(\xi_1, \xi_2) - m_1^A [m_2^A(\xi_1) + m_2^A(\xi_2) + m_2^A(\xi_2 - \xi_1)] + 2 (m_1^A)^3  \\
c_4^A (\xi_1, \xi_2, \xi_3) &= m_4^A(\xi_1, \xi_2, \xi_3) -m_2^A(\xi_1) m_2^A(\xi_3 -\xi_2) - m_2^A(\xi_2) m_2^A(\xi_3 -\xi_1) - m_2^A(\xi_3)m_2^A(\xi_2 - \xi_1) \\
& - m_1^A[m_3^A(\xi_2 - \xi_1, \xi_3 - \xi_1) + m_3^A(\xi_2,\xi_3) + m_3^A(\xi_2, \xi_4) + m_3^A(\xi_1,\xi_2)] \\
& + (m_1^A)^2 [ m_2^A(\xi_1) + m_2^A(\xi_2) + m_2^A(\xi_3) + m_2^A(\xi_3-\xi_1) + m_2^A(\xi_3 - \xi_2) + m_2^A(\xi_2-\xi_1)] \\
& + 6(m_1^A)^4 \\
\vdots \\
\end{split}
\end{align}
Notice that when $A(x)$ is a zero mean process $(m_1^A=0)$, the first three moments and cumulants are equivalent, but differ for orders ($n\geq4$). As for random vectors, non-zero higher-order cumulants indicate non-Gaussianity. In particular, odd-order cumulants give rise to asymmetric non-linearities whereas even-order cumulants give rise to symmetric non-linearities. Further details on the moment and cumulant properties of fields can be found in \cite{Brillinger1965,Shields2017}.


\subsection{Spectral Properties of Random Fields}
\label{S:3}

As discussed in \cite{Brillinger1965,Nikias1987}, it is common and advantageous to work with random fields in the Fourier space. For our purposes, the Fourier domain provides a convenient setting for a nonlinear expansion of random fields that can be derived directly from the third-order spectra. Here we review the spectral quantities necessary for the third-order expansion proposed herein.


\subsubsection{Polyspectra}
\label{S:3.1}

The $n^{th}$-order polyspectrum of a random field $A(x)$ is defined as the Fourier transform of its $n^{th}$-order cumulant \cite{Brillinger1965}
\begin{equation}
\begin{split}
C^{A}_{n}(\kappa_{1}, \kappa_{2}, \dots, \kappa_{n-1}) = &\frac{1}{(2 \pi)^{n-1}} \int_{-\infty}^{\infty} \dots \int_{-\infty}^{\infty} c^{A}_{n}(\xi_{1}, \xi_{2}, \dots , \xi_{n-1})\\
&e^{\iota\kappa_{1}\xi_{1}+\kappa_{2}\xi_{2}+ \dots + \kappa_{n-1}\xi_{n-1}}d\xi_{1} d\xi_{2} \dots d\xi_{n-1} 
\end{split}
\label{eqn:polyspectrum}
\end{equation}
The $2^{nd}$-order polyspectrum, also called the power spectrum, is by far the most widely studied and understood spectral quantity for random fields. The power spectrum, as its name suggests, expresses the power associated with each frequency component of a random field and is defined as
\begin{equation}
S^{A}(\kappa) = C^A_2(\kappa) = \frac{1}{2\pi}\int_{-\infty}^{\infty}c^{A}_{2}(\xi)e^{-\iota(\kappa\xi)}d\xi
\label{eqn:power_spectrum}
\end{equation}
The $3^{rd}$-order polyspectrum is referred to as the bispectrum. The bispectrum describes the third-order nonlinear wave interactions (i.e.\ interactions between wave pairs) in a $3^{rd}$-order stationary random field. It is defined in terms of the $3^{rd}$-order cumulant as
\begin{equation}
    B^A(\kappa_{1}, \kappa_{2}) = C^A_3(\kappa_1,\kappa_2) = \frac{1}{(2 \pi)^{2}}\int_{-\infty}^{\infty}\int_{-\infty}^{\infty}c^A_{3}(\xi_{1}, \xi_{2})e^{-\iota(\kappa_{1}\xi_{1} + \kappa_{2}\xi_{2})}d\xi_{1}d\xi_{2}
    \label{eqn:bispectrum}
\end{equation}

The power spectrum is a real quantity while a bispectrum can have both real and imaginary parts. The real part of the bispectrum corresponds to the Fourier transform of the symmetric part of the third-order cumulant, whereas the imaginary part corresponds to the Fourier transform of the antisymmetric part. As discussed by Lii et al.\ \cite{lii1976bispectral} and Elgar and Guza \cite{elgar1985observations}, the real component relates to the skewness of the field, while the imaginary component relates to the skewness of the derivative of the field. Meanwhile, the amplitude of the bispectrum represents the degree of quadratic phase coupling between the wave-numbers $\kappa_{1}$ and $\kappa_{2}$. A more detailed discussion can be found in \cite{Shields2017} and \cite{kim_thesis_2018}.

\subsubsection{Polycoherence}

For practical purposes, it is useful to induce a normalization of the polyspectrum, which introduces the notion of a polycoherence. Although several normalizations have been proposed \cite{Kim1979,mccomas1980bispectra,hinich2005normalizing}, the $n^{th}$-order squared polycoherence is a standard measure of higher-order polyspectra, and is defined here for stationary random fields as 
\begin{equation}
    |\rho_{A}^{(n)}(\boldsymbol{\kappa})|^2 = \frac{\left|\mathbb{E}\left[\prod_{k=1}^{n-1} dZ(\kappa_{k})dZ^*(\sum_{m=1}^{n-1}\kappa_{m})\right]\right|^2}{\mathbb{E}\left[\prod_{k=1}^{n-1}\left| dZ(\kappa_{k})\right|^2\right]\mathbb{E}\left[\left|dZ(\sum_{m=1}^{n-1}\kappa_{m})\right|^2\right]}
    \label{eqn:polycoherence}
\end{equation}
where $dZ(\kappa)$ are the Fourier coefficients of the generalized random field and $*$ denotes the complex conjugate. Of particular interest here is the third-order polycoherence, or bicoherence which can be derived from Eq.\ \eqref{eqn:polycoherence} and is given by \cite{Kim1979}:
\begin{align}
\label{eq:5_1}
\begin{split}
b_A^2(\kappa_1,\kappa_2) &= \frac{|B^A(\kappa_1, \kappa_2)|^2}{E[|dZ(\kappa_1)dZ(\kappa_2)|^2] S^A(\kappa_1 + \kappa_2)}\\
\end{split}
\end{align}
where $dZ(\kappa)$ are the Fourier coefficients of $A(x)$, $B^A(\kappa_1,\kappa_2)$ is the bispectrum given by Eq.\ \eqref{eqn:bispectrum}, and $S^A(\kappa)$ is the power spectrum from Eq.\ \eqref{eqn:power_spectrum}. By Schwartz' inequality, this definition of the bicoherence is bounded on $[0,1]$ which provides a convenient interpretation of the fraction of energy associated with phase coupling. Further interpretation of the bicoherence can be found in \cite{hinich2005normalizing,mccomas1980bispectra,Shields2017}.

Interestingly, the polycoherence also plays a crucial role in discriminating between non-linearity and non-stationarity in random fields\cite{Hanssen2003}. Here, however, we consider only third-order stationary processes.






\section{Spectral Representation Theorem}
\label{S:4}

Cramer's spectral representation \cite{cramer1967} states that any zero-mean, wide-sense stationary random field $A(x)$ can be expressed in terms of a spectral process $z(\kappa)$ through the following Fourier-Stiltjes integral
\begin{equation}
\label{eq:fourier_stiltjes_integral}
    A(x) = \int_{-\infty}^{\infty}e^{\iota\kappa x}dz(\kappa)
\end{equation}
where the spectral process $z(\kappa)$ has orthogonal increments that satisfy the following properties
\begin{equation}
\begin{aligned}
    &\mathbb{E}[dz(\kappa)] = 0\\
    &\mathbb{E}[z(\kappa)] = 0\\
    &\mathbb{E}[|z(\kappa)|^{2}] = F(\kappa)\\
    &\mathbb{E}[|dz(\kappa)|^{2}] = dF(\kappa)\\
\end{aligned}
\end{equation}
in which $F(\kappa)$ is the spectral distribution function of $z(\kappa)$ and $dF(\kappa)$ is the spectral density function. The power spectrum $S(\kappa)$ of the random field can be expressed in terms of the spectral density $dF(\kappa)$ by $dF(\kappa)= S(\kappa)d\kappa$.

\subsection{$k^{th}$-order Spectral Representation Theorem}

In general, for a zero-mean, $k^{th}$-order stationary random field $A(x)$, a spectral process $z(\kappa)$ can be assigned which satisfies Eq.\ \eqref{eq:fourier_stiltjes_integral}, but possesses additional $k^{th}$-order orthogonality properties
\begin{equation}
\begin{aligned}
    &\mathbb{E}[dz(\kappa)] = 0\\
    &\mathbb{E}[z(\kappa)] = 0\\
    &\mathbb{E}[|z(\kappa)|^{2}] = F(\kappa)\\
    &\mathbb{E}[|dz(\kappa)|^{2}] = dF(\kappa)\\
    &\mathbb{E}[z(\kappa_{1})z(\kappa_{2})z^{*}(\kappa_{3})] = \delta (\kappa_1 + \kappa_2 - \kappa_3)G(\kappa_{1}, \kappa_{2})\\
    &\mathbb{E}[dz(\kappa_{1})dz(\kappa_{2})dz^{*}(\kappa_{3})] = \delta (\kappa_1 + \kappa_2 - \kappa_3)dG(\kappa_{1}, \kappa_{2})\\
    &\vdots\\
    &\mathbb{E}[z(\kappa_{1})z(\kappa_{2}) \hdots z^{*}(\kappa_{k})] = \delta (\kappa_1 + \kappa_2 + \kappa_3 \hdots - \kappa_k) F_{k}(\kappa_1, \kappa_2, \kappa_3 \hdots \kappa_{k-1})\\
    &\mathbb{E}[dz(\kappa_{1})dz(\kappa_{2}) \hdots dz^{*}(\kappa_{k})] = \delta (\kappa_1 + \kappa_2 + \kappa_3 \hdots - \kappa_k) dF_{k}(\kappa_1, \kappa_2, \kappa_3 \hdots \kappa_{k-1})
    \label{eqn:orthogonality}
\end{aligned}
\end{equation}
where $F(\kappa)$ and $dF(\kappa)$ follow from above. $G(\kappa_{1}, \kappa_{2})$ is defined as the bispectral distribution function of the spectral process $z(\kappa)$ and $dG(\kappa_{1}, \kappa_{2})$ is defined as the bispectral density function. The bispectrum in Eq.\ \eqref{eqn:bispectrum} relates with the bispectral density $dG(\kappa_{1}, \kappa_{2})$ through $dG(\kappa_{1}, \kappa_{2}) = B(\kappa_{1}, \kappa_{2})d\kappa_{1}d\kappa_{2}$. Similarly, $F_{k}(\kappa_{1}, \kappa_{2}, \hdots, \kappa_{k-1})$ and $dF_{k}(\kappa_{1}, \kappa_{2}, \hdots, \kappa_{k-1})$ are $k^{th}$-order spectral distribution and density functions respectively. Generalizing, the $k^{th}$-order spectral density function relates to the $k^{th}$-order polyspectrum in Eq.\ \eqref{eqn:polyspectrum} through $dF_{k}(\kappa_{1}, \kappa_{2}, \hdots, \kappa_{k-1})=C_{k}(\kappa_{1}, \kappa_{2}, \hdots, \kappa_{k-1})d\kappa_{1}d\kappa_{2} \hdots d\kappa_{k-1}$.

\subsection{Bispectral Representation Theorem}

Following the $k^{th}$-order Spectral Representation Theorem, we are specifically interested in third-order stationary random fields, for which the orthogonality conditions in Eq.\ \eqref{eqn:orthogonality} hold up to order three. For such random fields, the process is fully stationary in its first, second, and third order properties ($3^{rd}$-order stationary) and can be expressed using the spectral representation in Eq.\ \eqref{eq:fourier_stiltjes_integral} -- referred to herein as the bispectral representation due to the third-order orthogonality and its expression in terms of stationary bispectrum.




\subsection{Real valued random fields}

The spectral representation theorems discussed in the previous sections provide general expressions for complex random fields. Here, we are specifically interested in real-valued random fields, for which the Cramer spectral representation can be written as
\begin{equation}
    A(x) = \int_{-\infty}^{\infty}\cos(\kappa x)du(\kappa) + sin(\kappa x)dv(\kappa)
    \label{eqn:Cramer_real}
\end{equation}
The components $du(\kappa)$ and $dv(\kappa)$ are the real and imaginary components of the orthogonal increments of $dz(\kappa)$ respectively. Both $du(\kappa)$ and $dv(\kappa)$ possess orthogonal properties similar to $dz(\kappa)$. A detailed description of the orthogonality conditions of these components can be found in \cite{Shields2017}.
\bigbreak

\section{Spectral Representation Methods}
\label{S:5}
Although the general form of the spectral representation was developed by Cramer, Rice \cite{Rice1945} was the first to exploit the spectral representation for the purposes of simulation, using its discretized form to simulate one-dimensional, uni-variate Gaussian random processes. Later formalized for second-order multi-dimensional, multi-variate, and non-stationary stochastic processes by Shinozuka \cite{Shinozuka1972_1,Shinozuka1972_2}, the method became known as the spectral representation method (SRM). Properties of stochastic processes simulated by the SRM were elaborated in a series of seminal papers on the method by Shinozuka and Deodatis \cite{Shinozuka1991,Deodatis1996,shinozuka1996simulation}. Recently, Shields and Kim \cite{Shields2017} extended the SRM for simulation of $3^{rd}$-order stationary stochastic processes. In this section, the $2^{nd}$- and $3^{rd}$-order formulations of the SRM are reviewed for one-dimensional, uni-variate random fields. 

\subsection{2nd-order Spectral Representation Method}
\label{S:5.1}


The SRM is used to simulate random fields by discretizing the Cramer spectral representation and simulating the orthogonal increments $du(\kappa)$ and $dv(\kappa)$ in Eq.\ \eqref{eqn:Cramer_real}. Two forms of the orthogonal increments have been proposed for $2^{nd}$-order random fields. One form suggests randomness in the amplitudes of independent increments
\begin{equation}
\begin{aligned}
    du(\kappa_{k}) &= X_{k}\\
    dv(\kappa_{k}) &= Y_{k}\\
\end{aligned}
\end{equation}
where $X_{k}$ and $Y_{k}$ are independent Gaussian Random variables with zero mean and variance equal to $S(\kappa_{k})\Delta \kappa_{k}$. In the other form, randomness is introduced through the phases of orthogonal harmonic functions 
\begin{equation}
\begin{aligned}
    du(\kappa_{k}) &= \sqrt{2}A_{k}\cos(\phi_{k})\\
    dv(\kappa_{k}) &= \sqrt{2}A_{k}\sin(\phi_{k})\\
\end{aligned}
\end{equation}
where $A_{k} = \sqrt{S(\kappa_k)\Delta\kappa}$ and $\phi_k$ are independent uniformly distributed random phase angles in the range $[0, 2\pi]$. Utilizing the second orthogonal increments gives the following form for the second-order SRM to simulate 1-dimensional, uni-variate random fields:
\begin{equation}
    A(x) = \sqrt{2} \sum_{k=0}^{\infty} \sqrt{2S(\kappa_{k})\Delta \kappa_{k}} \cos(\kappa_{k}x - \phi_{k})
\end{equation}
Simulation is then conducted by truncating the summation at an acceptable level, say $m$ terms.


\subsection{$3^{rd}$-order Spectral Representation Method}

The $2^{nd}$-order Spectral Representation Method has been extensively studied and applied over the past several decades and it is a powerful method for simulation of stationary and non-stationary Gaussian random fields. However, its extension to non-Gaussian, or higher-order random fields is not trivial. Thus, attempts to utilize it for such purposes have concentrated on coupling it with other theories (most notably Grigoriu's translation process theory \cite{grigoriu1995applied}). For nearly five decades, it went unrecognized that the spectral representation theory extends beyond $2^{nd}$-order and provides the mathematical framework for higher-order extension of the SRM.


In 2017, Shields and Kim derived this extension, providing a third-order SRM for the simulation of one-dimensional, uni-variate stochastic processes with asymmetric non-linearities \cite{Shields2017}. The authors established that in the presence of third-order spectral information, power associated with a particular frequency $S(\kappa_{k})$ can be decomposed into two components, a pure component $S_{p}(\kappa_{k})$ and an interactive component which arises from coupling of lower frequencies. The authors decompose the orthogonal increments $du(\kappa)$ and $dv(\kappa)$ into the pure components ($du_{p}$, $dv_{p}$) and interactive components ($du_{i}$, $dv_{i}$) as
\begin{equation}
\begin{split}
    du(\kappa_{k}) &= du_{p}(\kappa_{k}) + du_{i}(\kappa_{k})\\
    dv(\kappa_{k}) &= dv_{p}(\kappa_{k}) + dv_{i}(\kappa_{k})
\end{split}
\end{equation}
where
\begin{equation}
\begin{split}
    du_{p}(\kappa_{k}) &= 2\sqrt{S_P(\kappa_{k}) \Delta \kappa_{k} }\cos(\phi_{k})\\
    dv_{p}(\kappa_{k}) &= 2\sqrt{S_P(\kappa_{k}) \Delta \kappa_{k} }\sin(\phi_{k})\\
    du_{i}(\kappa_{k}) &= 2\sqrt{S(\kappa_{k}) \Delta \kappa_{k}} \sum_{i+j=k}^{i\geq j \geq 0} b_{p}(\kappa_{i}, \kappa_{j})\cos(\phi_{i} + \phi_{j} + \beta(\kappa_{i}, \kappa_{j}))\\
    dv_{i}(\kappa_{k}) &= 2\sqrt{S(\kappa_{k}) \Delta \kappa_{k}} \sum_{i+j=k}^{i\geq j \geq 0} b_{p}(\kappa_{i}, \kappa_{j})\sin(\phi_{i} + \phi_{j} + \beta(\kappa_{i}, \kappa_{j}))\\
\end{split}
\end{equation}
and $b_p(\kappa_i,\kappa_j)$ is the partial bicoherence defined as:
\begin{equation}
    b_{p}^{2}(\kappa_{i}, \kappa_{j}) = \frac{\mid B(\kappa_{i}, \kappa_{j}) \mid ^{2}}{S_{P}(\kappa_{i})S_{P}(\kappa_{j})S({\kappa_{i} + \kappa_{j}})}
\end{equation}
$S_P(\kappa)$ is the pure power spectrum (i.e.\ the component of the power spectrum remaining after wave interactions are removed) given by:
\begin{equation}
    S_{P}(\kappa_{k}) = S(\kappa_{k})\left[1 - \sum_{i+j=k}^{i \geq j \geq 0} b_{p}^{2}(\kappa_{i}, \kappa_{j})\right]
\end{equation}
and $\beta(\kappa_i,\kappa_j)$ is the biphase given by:
\begin{equation}
    \beta(\kappa_i,\kappa_j) = \arctan\dfrac{\Im [B(\kappa_i,\kappa_j)]}{\Re[B(\kappa_i,\kappa_j)]}
\end{equation}

Incorporating the proposed orthogonal increments into the discretized spectral representation yields the $3^{rd}$-order form of the SRM
\begin{equation}
\begin{aligned}
    A(x) = &  \sqrt{2}\sum_{k=0}^{\infty}\sqrt{2S_{P}(\kappa_{k}) \Delta \kappa_{k}} \cos (\kappa_{k}x - \phi_{k}) \\
    & + \sqrt{2}\sum_{k=0}^{\infty}\sum_{i+j=k}^{i \geq j \geq 0}\sqrt{2S(\kappa_{i} + \kappa_{j})\Delta (\kappa_{i} + \kappa_{j})} |b_{p}(\kappa_{i}, \kappa_{j})| \\
    & \cos \big((\kappa_{i} + \kappa_{j})x - (\phi_{i} + \phi_{j} + \beta (\kappa_{i},\kappa_{j}))\big)
\end{aligned}
\label{eqn:BSRM}
\end{equation}
where the first term corresponds to the classical $2^{nd}$-order SRM on the pure power spectrum and the second term models $3^{rd}$-order wave interactions. It has been shown in \cite{Shields2017} that simulations generated using Eq.\ \eqref{eqn:BSRM}, again using a suitable truncation of $m$ terms in the summation, match both the power spectrum and the bispectrum of the random field.

\section{Multi-dimensional $3^{rd}$-order Spectral Representation Method}

The form of the $3^{rd}$-order SRM given in Eq.\ \eqref{eqn:BSRM} can be used for the simulation of one-dimensional, uni-variate (1D-1V) random fields. In this section, we derive the expression for the simulation of general $d$-dimensional ($d$D-1V) third-order random fields. We first derive the simulation formula for two-dimensional random fields as this case is the most practical to show and is of particular relevance for many applications. We then extend it for three-dimensional and general $d$-dimensional random fields.

\subsection{Simulation of 2-dimensional random fields}

Let $A(x_{1}, x_{2})$ be a two-dimensional uni-variate random field with zero mean, power spectrum $S(\kappa_{1}, \kappa_{2})$, $2^{nd}$-order autocorrelation function $R_{2}(\xi_{1}, \xi_{2})$, Bispectrum $B(\kappa_{11}, \kappa_{21}, \kappa_{12}, \kappa_{22})$, and $3^{rd}$-order autocorrelation function $R_{3}(\xi_{11}, \xi_{21}, \xi_{12}, \xi_{22})$ satisfying:
\begin{equation}
    \mathbb{E}[A(x_{1}, x_{2})] = 0
\end{equation}
\begin{equation}
    \mathbb{E}[A(x_{1} + \xi_{1}, x_{2} + \xi_{2})A(x_{1}, x_{2})] = R_{2}(\xi_{1}, \xi_{2})
\end{equation}
\begin{equation}
    \mathbb{E}[A(x_{1} + \xi_{11}, x_{2} + \xi_{12})A(x_{1} + \xi_{21}, x_{2} + \xi_{22})A(x_{1}, x_{2})] = R_{3}(\xi_{11}, \xi_{21}, \xi_{12}, \xi_{22})
\end{equation}
\begin{equation}
    S(\kappa_{1}, \kappa_{2}) = \frac{1}{(2\pi)^{2}}\int_{-\infty}^{\infty}\int_{-\infty}^{\infty}R_{2}(\xi_{1}, \xi_{2})e^{-\iota(\kappa_{1}\xi_{1} + \kappa_{2}\xi_{2})}d\xi_{1}d\xi_{2}
\label{eqn:WK1}
\end{equation}
\begin{equation}
\begin{aligned}
    B(\kappa_{11}, \kappa_{21}, \kappa_{12}, \kappa_{22}) &= \frac{1}{(2\pi)^{4}}\int_{-\infty}^{\infty}\int_{-\infty}^{\infty}\int_{-\infty}^{\infty}\int_{-\infty}^{\infty}R_{3}(\xi_{11}, \xi_{21}, \xi_{12}, \xi_{22})\\
    &e^{-\iota(\kappa_{11}\xi_{11} + \kappa_{21}\xi_{21} + \kappa_{12}\xi_{12} + \kappa_{22}\xi_{22})}d\xi_{11}d\xi_{21}d\xi_{12}d\xi_{22}
\end{aligned}
\label{eqn:WK2}
\end{equation}
\begin{equation}
    R_{2}(\xi_{1}, \xi_{2}) = \int_{-\infty}^{\infty}\int_{-\infty}^{\infty}S(\kappa_{1}, \kappa_{2})e^{\iota(\kappa_{1}\xi_{1} + \kappa_{2}\xi_{2})}d\kappa_{1}d\kappa_{2}
\label{eqn:WK3}
\end{equation}
\begin{equation}
\begin{aligned}
     R_{3}(\xi_{11}, \xi_{21}, \xi_{12}, \xi_{22}) &= \int_{-\infty}^{\infty}\int_{-\infty}^{\infty}\int_{-\infty}^{\infty}\int_{-\infty}^{\infty}B(\kappa_{11}, \kappa_{21}, \kappa_{12}, \kappa_{22})\\
     &e^{\iota(\kappa_{11}\xi_{11} + \kappa_{21}\xi_{21} + \kappa_{12}\xi_{12} + \kappa_{22}\xi_{22})}d\kappa_{11}d\kappa_{21}d\kappa_{12}d\kappa_{22}
\end{aligned}
\label{eqn:WK4}
\end{equation}
where Eq.\ \eqref{eqn:WK1}, Eq.\ \eqref{eqn:WK3} and Eq.\ \eqref{eqn:WK2}, Eq.\ \eqref{eqn:WK4} constitute the $2^{nd}$ and $3^{rd}$ order Weiner-Khintchine transform pairs respectively.


Since we are interested in the simulation of real-valued random fields, the power spectrum is symmetric about the origin, i.e.
\begin{equation}
\label{eqn:power_spectrum_symmetry}
\begin{aligned}
	&S(\boldsymbol{\kappa}) = S(-\boldsymbol{\kappa}),
\end{aligned}
\end{equation}
which equates to the following two equations for 2-dimensional random fields
\begin{equation}
\begin{aligned}
    & S(\kappa_{1}, \kappa_{2}) = S(-\kappa_{1}, -\kappa_{2})\\
    & S(\kappa_{1}, -\kappa_{2}) = S(-\kappa_{1}, \kappa_{2})\\
\end{aligned}
\end{equation}
Furthermore, the following symmetries exist in the bispectrum \citep{Shields2017}
\begin{equation}
\label{eqn:bispectrum_symmetry_1}
\begin{aligned}
		& B(\boldsymbol{\kappa_{1}}, \boldsymbol{\kappa_{2}}) = B(\boldsymbol{\kappa_{2}}, \boldsymbol{\kappa_{1}})\\
\end{aligned}
\end{equation}
\begin{equation}
\label{eqn:bispectrum_symmetry_2}
\begin{aligned}
	& B(\boldsymbol{\kappa_{1}}, \boldsymbol{\kappa_{2}}) = B(-\boldsymbol{\kappa_{1}}, -\boldsymbol{\kappa_{2}})\\
\end{aligned}
\end{equation}
\begin{equation}
\label{eqn:bispectrum_symmetry_3}
\begin{aligned}
	& B(\boldsymbol{\kappa_{1}}, \boldsymbol{\kappa_{2}}) = B(-\boldsymbol{\kappa_{1}}-\boldsymbol{\kappa_{2}}, \boldsymbol{\kappa_{2}})\\
\end{aligned}
\end{equation}
Eqs.\ \eqref{eqn:bispectrum_symmetry_2} and  \eqref{eqn:bispectrum_symmetry_3} describe two different axes of symmetry along the origin. These symmetries for the 2-dimensional bispectrum are as follows
\begin{equation}
\begin{aligned}
    & B(\kappa_{11}, \kappa_{21}, \kappa_{12}, \kappa_{22}) =  B(\kappa_{11}, -\kappa_{21}, \kappa_{12}, \kappa_{22}) = B(-\kappa_{11}, -\kappa_{21}, -\kappa_{12}, -\kappa_{22}) = B(-\kappa_{11}, \kappa_{21}, -\kappa_{12}, -\kappa_{22})\\
    & B(\kappa_{11}, \kappa_{21}, \kappa_{12}, -\kappa_{22}) = B(\kappa_{11}, -\kappa_{21}, \kappa_{12}, -\kappa_{22}) = B(-\kappa_{11}, -\kappa_{21}, -\kappa_{12}, \kappa_{22}) = B(-\kappa_{11}, \kappa_{21}, -\kappa_{12}, \kappa_{22})\\
    & B(\kappa_{11}, \kappa_{21}, -\kappa_{12}, \kappa_{22}) = B(\kappa_{11}, -\kappa_{21}, -\kappa_{12}, \kappa_{22}) = B(-\kappa_{11}, -\kappa_{21}, \kappa_{12}, -\kappa_{22}) = B(-\kappa_{11}, \kappa_{21}, \kappa_{12}, -\kappa_{22})\\
    & B(\kappa_{11}, \kappa_{21}, -\kappa_{12}, -\kappa_{22}) = B(\kappa_{11}, -\kappa_{21}, -\kappa_{12}, -\kappa_{22}) = B(-\kappa_{11}, -\kappa_{21}, \kappa_{12}, \kappa_{22}) = B(-\kappa_{11}, \kappa_{21}, \kappa_{12}, \kappa_{22})\\
\end{aligned}
\end{equation}

Exploiting these symmetries allows us to replace the power spectrum, $S(\kappa_{1}, \kappa_{2})$ defined on the range $(-\infty \leq \kappa_{1} \leq \infty, -\infty \leq \kappa_{2} \leq \infty)$ by $2S(\kappa_{1}, \kappa_{2})$ defined on the range $(0 \leq \kappa_{1} \leq \infty, -\infty \leq \kappa_{2} \leq \infty)$ and replace the bispectrum $B(\kappa_{11}, \kappa_{12}, \kappa_{21}, \kappa_{22})$ defined on the range $(-\infty \leq \kappa_{11} \leq \infty, -\infty \leq \kappa_{21} \leq \infty, -\infty \leq \kappa_{12} \leq \infty, -\infty \leq \kappa_{22} \leq \infty)$ by $4B(\kappa_{11}, \kappa_{12}, \kappa_{21}, \kappa_{22})$ defined on the range $(0 \leq \kappa_{11} \leq \infty, 0 \leq \kappa_{21} \leq \infty, -\infty \leq \kappa_{12} \leq \infty, -\infty \leq \kappa_{22} \leq \infty)$

With these symmetries in place, along with the orthogonality conditions presented in Eq.\ \eqref{eqn:orthogonality}, any real valued 2-dimensional random field $A(x_{1}, x_{2})$ can be expressed in the form
\begin{equation}
\label{eqn:simulation_orthogonal_increments_form}
    A(x_{1}, x_{2}) = \int_{-\infty}^{\infty}\int_{0}^{\infty}[\cos(\kappa_{1}x_{1} + \kappa_{2}x_{2})du(\kappa_{1}, \kappa_{2}) + \sin(\kappa_{1}x_{1} + \kappa_{2}x_{2})dv(\kappa_{1}, \kappa_{2})]
\end{equation}
where processes $u(\kappa_{1}, \kappa_{2})$ and $v(\kappa_{1}, \kappa_{2})$ are defined on the domain $0 < \kappa_{1} < \infty, -\infty < \kappa_{2} < \infty $ and obey the following the orthogonality conditions:
\begin{equation}
\label{eqn:1_order_orthogonality}
    \mathbb{E}[u(\kappa_{1}, \kappa_{2})] = \mathbb{E}[v(\kappa_{1}, \kappa_{2})] = 0
\end{equation}

\begin{equation}
\label{eqn:1_order_orthogonal_increments}
    \mathbb{E}[du(\kappa_{1}, \kappa_{2})] = \mathbb{E}[dv(\kappa_{1}, \kappa_{2})] = 0
\end{equation}

\begin{equation}
\begin{aligned}
\label{eqn:2_order_orthogonality}
    \mathbb{E}[u^{2}(\kappa_{1}, \kappa_{2})] = \mathbb{E}[v^{2}(\kappa_{1}, \kappa_{2})] & = F_{1}(\kappa_{1}, \kappa_{2})\\
    \mathbb{E}[u(\kappa_{11}, \kappa_{21})u(\kappa_{12}, \kappa_{22})u(\kappa_{11} + \kappa_{12}, \kappa_{21} + \kappa_{22})] &=\\
    \mathbb{E}[v(\kappa_{11}, \kappa_{21})v(\kappa_{12}, \kappa_{22})v(\kappa_{11} + \kappa_{12}, \kappa_{21} + \kappa_{22})] &= G_{1}(\kappa_{11}, \kappa_{21}, \kappa_{12}, \kappa_{22})\\
    \mathbb{E}[u(\kappa_{1}, \kappa_{2})v(\kappa'_{1}, \kappa'_{2})] & = 0\\
    \mathbb{E}[u(\kappa_{1}, \kappa_{2})v(\kappa'_{1}, \kappa'_{2})v(\kappa''_{1}, \kappa''_{2})] & = 0\\
    \mathbb{E}[u(\kappa_{1}, \kappa_{2})u(\kappa'_{1}, \kappa'_{2})v(\kappa''_{1}, \kappa''_{2})] & = 0\\
\end{aligned}
\end{equation}

\begin{equation}
\begin{aligned}
\label{eqn:2_order_orthogonal_increments}
    \mathbb{E}[du^{2}(\kappa_{1}, \kappa_{2})] = \mathbb{E}[dv^{2}(\kappa_{1}, \kappa_{2})] &= S_{1}(\kappa_{1}, \kappa_{2})d\kappa_{1}d\kappa_{2}\\
    \mathbb{E}[du(\kappa_{1}, \kappa_{2})du(\kappa'_{1}, \kappa'_{2})] &= 0 \: if \kappa_{1} \neq \kappa'_{1} \: or \kappa_{2} \neq \kappa'_{2}\\
    \mathbb{E}[dv(\kappa_{1}, \kappa_{2})dv(\kappa'_{1}, \kappa'_{2})] & = 0 \: if \kappa_{1} \neq \kappa'_{1} \: or \kappa_{2} \neq \kappa'_{2}\\
    \mathbb{E}[du(\kappa_{1}, \kappa_{2})dv(\kappa'_{1}, \kappa'_{2})] &= 0\\
\end{aligned}
\end{equation}


\begin{equation}
\label{eqn:3_order_orthogonal_increments}
\begin{aligned}
	& \mathbb{E}[du(\kappa_{1}, \kappa_{2})du(\kappa'_{1}, \kappa'_{2})du(\kappa''_{1}, \kappa''_{2})] = 2\Re B(\kappa'_{1}, \kappa''_{1}, \kappa'_{2}, \kappa''_{2})\\
	& \mathbb{E}[du(\kappa_{1}, \kappa_{2})du(\kappa'_{1}, \kappa'_{2})dv(\kappa''_{1}, \kappa''_{2})] = -2\Im B(\kappa'_{1}, \kappa''_{1}, \kappa'_{2}, \kappa''_{2})\\
	& \mathbb{E}[du(\kappa_{1}, \kappa_{2})dv(\kappa'_{1}, \kappa'_{2})du(\kappa''_{1}, \kappa''_{2})] = -2\Im B(\kappa'_{1}, \kappa''_{1}, \kappa'_{2}, \kappa''_{2})\\
	& \mathbb{E}[du(\kappa_{1}, \kappa_{2})dv(\kappa'_{1}, \kappa'_{2})dv(\kappa''_{1}, \kappa''_{2})] = -2\Re B(\kappa'_{1}, \kappa''_{1}, \kappa'_{2}, \kappa''_{2})\\
	& \mathbb{E}[dv(\kappa_{1}, \kappa_{2})du(\kappa'_{1}, \kappa'_{2})du(\kappa''_{1}, \kappa''_{2})] = 2\Im B(\kappa'_{1}, \kappa''_{1}, \kappa'_{2}, \kappa''_{2})\\
	& \mathbb{E}[dv(\kappa_{1}, \kappa_{2})du(\kappa'_{1}, \kappa'_{2})dv(\kappa''_{1}, \kappa''_{2})] = 2\Re B(\kappa'_{1}, \kappa''_{1}, \kappa'_{2}, \kappa''_{2})\\
	& \mathbb{E}[dv(\kappa_{1}, \kappa_{2})dv(\kappa'_{1}, \kappa'_{2})du(\kappa''_{1}, \kappa''_{2})] = 2\Re B(\kappa'_{1}, \kappa''_{1}, \kappa'_{2}, \kappa''_{2})\\
	& \mathbb{E}[dv(\kappa_{1}, \kappa_{2})dv(\kappa'_{1}, \kappa'_{2})dv(\kappa''_{1}, \kappa''_{2})] = -2\Im B(\kappa'_{1}, \kappa''_{1}, \kappa'_{2}, \kappa''_{2})\\
	& \text{if} \ \kappa_{1} = \kappa'_{1} + \kappa''_{1} \\ & \text{otherwise} \ 0
\end{aligned}
\end{equation}
where $\Re $ and $\Im $ denote the real and imaginary components respectively. It is shown in \ref{A:1} that Eq.\ \eqref{eqn:simulation_orthogonal_increments_form} does, indeed represent a stochastic field with zero mean and $2^{nd}$-order and $3^{rd}$-order Autocorrelation Functions $R_{2}(\xi_{1}, \xi_{2})$ and $R_{3}(\xi_{11}, \xi_{21}, \xi_{12}, \xi_{22})$ respectively.



Discretizing Eq.\ \eqref{eqn:simulation_orthogonal_increments_form}, gives
\begin{equation}
    A(x_{1}, x_{2}) = \sum_{n_{2}=-\infty}^{\infty}\sum_{n_{1}=0}^{\infty}[\cos(\kappa_{1n_{1}}x_{1} + \kappa_{2n_{2}}x_{2})du(\kappa_{1n_{1}}, \kappa_{2n_{2}}) + \sin(\kappa_{1n_{1}}x_{1} + \kappa_{2n_{2}}x_{2})dv(\kappa_{1n_{1}}, \kappa_{2n_{2}})]\\
\end{equation}
where
\begin{equation}
    \kappa_{1n_{1}} = n_{1}\Delta\kappa_{1}
\end{equation}
\begin{equation}
    \kappa_{2n_{2}} = n_{2}\Delta\kappa_{2}
\end{equation}
with sufficiently small finite $\Delta\kappa_{1}$ and $\Delta\kappa_{2}$. If $du(\kappa_{1n_{1}}, \kappa_{2n_{2}})$ and $dv(\kappa_{1n_{1}}, \kappa_{2n_{2}})$ are defined as
\begin{equation}
\begin{aligned}
    du(\kappa_{1n_{1}}, \kappa_{2n_{2}}) &= \sqrt{2}A_{pn_{1}n_{2}}\cos\Phi_{n_{1}n_{2}}\\
    &+ \sum_{i_{1} + j_{1} = n_{1}}^{i_{1} \geq j_{1} \geq 0}\sum_{i_{2} + j_{2} = n_{2}}^{|i_{2}| \geq |j_{2}| \geq 0}\sqrt{2}A_{n_{1}n_{2}}b_{p}(\kappa_{1i_{1}}, \kappa_{1j_{1}}, \kappa_{2i_{2}}, \kappa_{2j_{2}})\\
    & \cos(\Phi_{i_{1}i_{2}} + \Phi_{j_{1}j_{2}} + \beta(\kappa_{1i_{1}}, \kappa_{1j_{1}}, \kappa_{2i_{2}}, \kappa_{2j_{2}}))
\end{aligned}
\end{equation}
\begin{equation}
\begin{aligned}
    dv(\kappa_{1n_{1}}, \kappa_{2n_{2}}) &= -\sqrt{2}A_{pn_{1}n_{2}}\sin\Phi_{n_{1}n_{2}}\\
    &- \sum_{i_{1} + j_{1} = n_{1}}^{i_{1} \geq j_{1} \geq 0}\sum_{i_{2} + j_{2} = n_{2}}^{|n_{2}| \geq |i_{2}| \geq |j_{2}| \geq 0}\sqrt{2}A_{n_{1}n_{2}}b_{p}(\kappa_{1i_{1}}, \kappa_{1j_{1}}, \kappa_{2i_{2}}, \kappa_{2j_{2}})\\
    & \sin(\Phi_{i_{1}i_{2}} + \Phi_{j_{1}j_{2}} + \beta(\kappa_{1i_{1}}, \kappa_{1j_{1}}, \kappa_{2i_{2}}, \kappa_{2j_{2}}))
\end{aligned}
\end{equation}
where
\begin{equation}
    A_{pn_{1}n_{2}} = \sqrt{2S_{p}(\kappa_{1n_{1}}, \kappa_{2n_{2}})\Delta\kappa_{1}\Delta\kappa_{2}}
\end{equation}
\begin{equation}
    A_{n_{1}n_{2}} = \sqrt{2S(\kappa_{1n_{1}}, \kappa_{2n_{2}})\Delta\kappa_{1}\Delta\kappa_{2}}
\end{equation}
\begin{equation}
    S_{p}(\kappa_{1n_{1}}, \kappa_{2n_{2}}) = S(\kappa_{1n_{1}}, \kappa_{2n_{2}})\Big(1 - \sum_{i_{1} + j_{1} = n_{1}}^{i_{1} \geq j_{1} \geq 0}\sum_{i_{2} + j_{2} = n_{2}}^{|i_{2}| \geq |j_{2}| \geq 0}b_{p}^{2}(\kappa_{1i_{1}}, \kappa_{1j_{1}}, \kappa_{2i_{2}}, \kappa_{2j_{2}})\Big)
\end{equation}
\begin{equation}
\begin{aligned}
	&b_{p}^{2}(\kappa_{1i_{1}}, \kappa_{1j_{1}}, \kappa_{2i_{2}}, \kappa_{2j_{2}}) = \frac{|B(\kappa_{1i_{1}}, \kappa_{1j_{1}}, \kappa_{2i_{2}}, \kappa_{2j_{2}}|^{2}\Delta\kappa_{1}\Delta\kappa_{2}}{S_{p}(\kappa_{1i_{1}}, \kappa_{2i_{2}})S_{p}(\kappa_{1j_{1}}, \kappa_{2j_{2}})S(\kappa_{1(i_{1} + j_{1})}, \kappa_{2(i_{2} + j_{2})})}\\
\end{aligned}
\end{equation}
and $\Phi_{n_{1}n_{2}}$ are independent random phase angles uniformly distributed in the range $[0, 2\pi]$, then the resulting 2-dimensional random field is third-order stationary possessing power spectrum $S(\kappa_1,\kappa_2)$ and bispectrum $B(\kappa_{11}, \kappa_{12}, \kappa_{21}, \kappa_{22})$. \ref{A:2} shows that the orthogonality requirements on $du(\kappa_{1n_{1}}, \kappa_{2n_{2}})$ and $dv(\kappa_{1n_{1}}, \kappa_{2n_{2}})$ are satisfied, and therefore that the process is third-order stationary possessing the prescribed spectra.


Using the above proposed increments, the following series representation is obtained
\begin{equation}
\label{eqn:infinite_sum_2d_1}
\begin{aligned}
    A(x_{1}, x_{2}) &= \sum_{n_{2}=-\infty}^{\infty}\sum_{n_{1}=0}^{\infty}[\sqrt{2}A_{pn_{1}n_{2}}\cos(\kappa_{1n_{1}}x_{1} + \kappa_{2n_{2}}x_{2} + \Phi_{n_{1}n_{2}})\\
    & + [\sum_{i_{1} + j_{1} = n_{1}}^{i_{1} \geq j_{1} \geq 0}\sum_{i_{2} + j_{2} = n_{2}}^{|i_{2}| \geq |j_{2}| \geq 0}\sqrt{2}A_{n_{1}n_{2}}b_{p}(\kappa_{1i_{1}}, \kappa_{1j_{1}}, \kappa_{2i_{2}}, \kappa_{2j_{2}})\\
    &\cos(\kappa_{1n_{1}}x_{1} + \kappa_{2n_{2}}x_{2} + \Phi_{i_{1}i_{2}} + \Phi_{j_{1}j_{2}} + \beta(\kappa_{1i_{1}}, \kappa_{1j_{1}}, \kappa_{2i_{2}}, \kappa_{2j_{2}}))]\\
\end{aligned}
\end{equation}
By rearranging the terms, we can express the series over only positive indices as
\begin{equation}
\label{eqn:infinite_sum_2d_2}
\begin{aligned}
    A(x_{1}, x_{2}) &= \sqrt{2}\sum_{n_{2}=0}^{\infty}\sum_{n_{1}=0}^{\infty}\Big[\sqrt{S_{p}(\kappa_{1n_{1}}, \kappa_{2n_{2}})\Delta\kappa_{1}\Delta\kappa_{2}}\cos(\kappa_{1n_{1}}x_{1} + \kappa_{2n_{2}}x_{2} + \Phi_{n_{1}n_{2}}^{(1)})\\
    & + \sqrt{S_{p}(\kappa_{n_{1}}, -\kappa_{n_{2}})\Delta\kappa_{1}\Delta\kappa_{2}}\cos(\kappa_{1n_{1}}x_{1} - \kappa_{2n_{2}}x_{2} + \Phi_{n_{1}n_{2}}^{(2)})\\
    & + \sum_{i_{1} + j_{1} = n_{1}}^{i_{1} \geq j_{1} \geq 0}\sum_{i_{2} + j_{2} = n_{2}}^{i_{2} \geq j_{2} \geq 0}\sqrt{2S(\kappa_{1n_{1}}, \kappa_{2n_{2}})}b_{p}(\kappa_{1i_{1}}, \kappa_{1j_{1}}, \kappa_{2i_{2}}, \kappa_{2j_{2}})\\
    &\cos(\kappa_{1n_{1}}x_{1} + \kappa_{2n_{2}}x_{2} + \Phi_{i_{1}i_{2}}^{(1)} + \Phi_{j_{1}j_{2}}^{(1)} + \beta(\kappa_{1i_{1}}, \kappa_{1j_{1}}, \kappa_{2i_{2}}, \kappa_{2j_{2}}))\\
    & + \sum_{i_{1} + j_{1} = n_{1}}^{i_{1} \geq j_{1} \geq 0}\sum_{i_{2} + j_{2} = n_{2}}^{i_{2} \geq j_{2} \geq 0}\sqrt{2S(\kappa_{1n_{1}}, -\kappa_{2n_{2}})}b_{p}(\kappa_{1i_{1}}, \kappa_{1j_{1}}, -\kappa_{2i_{2}}, -\kappa_{2j_{2}})\\
    &\cos(\kappa_{1n_{1}}x_{1} - \kappa_{2n_{2}}x_{2} + \Phi_{i_{1}i_{2}}^{(2)} + \Phi_{j_{1}j_{2}}^{(2)} + \beta(\kappa_{1i_{1}}, \kappa_{1j_{1}}, -\kappa_{2i_{2}}, -\kappa_{2j_{2}}))\\
    & + \sum_{i_{1} + j_{1} = n_{1}}^{i_{1} \geq j_{1} \geq 0}\sum_{i_{2} + j_{2} = n_{2}}^{i_{2} \geq j_{2} \geq 0}\sqrt{2S(\kappa_{1n_{1}}, -\kappa_{2i_{2}} + \kappa_{2j_{2}})}b_{p}(\kappa_{1i_{1}}, \kappa_{1j_{1}}, -\kappa_{2i_{2}}, \kappa_{2j_{2}})\\
    &\cos(\kappa_{1n_{1}}x_{1} - \kappa_{2i_{2}}x_{2} + \kappa_{2j_{2}}x_{2} + \Phi_{i_{1}i_{2}}^{(2)} + \Phi_{j_{1}j_{2}}^{(1)} + \beta(\kappa_{1i_{1}}, \kappa_{1j_{1}}, -\kappa_{2i_{2}}, \kappa_{2j_{2}}))\\
    & + \sum_{i_{1} + j_{1} = n_{1}}^{i_{1} \geq j_{1} \geq 0}\sum_{i_{2} + j_{2} = n_{2}}^{i_{2} \geq j_{2} \geq 0}\sqrt{2S(\kappa_{1n_{1}}, +\kappa_{2i_{2}} - \kappa_{2j_{2}})}b_{p}(\kappa_{1i_{1}}, \kappa_{1j_{1}}, \kappa_{2i_{2}}, -\kappa_{2j_{2}})\\
    &\cos(\kappa_{1n_{1}}x_{1} + \kappa_{2i_{2}}x_{2} - \kappa_{2j_{2}}x_{2} + \Phi_{i_{1}i_{2}}^{(1)} + \Phi_{j_{1}j_{2}}^{(2)} + \beta(\kappa_{1i_{1}}, \kappa_{1j_{1}}, \kappa_{2i_{2}}, -\kappa_{2j_{2}}))\Big]\\
\end{aligned}
\end{equation}
While Eq.\ \eqref{eqn:infinite_sum_2d_1} provides a compact notation, Eq.\ \eqref{eqn:infinite_sum_2d_2} sums only over positive indices which may be beneficial for practical implementation. Note that, since the formula sums over the positive and negative range of $\kappa_{2}$ simultaneously, we need to use two different sets of random phase angles which are differentiated using superscripts $\Phi^{(1)}$ and $\Phi^{(2)}$. 

\subsubsection{Simulation formula for general 2-dimensional random fields}
\label{general_2d}

While Eqs.\ \eqref{eqn:infinite_sum_2d_1} - \eqref{eqn:infinite_sum_2d_2} provide a theoretical framework for the simulation of 2-dimensional third-order stationary random fields, the infinite series representation of  Eq.\ \eqref{eqn:infinite_sum_2d_1} cannot be implemented in practice. A practical implementation truncates these summations as
\begin{equation}
\label{eqn:finite_sum_2d}
\begin{aligned}
    A(x_{1}, x_{2}) &= \sum_{n_{2}=-N_{2}}^{N_{2}}\sum_{n_{1}=0}^{N_1}[\sqrt{2}A_{pn_{1}n_{2}}\cos(\kappa_{1n_{1}}x_{1} + \kappa_{2n_{2}}x_{2} + \Phi_{n_{1}n_{2}})\\
    & + [\sum_{i_{1} + j_{1} = n_{1}}^{i_{1} \geq j_{1} \geq 0}\sum_{i_{2} + j_{2} = n_{2}}^{|n_{2}| \geq |i_{2}| \geq |j_{2}| \geq 0}\sqrt{2}A_{n_{1}n_{2}}b_{p}(\kappa_{1i_{1}}, \kappa_{1j_{1}}, \kappa_{2i_{2}}, \kappa_{2j_{2}})\\
    &\cos(\kappa_{1n_{1}}x_{1} + \kappa_{2n_{2}}x_{2} + \Phi_{i_{1}i_{2}} + \Phi_{j_{1}j_{2}} + \beta(\kappa_{1i_{1}}, \kappa_{1j_{1}}, \kappa_{2i_{2}}, \kappa_{2j_{2}}))]\\
\end{aligned}
\end{equation}
where
\begin{equation}
\begin{aligned}
    & A_{pn_{1}n_{2}} = \sqrt{2S_{p}(\kappa_{1n_{1}}, \kappa_{2n_{2}})\Delta\kappa_{1}\Delta\kappa_{2}}\\
    & A_{n_{1}n_{2}} = \sqrt{2S(\kappa_{1n_{1}}, \kappa_{2n_{2}})\Delta\kappa_{1}\Delta\kappa_{2}}\\
    & S_{p}(\kappa_{1n_{1}}, \kappa_{2n_{2}}) = S(\kappa_{1n_{1}}, \kappa_{2n_{2}})\Big(1 - \sum_{i_{1} + j_{1} = n_{1}}^{i_{1} \geq j_{1} \geq 0}\sum_{i_{2} + j_{2} = n_{2}}^{|n_{2}| \geq |i_{2}| \geq |j_{2}| \geq 0}b_{p}^{2}(\kappa_{1i_{1}}, \kappa_{1j_{1}}, \kappa_{2i_{2}}, \kappa_{2j_{2}})\Big)\\
	& b_{p}^{2}(\kappa_{1i_{1}}, \kappa_{2j_{1}}, \kappa_{1i_{2}}, \kappa_{2j_{2}}) = \frac{|B(\kappa_{1i_{1}}, \kappa_{2j_{1}}, \kappa_{1i_{2}}, \kappa_{2j_{2}}|^{2}\Delta\kappa_{1}\Delta\kappa_{2}}{S_{p}(\kappa_{1i_{1}}, \kappa_{2i_{2}})S_{p}(\kappa_{1j_{1}}, \kappa_{2j_{2}})S(\kappa_{1(i_{1} + j_{1})}, \kappa_{2(i_{2} + j_{2})})}\\
	& \kappa_{1n_{1}} = n_{1}\Delta\kappa_{1} \ ; \  \kappa_{2n_{2}} = n_{2}\Delta\kappa_{2} \\
	& \Delta\kappa_{1} = \frac{\kappa_{1u}}{N_{1}} \ ; \  \Delta\kappa_{2} = \frac{\kappa_{2u}}{N_{2}} \\
\end{aligned}
\end{equation}
and
\begin{equation}
	S(\kappa_{1}, 0) = S(0, \kappa_{2}) = 0 \ for -\infty \leq \kappa_{1} \leq \infty \ and -\infty \leq \kappa_{2} \leq \infty
\end{equation}
\begin{equation}
\begin{aligned}
	& B(\kappa_{11}, \kappa_{12}, \kappa_{21}, 0) = B(\kappa_{11}, \kappa_{12}, 0, \kappa_{22}) = B(\kappa_{11}, 0, \kappa_{21}, \kappa_{22}) = B(0, \kappa_{12}, \kappa_{21}, \kappa_{22}) = 0 \\
	&  \text{for} -\infty \leq \kappa_{11} \leq \infty \ ; -\infty \leq \kappa_{12} \leq \infty \ \text{and} -\infty \leq \kappa_{21} \leq \infty ; -\infty \leq \kappa_{22} \leq \infty \\
\end{aligned}
\end{equation}
$\kappa_{1u}$ and $\kappa_{2u}$ are the cutoff wave-numbers for the $x_{1}$ and $x_{2}$ axes respectively. These cutoff wave-numbers are chosen to satisfy the condition
\begin{equation}
\begin{aligned}
    & \int_{0}^{\kappa_{1u}}\int_{-\kappa_{2u}}^{\kappa_{2u}}S(\kappa_{1}, \kappa_{2})d\kappa_{1}d\kappa_{2} = (1 - \epsilon)\int_{0}^{\infty}\int_{-\infty}^{\infty}S(\kappa_{1}, \kappa_{2})d\kappa_{1}d\kappa_{2}\\
\end{aligned}
\end{equation}

\begin{equation}
\begin{aligned}
    & \int_{0}^{\kappa_{1u}}\int_{0}^{\kappa_{1u}}\int_{-\kappa_{2u}}^{\kappa_{2u}}\int_{-\kappa_{2u}}^{\kappa_{2u}}B(\kappa_{11}, \kappa_{12}, \kappa_{21}, \kappa_{22})d\kappa_{11}d\kappa_{12}d\kappa_{21}d\kappa_{22}\\
    & = (1 - \epsilon)\int_{0}^{\kappa_{1u}}\int_{0}^{\kappa_{1u}}\int_{-\kappa_{2u}}^{\kappa_{2u}}\int_{-\kappa_{2u}}^{\kappa_{2u}}B(\kappa_{11}, \kappa_{12}, \kappa_{21}, \kappa_{22})d\kappa_{11}d\kappa_{12}d\kappa_{21}d\kappa_{22}\\
\end{aligned}
\end{equation}
where $\epsilon \ll 1$. This effectively means that the power spectrum and the bispectrum above the cutoff wave-numbers are mathematically or physically insignificant. 

It is further straightforward to show that the the simulated random fields are periodic along the $x_{1}$ and $x_{2}$ axes with period
\begin{equation}
\begin{aligned}
    & L_{x_{1}} = \frac{2\pi}{\Delta \kappa_{1}} \\
    & L_{x_{2}} = \frac{2\pi}{\Delta \kappa_{2}} \\
\end{aligned}    
\end{equation}

Additionally, the following conditions are imposed on the spatial increments to prevent aliasing while generating the samples
\begin{equation}
\begin{aligned}
    & \Delta x_{1} \leq \frac{2\pi}{2\kappa_{1u}} \\
    & \Delta x_{2} \leq \frac{2\pi}{2\kappa_{2u}} \\
\end{aligned}    
\end{equation}

\subsubsection{Simulation formula for 2-dimensional quadrant random fields}
\label{quadrant_2d}

Quadrant random fields have additional symmetry beyond the symmetry presented in  Eq.\ \eqref{eqn:power_spectrum_symmetry} and  Eq.\ \eqref{eqn:bispectrum_symmetry_2}-\eqref{eqn:bispectrum_symmetry_3}. Specifically,
\begin{equation}
\begin{aligned}
	&S(\kappa_{1}, \kappa_{2}) = S(I_{1}\kappa_{1}, I_{2}\kappa_{2}) \ \text{for} \ I_{1},  I_{2} = \pm 1\\
	&B(\kappa_{11}, \kappa_{12}, \kappa_{21}, \kappa_{22}) = B(I_{11}\kappa_{11}, I_{12}\kappa_{12}, I_{21}\kappa_{21}, I_{22}\kappa_{22}) \ \text{for} \ I_{11},  I_{12}, I_{21}, I_{22} = \pm 1\\
\end{aligned}
\end{equation}

As a result of these additional symmetries, the simulation formula for 2D-1V third-order quadrant random fields simplifies to
\begin{equation}
\label{eqn:2d_quadrant}
\begin{aligned}
    A(x_{1}, x_{2}) &= \sqrt{2}\sum_{n_{2}=0}^{N_{2}}\sum_{n_{1}=0}^{N_{1}}\sqrt{S_{p}(\kappa_{1n_{1}}, \kappa_{2n_{2}})\Delta\kappa_{1}\Delta\kappa_{2}}\Bigg[\cos(\kappa_{1n_{1}}x_{1} + \kappa_{2n_{2}}x_{2} + \Phi_{n_{1}n_{2}}^{(1)})\\
    & + \cos(\kappa_{1n_{1}}x_{1} - \kappa_{2n_{2}}x_{2} + \Phi_{n_{1}n_{2}}^{(2)})\\
    & + \sqrt{2}\sum_{i_{1} + j_{1} = n_{1}}^{i_{1} \geq j_{1} \geq 0}\sum_{i_{2} + j_{2} = n_{2}}^{i_{2} \geq j_{2} \geq 0}\sqrt{S(\kappa_{1n_{1}}, \kappa_{2n_{2}})}b_{p}(\kappa_{1i_{1}}, \kappa_{1j_{1}}, \kappa_{2i_{2}}, \kappa_{2j_{2}})\\
    &\Big[ \cos(\kappa_{1n_{1}}x_{1} + \kappa_{2n_{2}}x_{2} + \Phi_{i_{1}i_{2}}^{(1)} + \Phi_{j_{1}j_{2}}^{(1)} + \beta(\kappa_{1i_{1}}, \kappa_{1j_{1}}, \kappa_{2i_{2}}, \kappa_{2j_{2}}))\\
    & + \cos(\kappa_{1n_{1}}x_{1} - \kappa_{2n_{2}}x_{2} + \Phi_{i_{1}i_{2}}^{(2)} + \Phi_{j_{1}j_{2}}^{(2)} + \beta(\kappa_{1i_{1}}, \kappa_{1j_{1}}, -\kappa_{2i_{2}}, -\kappa_{2j_{2}}))\\
    & + \cos(\kappa_{1n_{1}}x_{1} - \kappa_{2i_{2}}x_{2} + \kappa_{2j_{2}}x_{2} + \Phi_{i_{1}i_{2}}^{(2)} + \Phi_{j_{1}j_{2}}^{(1)} + \beta(\kappa_{1i_{1}}, \kappa_{1j_{1}}, -\kappa_{2i_{2}}, \kappa_{2j_{2}}))\\
    & + \cos(\kappa_{1n_{1}}x_{1} + \kappa_{2i_{2}}x_{2} - \kappa_{2j_{2}}x_{2} + \Phi_{i_{1}i_{2}}^{(1)} + \Phi_{j_{1}j_{2}}^{(2)} + \beta(\kappa_{1i_{1}}, \kappa_{1j_{1}}, \kappa_{2i_{2}}, \kappa_{2j_{2}})) \Big] \Bigg]\\
\end{aligned}
\end{equation}
where
\begin{equation}
\begin{aligned}
	& S_{p}(\kappa_{1n_{1}}, \kappa_{2n_{2}}) = S(\kappa_{1n_{1}}, \kappa_{2n_{2}})\Big(1 - \sum_{i_{1} + j_{1} = n_{1}}^{i_{1} \geq j_{1} \geq 0}\sum_{i_{2} + j_{2} = n_{2}}^{i_{2} \geq j_{2} \geq 0}b_{p}^{2}(\kappa_{1i_{1}}, \kappa_{1j_{1}}, \kappa_{2i_{2}}, \kappa_{2j_{2}})\Big)\\
	& b_{p}^{2}(\kappa_{1i_{1}}, \kappa_{2j_{1}}, \kappa_{1i_{2}}, \kappa_{2j_{2}}) = \frac{|B(\kappa_{1i_{1}}, \kappa_{2j_{1}}, \kappa_{1i_{2}}, \kappa_{2j_{2}}|^{2}\Delta\kappa_{1}\Delta\kappa_{2}}{S_{p}(\kappa_{1i_{1}}, \kappa_{2i_{2}})S_{p}(\kappa_{1j_{1}}, \kappa_{2j_{2}})S(\kappa_{1(i_{1} + j_{1})}, \kappa_{2(i_{2} + j_{2})})}\\
\end{aligned}
\end{equation}

\subsection{Simulation of 3-dimensional random fields}

Let $A(x_{1}, x_{2}, x_{3})$ be a three-dimensional uni-variate (3d-1v) random field with zero mean and power spectrum $S(\kappa_{1}, \kappa_{2}, \kappa_{3})$, $2^{nd}$-order autocorrelation function $R_{2}(\xi_{1}, \xi_{2}, \xi_{3})$, bispectrum $B(\kappa_{11}, \kappa_{12},  \kappa_{21}, \kappa_{22}, \kappa_{31}, \kappa_{32})$, and $3^{rd}$-order autocorrelation function $R_{3}(\xi_{11}, \xi_{12}, \xi_{21}, \xi_{22}, \xi_{31},  \xi_{32})$ satisfying:
\begin{equation}
    \mathbb{E}[A(x_{1}, x_{2}, x_{3})] = 0
\end{equation}
\begin{equation}
    \mathbb{E}[A(x_{1} + \xi_{1}, x_{2} + \xi_{2}, x_{3} + \xi_{3})A(x_{1}, x_{2}, x_{3})] = R_{2}(\xi_{1}, \xi_{2}, \xi_{3})
\end{equation}
\begin{equation}
    \mathbb{E}[A(x_{1} + \xi_{11}, x_{2} + \xi_{21}, x_{3} + \xi_{31})A(x_{1} + \xi_{12}, x_{2} + \xi_{22}, x_{3} + \xi_{32})A(x_{1}, x_{2})] = R_{3}(\xi_{11}, \xi_{12}, \xi_{21}, \xi_{22}, \xi_{31}, \xi_{32})
\end{equation}
\begin{equation}
    S(\kappa_{1}, \kappa_{2}, \kappa_{3}) = \frac{1}{(2\pi)^{3}}\int_{-\infty}^{\infty}\int_{-\infty}^{\infty}R_{2}(\xi_{1}, \xi_{2}, \xi_{3})e^{-\iota(\kappa_{1}\xi_{1} + \kappa_{2}\xi_{2} + \kappa_{3}\xi_{3})}d\xi_{1}d\xi_{2}d\xi_{3}
\label{eqn:WK11}
\end{equation}
\begin{equation}
\begin{aligned}
    & B(\kappa_{11}, \kappa_{12}, \kappa_{21}, \kappa_{22}, \kappa_{31}, \kappa_{32}) = \frac{1}{(2\pi)^{6}}\int_{-\infty}^{\infty}\int_{-\infty}^{\infty}\int_{-\infty}^{\infty}\int_{-\infty}^{\infty}\int_{-\infty}^{\infty}\int_{-\infty}^{\infty}\\
    & R_{3}(\xi_{11}, \xi_{12}, \xi_{21}, \xi_{22}, \xi_{31}, \xi_{32}) e^{-\iota(\kappa_{11}\xi_{11} + \kappa_{12}\xi_{12} + \kappa_{21}\xi_{21} + \kappa_{22}\xi_{22} + \kappa_{31}\xi_{31} + \kappa_{32}\xi_{32})}\\
    & d\xi_{11}d\xi_{12}d\xi_{21}d\xi_{22}d\xi_{31}d\xi_{32}\\
\end{aligned}
\label{eqn:WK21}
\end{equation}
\begin{equation}
    R_{2}(\xi_{1}, \xi_{2}, \xi_{3}) = \int_{-\infty}^{\infty}\int_{-\infty}^{\infty}S(\kappa_{1}, \kappa_{2}, \kappa_{3})e^{\iota(\kappa_{1}\xi_{1} + \kappa_{2}\xi_{2} + \kappa_{3}\xi_{3})}d\kappa_{1}d\kappa_{2}d\kappa_{3}
\label{eqn:WK31}
\end{equation}
\begin{equation}
\begin{aligned}
     & R_{3}(\xi_{11}, \xi_{12}, \xi_{21}, \xi_{22}, \xi_{31}, \xi_{32}) = \int_{-\infty}^{\infty}\int_{-\infty}^{\infty}\int_{-\infty}^{\infty}\int_{-\infty}^{\infty}\int_{-\infty}^{\infty}\int_{-\infty}^{\infty}\\
     & B(\kappa_{11}, \kappa_{12}, \kappa_{21}, \kappa_{22}, \kappa_{31}, \kappa_{32})e^{\iota(\kappa_{11}\xi_{11} + \kappa_{12}\xi_{12} + \kappa_{21}\xi_{21} + \kappa_{22}\xi_{22} + \kappa_{31}\xi_{31} + \kappa_{32}\xi_{32})}\\
     & d\kappa_{11}d\kappa_{12}d\kappa_{21}d\kappa_{22}d\kappa_{31}d\kappa_{32}
\end{aligned}
\label{eqn:WK41}
\end{equation}
where Eq.\ \eqref{eqn:WK11}, Eq.\ \eqref{eqn:WK31} and Eq.\ \eqref{eqn:WK21}, Eq.\ \eqref{eqn:WK41} constitute the $2^{nd}$ and $3^{rd}$ order Weiner-Khintchine transform pairs respectively.

The spectra for 3D random fields have the same symmetries presented in Eqs.\ \eqref{eqn:power_spectrum_symmetry} and \eqref{eqn:bispectrum_symmetry_1}--\eqref{eqn:bispectrum_symmetry_3}.

\subsubsection{Simulation formula for general 3-dimensional random fields}
\label{general_3d}

The simulation formula in this case is a straightforward extension of Eq.\ \eqref{eqn:finite_sum_2d} as follows
\begin{equation}
\label{eqn:finite_sum_3d}
\begin{aligned}
    A(x_{1}, x_{2}, x_{3}) &= \sum_{n_{3}=-N_{3}}^{N_{3}}\sum_{n_{2}=-N_{2}}^{N_{2}}\sum_{n_{1}=0}^{N_{1}}[\sqrt{2}A_{pn_{1}n_{2}n_{3}}\cos(\kappa_{1n_{1}}x_{1} + \kappa_{2n_{2}}x_{2} + \kappa_{3n_{3}}x_{3} + \Phi_{n_{1}n_{2}n_{3}})\\
    & + [\sum_{i_{1} + j_{1} = n_{1}}^{i_{1} \geq j_{1} \geq 0}\sum_{i_{2} + j_{2} = n_{2}}^{|n_{2}| \geq|i_{2}| \geq |j_{2}| \geq 0}\sum_{i_{3} + j_{3} = n_{3}}^{|n_{3}| \geq|i_{3}| \geq |j_{3}| \geq 0}\sqrt{2}A_{n_{1}n_{2}n_{3}}b_{p}(\kappa_{1i_{1}}, \kappa_{1j_{1}}, \kappa_{2i_{2}}, \kappa_{2j_{2}}, \kappa_{3i_{3}}, \kappa_{3j_{3}})\\
    &\cos(\kappa_{1n_{1}}x_{1} + \kappa_{2n_{2}}x_{2} + \kappa_{3n_{3}}x_{3} + \Phi_{i_{1}i_{2}i_{3}} + \Phi_{j_{1}j_{2}j_{3}} + \beta(\kappa_{1i_{1}}, \kappa_{1j_{1}}, \kappa_{2i_{2}}, \kappa_{2j_{2}}, \kappa_{3i_{3}}, \kappa_{3j_{3}}))]\\
\end{aligned}
\end{equation}
where
\begin{equation}
\begin{aligned}
    & A_{pn_{1}n_{2}n_{3}} = \sqrt{2S_{p}(\kappa_{1n_{1}}, \kappa_{2n_{2}}, \kappa_{3n_{3}})\Delta\kappa_{1}\Delta\kappa_{2}\Delta\kappa_{3}}\\
    & A_{n_{1}n_{2}n_{3}} = \sqrt{2S(\kappa_{1n_{1}}, \kappa_{2n_{2}}, \kappa_{3n_{3}})\Delta\kappa_{1}\Delta\kappa_{2}}\Delta\kappa_{3}\\
    & S_{p}(\kappa_{1n_{1}}, \kappa_{2n_{2}}, \kappa_{3n_{3}}) = S(\kappa_{1n_{1}}, \kappa_{2n_{2}}, \kappa_{3n_{3}})\Big(1 - \sum_{i_{1} + j_{1} = n_{1}}^{i_{1} \geq j_{1} \geq 0}\sum_{i_{2} + j_{2} = n_{2}}^{|n_{1}| \geq |i_{1}| \geq |j_{1}| \geq 0}\\
    & \sum_{i_{3} + j_{3} = n_{3}}^{|n_{3}| \geq |i_{3}| \geq |j_{3}| \geq 0}b_{p}^{2}(\kappa_{1i_{1}}, \kappa_{1j_{1}}, \kappa_{2i_{2}}, \kappa_{2j_{2}}, \kappa_{3i_{3}}, \kappa_{3j_{3}})\Big)\\
	& b_{p}^{2}(\kappa_{1i_{1}}, \kappa_{2j_{1}}, \kappa_{1i_{2}}, \kappa_{2j_{2}}, \kappa_{3i_{3}}, \kappa_{3j_{3}}) = \frac{|B(\kappa_{1i_{1}}, \kappa_{2j_{1}}, \kappa_{1i_{2}}, \kappa_{2j_{2}}, \kappa_{3i_{3}}, \kappa_{3j_{3}}|^{2}\Delta\kappa_{1}\Delta\kappa_{2}\Delta\kappa_{3}}{S_{p}(\kappa_{1i_{1}}, \kappa_{2i_{2}}, \kappa_{3i_{3}})S_{p}(\kappa_{1j_{1}}, \kappa_{2j_{2}}, \kappa_{3j_{3}})S(\kappa_{1(i_{1} + j_{1})}, \kappa_{2(i_{2} + j_{2})}, \kappa_{3(i_{3} + j_{3})})}\\
	& \kappa_{1n_{1}} = n_{1}\Delta\kappa_{1} \ ; \  \kappa_{2n_{2}} = n_{2}\Delta\kappa_{2} ; \  \kappa_{2n_{2}} = n_{2}\Delta\kappa_{3} \\
	& \Delta\kappa_{1} = \frac{\kappa_{1u}}{N_{1}} \ ; \  \Delta\kappa_{2} = \frac{\kappa_{2u}}{N_{2}} ; \  \Delta\kappa_{3} = \frac{\kappa_{3u}}{N_{3}}\\
\end{aligned}
\end{equation}
and
\begin{equation}
	S(\kappa_{1}, \kappa_{2}, 0) = S(\kappa_{1}, 0, \kappa_{3}) = S(0, \kappa_{2}, \kappa_{3}) = 0 \ \text{for} -\infty \leq \kappa_{1} \leq \infty \; -\infty \leq \kappa_{2} \leq \infty\; -\infty \leq \kappa_{3} \leq \infty
\end{equation}
\begin{equation}
\begin{aligned}
	& B(\kappa_{11}, \kappa_{12}, \kappa_{21}, \kappa_{22}, \kappa_{31}, 0) = B(\kappa_{11}, \kappa_{12}, \kappa_{21}, \kappa_{22}, 0, \kappa_{32}) = B(\kappa_{11}, \kappa_{12}, \kappa_{21}, 0, \kappa_{31}, \kappa_{32})\\
	& = B(\kappa_{11}, \kappa_{12}, 0, \kappa_{22}, \kappa_{31}, \kappa_{32}) = B(\kappa_{11}, 0, \kappa_{21}, \kappa_{22}, \kappa_{31}, \kappa_{32}) = B(0, \kappa_{12}, \kappa_{21}, \kappa_{22}, \kappa_{31}, \kappa_{32}) = 0 \\
	&  \text{for} -\infty \leq \kappa_{11} \leq \infty \ ; -\infty \leq \kappa_{12} \leq \infty \ -\infty \leq \kappa_{21} \leq \infty ; -\infty \leq \kappa_{22} \leq \infty \ \text{and} \\
	&  -\infty \leq \kappa_{31} \leq \infty ; -\infty \leq \kappa_{32} \leq \infty \\
\end{aligned}
\end{equation}
$\kappa_{1u}$, $\kappa_{2u}$ and $\kappa_{3u}$ are the cutoff wave-numbers for the $x_{1}$,$x_{2}$ and $x_{3}$ axis respectively. These cutoff wave-numbers are chosen to satisfy the conditions
\begin{equation}
\begin{aligned}
    & \int_{0}^{\kappa_{1u}}\int_{-\kappa_{2u}}^{\kappa_{2u}}\int_{-\kappa_{3u}}^{\kappa_{3u}}S(\kappa_{1}, \kappa_{2}, \kappa_{3})d\kappa_{1}d\kappa_{2}d\kappa_{3} = (1 - \epsilon)\int_{0}^{\infty}\int_{-\infty}^{\infty}\int_{-\infty}^{\infty}S(\kappa_{1}, \kappa_{2}, \kappa_{3})d\kappa_{1}d\kappa_{2}d\kappa_{3}\\
\end{aligned}
\end{equation}
\begin{equation}
\begin{aligned}
    & \int_{0}^{\kappa_{1u}}\int_{0}^{\kappa_{1u}}\int_{-\kappa_{2u}}^{\kappa_{2u}}\int_{-\kappa_{2u}}^{\kappa_{2u}}\int_{-\kappa_{3u}}^{\kappa_{3u}}\int_{-\kappa_{3u}}^{\kappa_{3u}}B(\kappa_{11}, \kappa_{12}, \kappa_{21}, \kappa_{22}, \kappa_{31}, \kappa_{32})d\kappa_{11}d\kappa_{12}d\kappa_{21}d\kappa_{22}d\kappa_{31}d\kappa_{32}\\
    & = (1 - \epsilon)\int_{0}^{\kappa_{1u}}\int_{0}^{\kappa_{1u}}\int_{-\kappa_{2u}}^{\kappa_{2u}}\int_{-\kappa_{2u}}^{\kappa_{2u}}\int_{-\kappa_{3u}}^{\kappa_{3u}}\int_{-\kappa_{3u}}^{\kappa_{3u}}B(\kappa_{11}, \kappa_{12}, \kappa_{21}, \kappa_{22}, \kappa_{31}, \kappa_{32})d\kappa_{11}d\kappa_{12}d\kappa_{21}d\kappa_{22}d\kappa_{31}d\kappa_{32}\\
\end{aligned}
\end{equation}
where $\epsilon ll 1$. 

The simulated random fields are periodic along $x_{1}$, $x_{2}$ and $x_{3}$ with period
\begin{equation}
\begin{aligned}
    & L_{x_{1}} = \frac{2\pi}{\Delta \kappa_{1}} \\
    & L_{x_{2}} = \frac{2\pi}{\Delta \kappa_{2}} \\
    & L_{x_{3}} = \frac{2\pi}{\Delta \kappa_{3}} \\
\end{aligned}    
\end{equation}
and the aliasing conditions on the spatial increments are as follows
\begin{equation}
\begin{aligned}
    & \Delta x_{1} \leq \frac{2\pi}{2\kappa_{1u}} \\
    & \Delta x_{2} \leq \frac{2\pi}{2\kappa_{2u}} \\
    & \Delta x_{3} \leq \frac{2\pi}{2\kappa_{3u}} \\
\end{aligned}    
\end{equation}

\subsubsection{Simulation formula for quadrant 3-dimensional random fields}
\label{quadrant_3d}

The development of this section mirrors that of Section \ref{quadrant_2d} to present the simulation formula for third-order 3D quadrant random fields. The symmetry in the polyspectra for 3D quadrant random fields are as follows
\begin{equation}
\begin{aligned}
	&S(\kappa_{1}, \kappa_{2}, \dots, \kappa_{3}) = S(I_{1}\kappa_{1}, I_{2}\kappa_{2}, I_{3}\kappa_{3}) \ for \ I_{1},  I_{2}, I_{3} = \pm 1\\
	&B(\kappa_{11}, \kappa_{12}, \kappa_{21}, \kappa_{22}, \kappa_{31}, \kappa_{32}) = B(I_{11}\kappa_{11}, I_{12}\kappa_{12}, I_{21}\kappa_{21}, I_{22}\kappa_{22}, I_{31}\kappa_{31}, I_{32}\kappa_{32})\\
	& \text{for} \ I_{11},  I_{12}, I_{21}, I_{22}, I_{31}, I_{32} = \pm 1\\
\end{aligned}
\end{equation}

Similar to Eq.\ \eqref{eqn:2d_quadrant}, the simulation formula for 3D-1V third-order quadrant random fields simplifies to

\begin{equation}
\label{eqn:3d_quadrant}
\begin{aligned}
    A(x_{1}, x_{2}, x_{3}) &= \sqrt{2}\sum_{n_{3}=0}^{N_{3}}\sum_{n_{2}=0}^{N_{2}}\sum_{n_{1}=0}^{N_{1}}\Bigg[\sqrt{S_{p}(\kappa_{1n_{1}}, \kappa_{2n_{2}}, \kappa_{3n_{3}})\Delta\kappa_{1}\Delta\kappa_{2}\Delta\kappa_{3}}\\
    & \sum_{I_{1}=1, I_{2}=\pm, I_{3}=\pm} \cos(I_{1}\kappa_{1n_{1}}x_{1} + I_{2}\kappa_{2n_{2}}x_{2} + I_{3}\kappa_{3n_{3}}x_{3} + \Phi_{n_{1}n_{2}n_{3}}^{I_{1}I_{2}I_{3}})\\
    & + \sqrt{2}\sum_{i_{1} + j_{1} = n_{1}}^{i_{1} \geq j_{1} \geq 0}\sum_{i_{2} + j_{2} = n_{2}}^{i_{2} \geq j_{2} \geq 0}\sum_{i_{3} + j_{3} = n_{3}}^{i_{3} \geq j_{3} \geq 0}\sqrt{S(\kappa_{1n_{1}}, \kappa_{2n_{2}}, \kappa_{3n_{3}})}b_{p}(\kappa_{1i_{1}}, \kappa_{1j_{1}}, \kappa_{2i_{2}}, \kappa_{2j_{2}}, \kappa_{3i_{3}}, \kappa_{3j_{3}})\\
    &\Big[ \sum_{I_{1}=1, I_{21} = \pm, I_{22} = \pm, I_{31} = \pm, I_{32} = \pm}\cos(I_{1}\kappa_{1n_{1}}x_{1} + I_{21}\kappa_{2i_{2}}x_{2} + I_{22}\kappa_{2j_{2}}x_{2}\\
    & + I_{31}\kappa_{3i_{3}}x_{3} + I_{32}\kappa_{3j_{3}}x_{3} + \Phi_{i_{1}i_{2}i_{3}}^{I_{1}I_{21}I_{31}} + \Phi_{j_{1}j_{2}j_{3}}^{I_{1}I_{22}I_{32}} + \beta(\kappa_{1i_{1}}, \kappa_{1j_{1}}, \kappa_{2i_{2}}, \kappa_{2j_{2}}, \kappa_{3i_{3}}, \kappa_{3j_{3}})) \Big] \Bigg]\\
\end{aligned}
\end{equation}
where
\begin{equation}
\begin{aligned}
	& S_{p}(\kappa_{1n_{1}}, \kappa_{2n_{2}}, \kappa_{3n_{3}}) = S(\kappa_{1n_{1}}, \kappa_{2n_{2}}, \kappa_{3n_{3}})\Big(1 - \sum_{i_{1} + j_{1} = n_{1}}^{i_{1} \geq j_{1} \geq 0}\sum_{i_{2} + j_{2} = n_{2}}^{i_{2} \geq j_{2} \geq 0}\sum_{i_{3} + j_{3} = n_{3}}^{i_{3} \geq j_{3} \geq 0}b_{p}^{2}(\kappa_{1i_{1}}, \kappa_{1j_{1}}, \kappa_{2i_{2}}, \kappa_{2j_{2}}, \kappa_{3i_{3}}, \kappa_{3j_{3}})\Big)\\
	& b_{p}^{2}(\kappa_{1i_{1}}, \kappa_{1j_{1}}, \kappa_{2i_{2}}, \kappa_{2j_{2}}, \kappa_{3i_{3}}, \kappa_{3j_{3}}) = \frac{|B(\kappa_{1i_{1}}, \kappa_{1j_{1}}, \kappa_{2i_{1}}, \kappa_{2j_{2}}, \kappa_{3i_{3}}, \kappa_{3j_{3}}|^{2}\Delta\kappa_{1}\Delta\kappa_{2}\Delta\kappa_{3}}{S_{p}(\kappa_{1i_{1}}, \kappa_{2i_{2}}, \kappa_{3i_{3}})S_{p}(\kappa_{1j_{1}}, \kappa_{2j_{2}}, \kappa_{3j_{3}})S(\kappa_{1(i_{1} + j_{1})}, \kappa_{2(i_{2} + j_{2})}, , \kappa_{3(i_{3} + j_{3})})}\\
\end{aligned}
\end{equation}

\subsection{Simulation of d-dimensional random fields}
\label{theory_dd}

Let $A(x_{1}, x_{2}, \dots, x_{d})$ be a $d$-dimensional uni-variate ($d$D-1V) random field with zero mean, power spectrum $S(\kappa_{1}, \kappa_{2}, \dots, \kappa_{d})$, $2^{nd}$-order autocorrelation function $R_{2}(\xi_{1}, \xi_{2}, \dots, \xi_{d})$, bispectrum $B(\kappa_{11}, \kappa_{12}, \kappa_{21}, \kappa_{22}, \kappa_{31}, , \kappa_{32}, \dots, \dots, \kappa_{n1}, \kappa_{n2})$, and $3^{rd}$-order Autocorrelation function $R_{3}(\xi_{11}, \xi_{12}, \xi_{21}, \xi_{22}, \xi_{31}, \xi_{32}, \dots, \dots, \xi_{n1}, \xi_{n2})$. For convenience, let us define the following new vector quantities:
\begin{equation}
\begin{aligned}
	&\text{Position vector:} \ \overline{x} = [x_{1}, x_{2}, \dots, x_{n}]^{T}\\
	&\text{Seperation vector:} \ \overline{\xi} = [\xi_{1}, \xi_{2}, \dots, \xi_{n}]^{T}\\
	&\text{Wave number vector:} \ \overline{\kappa} = [\kappa_{1}, \kappa_{2}, \dots, \kappa_{n}]^{T}\\
\end{aligned}
\label{eqn:vector_notation}
\end{equation}

The random field $A(\overline{x})$ is third-order stationary satisfying the following conditions:
\begin{equation}
    \mathbb{E}[A(\overline{x})] = 0
\end{equation}
\begin{equation}
    \mathbb{E}[A(\overline{x} + \overline{\xi})A(\overline{x})] = R_{2}(\overline{\xi})
\end{equation}
\begin{equation}
    \mathbb{E}[A(\overline{x} + \overline{\xi}_{1})A(\overline{x} + \overline{\xi}_{2})A(\overline{x} + \overline{\xi})] = R_{3}(\overline{\xi}_{1}, \overline{\xi}_{2})
\end{equation}
\begin{equation}
    S(\overline{\kappa}) = \frac{1}{(2\pi)^{n}}\int_{-\infty}^{\infty}R_{2}(\overline{\xi})e^{-\iota(\overline{\kappa}\overline{\xi})}d\overline{\xi}
\label{eqn:WK1n}
\end{equation}
\begin{equation}
\begin{aligned}
    & B(\overline{\kappa}_{1}, \overline{\kappa}_{2}) = \frac{1}{(2\pi)^{2n}}\int_{-\infty}^{\infty}R_{3}(\overline{\xi}_{1}, \overline{\xi}_{2}) e^{-\iota(\overline{\kappa}_{1}\overline{\xi}_{1} + \overline{\kappa}_{2}\overline{\xi}_{2}})d\overline{\xi}_{1}d\overline{\xi}_{2}\\
\end{aligned}
\label{eqn:WK2n}
\end{equation}
\begin{equation}
    R_{2}(\overline{\xi}) = \int_{-\infty}^{\infty}S(\overline{\kappa})e^{\iota(\overline{\kappa}\overline{\xi})}d\overline{\kappa}
\label{eqn:WK3n}
\end{equation}
\begin{equation}
\begin{aligned}
     & R_{3}(\overline{\kappa}_{1}, \overline{\kappa}_{2}) = \int_{-\infty}^{\infty}\int_{-\infty}^{\infty}B(\overline{\kappa}_{1}, \overline{\kappa}_{2})e^{\iota(\overline{\kappa}_{1} \cdot \overline{\xi}_{1} + \overline{\kappa}_{2} \cdot \overline{\xi}_{2})}d\overline{\kappa}_{1}\overline{\kappa}_{2}
\end{aligned}
\label{eqn:WK4n}
\end{equation}
where Eq.\ \eqref{eqn:WK1n}, Eq.\ \eqref{eqn:WK3n} and Eq.\ \eqref{eqn:WK2n}, Eq.\ \eqref{eqn:WK4n} constitute the $2^{nd}$ and $3^{rd}$ order Weiner-Khintchine transform pairs respectively. The symmetries in Eqs.\ \eqref{eqn:power_spectrum_symmetry} and \eqref{eqn:bispectrum_symmetry_1} -- \eqref{eqn:bispectrum_symmetry_3} still hold here.

\subsubsection{Simulation formula for general d-dimensional random fields}
\label{general_dd}

The formula for the simulation of general $d$-dimensional random fields follows Sec. \ref{general_2d} and \ref{general_3d} closely as
\begin{equation}
\label{eqn:finite_sum_dd}
\begin{aligned}
    A(\overline{x}) &= \sum_{n_{d}=-N_{d}}^{N_{d}} \dots \sum_{n_{2}=-N_{2}}^{N_{2}}\sum_{n_{1}=0}^{N_{1}}[\sqrt{2}A_{p\overline{n}}\cos(\overline{\kappa}\cdot \overline{x} + \Phi_{\overline{n}})\\
    & + \sum_{i_{1} + j_{1} = n_{1}}^{i_{1} \geq j_{1} \geq 0}\sum_{i_{2} + j_{2} = n_{2}}^{|n_{2}| \geq|i_{2}| \geq |j_{2}| \geq 0} \dots \sum_{i_{d} + j_{d} = n_{d}}^{|n_{d}| \geq|i_{d}| \geq |j_{d}| \geq 0}\sqrt{2}A_{\overline{n}}b_{p}(\overline{\kappa}_{\overline{i}}, \overline{\kappa}_{\overline{j}}) \cos(\overline{\kappa} \cdot \overline{x} + \Phi_{\overline{i}} + \Phi_{\overline{j}} + \beta(\overline{\kappa}_{\overline{i}}, \overline{\kappa}_{\overline{j}}))]\\
\end{aligned}
\end{equation}
where
\begin{equation}
\begin{aligned}
    & A_{p\overline{n}} = \sqrt{2S_{p}(\overline{\kappa}_{\overline{n}})\Delta\kappa_{1}\Delta\kappa_{2} \dots \Delta\kappa_{d}}\\
    & A_{\overline{n}} = \sqrt{2S(\overline{\kappa}_{\overline{n}})\Delta\kappa_{1}\Delta\kappa_{2} \dots \Delta\kappa_{d}}\\
    & S_{p}(\overline{\kappa}_{\overline{n}}) = S(\overline{\kappa}_{\overline{n}})\Big(1 - \sum_{i_{1} + j_{1} = n_{1}}^{i_{1} \geq j_{1} \geq 0}\sum_{i_{2} + j_{2} = n_{2}}^{|n_{1}| \geq |i_{1}| \geq |j_{1}| \geq 0} \dots \sum_{i_{d} + j_{d} = n_{d}}^{|n_{d}| \geq |i_{d}| \geq |j_{d}| \geq 0} b_{p}^{2}(\overline{\kappa}_{i}, \overline{\kappa}_{j})\Big)\\
	& b_{p}^{2}(\overline{\kappa}_{\overline{i}}, \overline{\kappa}_{\overline{j}}) = \frac{|B(\overline{\kappa}_{\overline{i}}, \overline{\kappa}_{\overline{j}})|^{2}\Delta\kappa_{1}\Delta\kappa_{2} \dots \Delta\kappa_{d}}{S_{p}(\overline{\kappa}_{\overline{i}})S_{p}(\overline{\kappa}_{\overline{j}})S(\overline{\kappa}_{\overline{n}})}\\
	& \kappa_{1n_{1}} = n_{1}\Delta\kappa_{1} \ ; \  \kappa_{2n_{2}} = n_{2}\Delta\kappa_{2} ; \dots;  \kappa_{dn_{d}} = n_{d}\Delta\kappa_{d} \\
	& \Delta\kappa_{1} = \frac{\kappa_{1u}}{N_{1}} \ ; \  \Delta\kappa_{2} = \frac{\kappa_{2u}}{N_{2}} ; \dots ;  \Delta\kappa_{d} = \frac{\kappa_{du}}{N_{d}}\\
\end{aligned}
\end{equation}
and
\begin{equation}
\begin{aligned}
	& S(0, \kappa_{2}, \dots, \kappa_{d}) = S(\kappa_{1}, 0,\dots, \kappa_{d}) = S(\kappa_{1}, \kappa_{2},\dots, 0) = 0 \\
	& \ \text{for} -\infty \leq \kappa_{1} \leq \infty \; -\infty \leq \kappa_{2} \leq \infty; \dots; -\infty \leq \kappa_{d} \leq \infty \\
\end{aligned}
\end{equation}

\begin{equation}
\begin{aligned}
	& B(0, \kappa_{12}, \dots \kappa_{d1}, \kappa_{d2}) = B(\kappa_{11}, 0, \dots \kappa_{d1}, \kappa_{d2}) = \dots = B(\kappa_{11}, \kappa_{12}, \dots 0, \kappa_{d2}) = B(\kappa_{11}, \kappa_{12}, \dots \kappa_{d1}, 0) = 0 \\
	&  \text{for} -\infty \leq \kappa_{11} \leq \infty \ ; -\infty \leq \kappa_{12} \leq \infty; \dots; -\infty \leq \kappa_{d1} \leq \infty ; -\infty \leq \kappa_{d2} \leq \infty \\
\end{aligned}
\end{equation}
In the above expressions, the overline subscripts denote the iterable index sets $\overline{n}=\{n_1,n_2,\dots,n_d\}$, $\overline{i}=\{i_1,i_2,\dots,i_d\}$, and $\overline{j}=\{j_1,j_2,\dots,j_d\}$. In particular, $\Phi_{\overline{n}}$ denotes the $d^{th}$-order tensor of random phase angles indexed as $\Phi_{n_1n_2\dots n_d}$ and $A_{p\overline{n}},A_{\overline{n}}$ denote $d^{th}$-order tensors of amplitudes having components $A_{pn_1n_2\dots n_d}, A_{n_1n_2\dots n_d}$. Indexing of the wave number combines the vector overline notations of Eq.\ \eqref{eqn:vector_notation} with the overline subscripts such that $\overline{\kappa}_{\overline{n}}$ denotes the wave number set $(\kappa_{1n_1},\kappa_{2n_2},\dots,\kappa_{dn_d})$. Finally, $\kappa_{1u}$, $\kappa_{2u}$, \dots and $\kappa_{du}$ are the cutoff wave-numbers for the $x_{1}$, $x_{2}$ \dots $x_{d}$ axes respectively, satisfying
\begin{equation}
\begin{aligned}
    & \int_{-\overline{\kappa}_{u}}^{\overline{\kappa}_{u}}S(\overline{\kappa})d\overline{\kappa} = (1 - \epsilon)\int_{-\overline{\kappa}_{u}}^{\overline{\kappa}_{u}}S(\overline{\kappa})d\overline{\kappa}\\
\end{aligned}
\end{equation}
\begin{equation}
\begin{aligned}
    & \int_{-\overline{\kappa}_{\overline{i}u}}^{\overline{\kappa}_{\overline{i}u}}\int_{-\overline{\kappa}_{\overline{j}u}}^{\overline{\kappa}_{\overline{j}u}}B(\overline{\kappa}_{\overline{i}}, \overline{\kappa}_{\overline{j}})d\overline{\kappa}_{\overline{i}}d\overline{\kappa}_{\overline{j}} = (1 - \epsilon)\int_{-\overline{\kappa}_{\overline{i}u}}^{\overline{\kappa}_{\overline{i}u}}\int_{-\overline{\kappa}_{\overline{j}u}}^{\overline{\kappa}_{\overline{j}u}}B(\overline{\kappa}_{\overline{i}}, \overline{\kappa}_{\overline{j}})d\overline{\kappa}_{\overline{i}}d\overline{\kappa}_{\overline{j}}\\
\end{aligned}
\end{equation}
where $\epsilon \ll 1$. 

The simulated random fields are periodic along the $x_{1}$, $x_{2}$ \dots $x_{d}$ with period
\begin{equation}
\begin{aligned}
    & L_{x_{1}} = \frac{2\pi}{\Delta \kappa_{1}} \\
    & L_{x_{2}} = \frac{2\pi}{\Delta \kappa_{2}} \\
    & \vdots \\
    & L_{x_{d}} = \frac{2\pi}{\Delta \kappa_{d}} \\
\end{aligned}    
\end{equation}
and the conditions to prevent aliasing are given as
\begin{equation}
\begin{aligned}
    & \Delta x_{1} \leq \frac{2\pi}{2\kappa_{1u}} \\
    & \Delta x_{2} \leq \frac{2\pi}{2\kappa_{2u}} \\
    & \vdots \\
    & \Delta x_{d} \leq \frac{2\pi}{2\kappa_{du}} \\
\end{aligned}    
\label{eqn:aliasing}
\end{equation}

\subsubsection{Simulation formula for $d$-dimensional quadrant random fields}
\label{quadrant_dd}

The development of this section mirrors that of Sections \ref{quadrant_2d} - \ref{quadrant_3d}. The symmetry in the polyspectra in the case of quadrant random fields is given by
\begin{equation}
\begin{aligned}
	&S(\kappa_{1}, \kappa_{2}, \dots, \kappa_{d}) = S(I_{1}\kappa_{1}, I_{2}\kappa_{2}, \dots, I_{d}\kappa_{d}) \ \text{for} \ I_{1},  I_{2}, \dots I_{d} = \pm 1\\
	&B(\kappa_{11}, \kappa_{12}, \kappa_{21}, \kappa_{22}, \dots, \dots, \kappa_{d1}, \kappa_{d2}) = B(I_{11}\kappa_{11}, I_{12}\kappa_{12}, I_{21}\kappa_{21}, I_{22}\kappa_{22},  \dots, \dots, I_{d1}\kappa_{d1}, I_{d2}\kappa_{d2})\\
	& \text{for} \ I_{11},  I_{12}, I_{21}, I_{22}, \dots, \dots, I_{d1}, I_{d2} = \pm 1\\
\end{aligned}
\end{equation}

As a result, the simulation formula for $d$D-1V third-order quadrant random fields simplifies to (similar to Eq.\ \eqref{eqn:2d_quadrant} and \eqref{eqn:3d_quadrant})
\begin{equation}
\label{eqn:dd_quadrant}
\begin{aligned}
    & A(x_{1}, x_{2},\dots, x_{d}) = \sqrt{2}\sum_{n_{d}=0}^{N_{d}} \dots \sum_{n_{2}=0}^{N_{2}}\sum_{n_{1}=0}^{N_{1}}\Bigg[\sqrt{S_{p}(\kappa_{1n_{1}}, \kappa_{2n_{2}}, \dots, \kappa_{dn_{d}})\Delta\kappa_{1}\Delta\kappa_{2} \dots \Delta\kappa_{d}}\\
    & \sum_{I_{1}=1, I_{2}=\pm1,\dots, I_{d}=\pm1} \cos(I_{1}\kappa_{1n_{1}}x_{1} + I_{2}\kappa_{2n_{2}}x_{2} + \dots + I_{d}\kappa_{dn_{d}}x_{3} + \Phi_{n_{1}n_{2} \dots n_{d}}^{I_{1}I_{2} \dots I_{d}})\\
    & + \sqrt{2}\sum_{i_{1} + j_{1} = n_{1}}^{i_{1} \geq j_{1} \geq 0}\sum_{i_{2} + j_{2} = n_{2}}^{i_{2} \geq j_{2} \geq 0} \dots \sum_{i_{d} + j_{d} = n_{d}}^{i_{d} \geq j_{d} \geq 0}\sqrt{S(\kappa_{1n_{1}}, \kappa_{2n_{2}}, \dots, \kappa_{dn_{d}})}b_{p}(\kappa_{1i_{1}}, \kappa_{1j_{1}}, \kappa_{2i_{2}}, \kappa_{2j_{2}},\dots,  \kappa_{di_{d}}, \kappa_{dj_{d}})\\
    &\Big[ \sum_{I_{1}=1, I_{21} = \pm1, I_{22} = \pm1, \dots I_{d1} = \pm1, I_{d2} = \pm1}\cos(I_{1}\kappa_{1i_{1}}x_{1} + I_{21}\kappa_{2i_{2}}x_{2} + I_{22}\kappa_{2j_{2}}x_{2} + \dots \\
    & + I_{d1}\kappa_{di_{d}}x_{d} + I_{d2}\kappa_{dj_{d}}x_{d} + \Phi_{i_{1}i_{2} \dots i_{d}}^{I_{1}I_{21} \dots I_{d1}} + \Phi_{j_{1}j_{2} \dots j_{d}}^{I_{1}I_{22} \dots I_{d2}} + \beta(\kappa_{1i_{1}}, \kappa_{1j_{1}}, \kappa_{2i_{2}}, \kappa_{2j_{2}}, \dots, \kappa_{di_{d}}, \kappa_{dj_{d}})) \Big] \Bigg]\\
\end{aligned}
\end{equation}
where
\begin{equation}
\begin{aligned}
	& S_{p}(\kappa_{1n_{1}}, \kappa_{2n_{2}}, \dots \kappa_{dn_{d}}) = S(\kappa_{1n_{1}}, \kappa_{2n_{2}}, \dots, \kappa_{dn_{d}})\Big(1 - \\
	& \sum_{i_{1} + j_{1} = n_{1}}^{i_{1} \geq j_{1} \geq 0}\sum_{i_{2} + j_{2} = n_{2}}^{i_{2} \geq j_{2} \geq 0} \dots \sum_{i_{d} + j_{d} = n_{d}}^{i_{d} \geq j_{d} \geq 0}b_{p}^{2}(\kappa_{1i_{1}}, \kappa_{1j_{1}}, \kappa_{2i_{2}}, \kappa_{2j_{2}}, \dots, \kappa_{di_{d}}, \kappa_{dj_{d}})\Big)\\
	& b_{p}^{2}(\kappa_{1i_{1}}, \kappa_{1j_{1}}, \kappa_{2i_{2}}, \kappa_{2j_{2}}, \dots, \kappa_{di_{d}}, \kappa_{dj_{d}}) = \\
	& \frac{|B(\kappa_{1i_{1}}, \kappa_{1j_{1}}, \kappa_{2i_{1}}, \kappa_{2j_{2}}, \dots,  \kappa_{di_{d}}, \kappa_{dj_{d}}|^{2}\Delta\kappa_{1}\Delta\kappa_{2} \dots \Delta\kappa_{d}}{S_{p}(\kappa_{1i_{1}}, \kappa_{2i_{2}}, \dots,  \kappa_{di_{d}})S_{p}(\kappa_{1j_{1}}, \kappa_{2j_{2}}, \dots,  \kappa_{dj_{d}})S(\kappa_{1(i_{1} + j_{1})}, \kappa_{2(i_{2} + j_{2})}, \dots, \kappa_{d(i_{d} + j_{d})})}\\
\end{aligned}
\end{equation}
Note that we forego the overline index set notation in lieu of the full indicial notation given the introduction of additional summations associated with the symmetries. For simulation purposes, we further note that the quadrant random fields require the generation of $2^d$ sets of $d^{th}$-order tensors of random phase angles.

\section{Third-order Spectral Representation Method with Fast Fourier Transform}

The simulation formulae presented up until now can be used for simulating random fields, but the simulations in their current form grow increasingly computational intensive with increasing dimension; so much so that simulating 3-dimensional random fields becomes impractical. In fully discretized form, the simulation formula for a 1D-1V third-order random field is given by
\begin{equation}
\begin{aligned}
    A(m \Delta x) = &\sqrt{2}\sum_{n=0}^{N}\sqrt{2S_{P}(n \Delta \kappa) \Delta \kappa} \cos \big((n \Delta \kappa)(m \Delta x) - \phi_{n}\big) \\
    &+ \sqrt{2}\sum_{n=0}^{\infty}\sum_{i+j=n}^{i \geq j \geq 0}\sqrt{2S(n \Delta\kappa)\Delta\kappa}|b(\kappa_{i}, \kappa_{j})| \\
    & \cos \big((n \Delta\kappa)(m \Delta x) - (\phi_{i} + \phi_{j} + \beta (\kappa_{i},\kappa_{j}))\big), \hspace{6pt} m=1,\dots,M
\end{aligned}
\label{eqn:1D-1v_discrete}
\end{equation}
Assuming that all required data such as partial bicoherences, biphase angles, etc. have been computed a priori, sample function generation has complexity $O(MN)$ for this 1-dimensional case. This complexity increases exponentially for multi-dimensional random fields to order $O((MN)^{d})$ where $d$ is the dimension of the random field. Here, we introduce a fast Fourier transform (FFT) based implementation to reduce the complexity of the simulations.

We first develop an FFT based implementation for simulation of 1D-1V third-order random fields and subsequently extend it to the $2$D-1V, $3$D-1V, and $d$D-1V cases. Let us start by writing Eq.\ \eqref{eqn:1D-1v_discrete} in its complete form
\begin{equation}
\begin{aligned}
A(m \Delta x) = &\sqrt{2}\sum_{n=0}^{\infty}\sqrt{2S(k \Delta \kappa) \Delta \kappa(1 - \sum_{i+j=n}^{i \geq j \geq 0} b_{p}^{2}(\kappa_{i}, \kappa_{j}))} \\
& \cos \big((n \Delta \kappa)(m \Delta x) - \phi_{k}\big) \\
&+ \sqrt{2}\sum_{n=0}^{\infty}\sum_{i+j=n}^{i \geq j \geq 0}\sqrt{2S(n \Delta \kappa)\Delta \kappa} |b_{p}(\kappa_{i}, \kappa_{j})| \\
& \cos \big((n \Delta \kappa)(m \Delta x) - (\phi_{i} + \phi_{j} + \beta (\kappa_{i},\kappa_{j}))\big)
\end{aligned}
\end{equation}
Simplifying the representation from two additive terms to only one term we get
\begin{equation}
\begin{aligned}
A(m \Delta x) = &\sqrt{2}\sum_{n=0}^{\infty}\sqrt{2S(n \Delta \kappa) \Delta \kappa }\\
& \Big[ \sqrt{(1 - \sum_{i+j=n}^{i \geq j \geq 0} b_{p}^{2}(\kappa_{i}, \kappa_{j}))} \cos \big((n \Delta \kappa) (m \Delta x) - \phi_{n}\big)\\
& + \sum_{i+j=n}^{i \geq j \geq 0} |b_{p}(\kappa_{i}, \kappa_{j})| \cos \big((n \Delta \kappa)(m \Delta x) - (\phi_{i} + \phi_{j} + \beta (\kappa_{i},\kappa_{j}))\big) \Big]
\end{aligned}
\end{equation}
From Euler's notation we have that $e^{i\phi} = \cos\phi + \iota\sin\phi$, hence $\cos\phi = \Re[e^{\iota\phi}]$. Applying Euler's notation, we have
\begin{equation}
\begin{aligned}
A(m \Delta x) = &\sqrt{2}\sum_{n=0}^{\infty}\sqrt{2S(n \Delta \kappa) \Delta \kappa }\\
& \Re{\Big[ \sqrt{(1 - \sum_{i+j=n}^{i \geq j \geq 0} b_{p}^{2}(\kappa_{i}, \kappa_{j}))} e^{\iota\big((n \Delta \kappa) (m \Delta x) - \phi_{n}\big)}} \\
& +\sum_{i+j=k}^{i \geq j \geq 0} |b_{p}(\kappa_{i}, \kappa_{j})| e^{\iota\big((n \Delta \kappa)(m \Delta x) - (\phi_{i} + \phi_{j} + \beta (\kappa_{i},\kappa_{j}))\big)} \Big]\
\end{aligned}
\end{equation}
which can be factored as
\begin{equation}
\begin{aligned}
A(m \Delta x) = &\sqrt{2}\sum_{n=0}^{\infty}\sqrt{2S(n \Delta \kappa) \Delta \kappa }\\
&\Re{\Big[\Big(e^{\iota\big((n \Delta \kappa) (m \Delta x)\big)}\Big)} \Big(\sqrt{(1 - \sum_{i+j=n}^{i \geq j \geq 0} b_{p}^{2}(\kappa_{i}, \kappa_{j}))} e^{-\iota\phi_{n}} \\
& +\sum_{i+j=n}^{i \geq j \geq 0} |b_{p}(\kappa_{i}, \kappa_{j})| e^{-\iota(\phi_{i} + \phi_{j} + \beta (\kappa_{i},\kappa_{j}))} \Big) \Big]\
\end{aligned}
\label{eqn:FFT_long}
\end{equation}

The standard form for implementation of the FFT is given by \cite{Cooley2006}:
\begin{equation}
\begin{aligned}
A_{m} = \sum_{n=0}^{N-1}B_{n}e^{-2 \pi \iota \frac{mn}{N}}
\end{aligned}
\label{eqn:FFT_standard}
\end{equation}
For illustration of the implementation here, we will adopt the following shorthand notation. Let $\bm{A}=\{A_m; m=0,\dots\,M-1\}$ where $A_m=A(m\Delta x)$ and $\bm{B}=\{B_n; n=0,\dots,N-1\}$ where $B_n=B(n\Delta\kappa)$, then the fast Fourier transform will be expressed as $\bm{A}=\text{FFT}(\bm{B})$. Similarly, the inverse FFT is denoted $\bm{A}=\text{IFFT}(\bm{B})$. Practically speaking, the FFT implementation involves typically a $\dfrac{1}{N}$ normalization term and therefore inverse FFT requires a multiplication by $N$.

By grouping terms in Eq.\ \eqref{eqn:FFT_long} as follows,
\begin{equation}
\begin{aligned}
A(m \Delta x) & = \Re{\Big[\sum_{n=0}^{\infty}B_{n}e^{i\big((n \Delta \kappa) (m \Delta x)\big)}\Big]}\\
& \text{where} \: B_{n} = \sqrt{2} C_{n}\Big[ \sqrt{(1 - \sum_{i+j=n}^{i \geq j \geq 0} b_{p}^{2}(\kappa_{i}, \kappa_{j}))} e^{i\phi_{n}} +\sum_{i+j=n}^{i \geq j \geq 0} |b_{p}(\kappa_{i}, \kappa_{j})| e^{i(\phi_{i} + \phi_{j} + \beta (\kappa_{i},\kappa_{j}))} \Big]\\
& C_{n} = \sqrt{2S(n \Delta \kappa) \Delta \kappa }
\end{aligned}
\end{equation}
we see that the simulation formula in Eq.\ \eqref{eqn:FFT_long} can be expressed in the compact form of the FFT operator in Eq.\ \eqref{eqn:FFT_standard}. More specifically,
\begin{equation}
    \bm{A} = \Re \{N\text{IFFT}(\bm{B})\}
\end{equation}



\subsection{FFT implementation for the simulation of 2-dimensional random fields}


Next, we develop the FFT implementation for 2D-1V third-order random fields and subsequently extend it to the simulation of $d$D-1V random fields. The discretized simulation formula for 2-dimensional random fields is given by
\begin{equation}
\begin{aligned}
    A(m_{1} \Delta x_{1}, m_{2} \Delta x_{2}) &= \sum_{n_{2}=-N_{2}}^{N_{2}}\sum_{n_{1}=0}^{N_1}[\sqrt{2}\sqrt{2S_{p}(n_{1} \Delta \kappa_{1}, n_{2} \Delta \kappa_{2})\Delta\kappa_{1}\Delta\kappa_{2}}\\
    & \cos(n_{1}\Delta \kappa_{1} m_{1} \Delta x_{1} + n_{2}\Delta \kappa_{2} m_{2} \Delta x_{1} + \Phi_{n_{1}n_{2}})\\
    & + [\sum_{i_{1} + j_{1} = n_{1}}^{i_{1} \geq j_{1} \geq 0}\sum_{i_{2} + j_{2} = n_{2}}^{|n_{2}| \geq |i_{2}| \geq |j_{2}| \geq 0}\sqrt{2}\sqrt{2S(n_{1} \Delta \kappa_{1}, n_{2} \Delta \kappa_{2})\Delta\kappa_{1}\Delta\kappa_{2}}\\
    & b_{p}(i_{1} \Delta \kappa_{1}, j_{1} \Delta \kappa_{1}, i_{2} \Delta \kappa_{2}, j_{2} \Delta \kappa_{2}) \cos(n_{1}\Delta \kappa_{1} m_{1} \Delta x_{1} + n_{2}\Delta \kappa_{2} m_{2} \Delta x_{2} + \Phi_{i_{1}i_{2}} + \Phi_{j_{1}j_{2}} + \\
    & \beta(i_{1} \Delta \kappa_{1}, j_{1} \Delta \kappa_{1}, i_{2} \Delta \kappa_{2}, j_{2} \Delta \kappa_{2}))]\\
\label{eqn:discrete_2d}
\end{aligned}
\end{equation}
where
\begin{equation}
\begin{aligned}
    & S_{p}(n_{1} \Delta \kappa_{1}, n_{2} \Delta \kappa_{2}) = S(n_{1} \Delta \kappa_{1}, n_{2} \Delta \kappa_{2})\Big(1 - \sum_{i_{1} + j_{1} = n_{1}}^{i_{1} \geq j_{1} \geq 0}\sum_{i_{2} + j_{2} = n_{2}}^{|n_{2}| \geq |i_{2}| \geq |j_{2}| \geq 0}b_{p}^{2}(i_{1} \Delta \kappa_{1}, j_{1} \Delta \kappa_{1}, i_{2} \Delta \kappa_{2}, j_{2} \Delta \kappa_{2})\Big)\\
	& b_{p}^{2}(i_{1} \Delta \kappa_{1}, j_{1} \Delta \kappa_{1}, i_{2} \Delta \kappa_{2}, j_{2} \Delta \kappa_{2}) = \frac{|B(i_{1} \Delta \kappa_{1}, j_{1} \Delta \kappa_{1}, i_{2} \Delta \kappa_{2}, j_{2} \Delta \kappa_{2})|^{2}\Delta\kappa_{1}\Delta\kappa_{2}}{S_{p}(i_{1} \Delta \kappa_{1}, i_{2} \Delta \kappa_{2})S_{p}(j_{1} \Delta\kappa_{1}, j_{2} \Delta \kappa_{2})S((i_{1} + j_{1}) \Delta \kappa_{1}, (i_{2} + j_{2}) \Delta \kappa_{2})}\\
\end{aligned}
\end{equation}

Following similar steps as in the 1D-1V FFT implementation and applying symmetry to sum over only positive indices we have
\begin{equation}
\begin{aligned}
    A(m_{1} \Delta x_{1}, m_{2} \Delta x_{2}) &=2  \sum_{n_{2}=0}^{N_{2}}\sum_{n_{1}=0}^{N_1}\sqrt{S(n_{1} \Delta \kappa_{1}, n_{2} \Delta \kappa_{2})\Delta\kappa_{1}\Delta\kappa_{2}}\\
    & \Big[ \sqrt{(1 - \sum_{i_{1} + j_{1} = n_{1}}^{i_{1} \geq j_{1} \geq 0}\sum_{i_{2} + j_{2} = n_{2}}^{n_{2} \geq |i_{2}| \geq |j_{2}| \geq 0}b_{p}^{2}(i_{1} \Delta \kappa_{1}, j_{1} \Delta \kappa_{1}, i_{2} \Delta \kappa_{2}, j_{2} \Delta \kappa_{2})}e^{\iota\Phi_{n_{1}n_{2}}^{(1)}}\\
    & + [\sum_{i_{1} + j_{1} = n_{1}}^{i_{1} \geq j_{1} \geq 0}\sum_{i_{2} + j_{2} = n_{2}}^{n_{2} \geq |i_{2}| \geq |j_{2}| \geq 0} b_{p}(i_{1} \Delta \kappa_{1}, j_{1} \Delta \kappa_{1}, i_{2} \Delta \kappa_{2}, j_{2} \Delta \kappa_{2})\\
    & e^{\iota(\Phi_{i_{1}i_{2}}^{(1|2)} + \Phi_{j_{1}j_{2}}^{(1|2)} +\beta(i_{1} \Delta \kappa_{1}, j_{1} \Delta \kappa_{1}, i_{2} \Delta \kappa_{2}, j_{2} \Delta \kappa_{2}))}] \Big] e^{\iota(n_{1}\Delta \kappa_{1} m_{1} \Delta x_{1} + n_{2}\Delta \kappa_{2} m_{2} \Delta x_{2})}\\
    &  + \sum_{n_{2}=0}^{N_{2}}\sum_{n_{1}=0}^{N_1}\sqrt{S(n_{1} \Delta \kappa_{1}, -n_{2} \Delta \kappa_{2})\Delta\kappa_{1}\Delta\kappa_{2}}\\
    & \Big[\sqrt{(1 - \sum_{i_{1} + j_{1} = n_{1}}^{i_{1} \geq j_{1} \geq 0}\sum_{i_{2} + j_{2} = -n_{2}}^{n_{2} \geq |i_{2}| \geq |j_{2}| \geq 0}b_{p}^{2}(i_{1} \Delta \kappa_{1}, j_{1} \Delta \kappa_{1}, i_{2} \Delta \kappa_{2}, j_{2} \Delta \kappa_{2})}e^{\iota\Phi_{n_{1}n_{2}}^{(2)}}\\
    & + [\sum_{i_{1} + j_{1} = n_{1}}^{i_{1} \geq j_{1} \geq 0}\sum_{i_{2} + j_{2} = -n_{2}}^{n_{2} \geq |i_{2}| \geq |j_{2}| \geq 0} b_{p}(i_{1} \Delta \kappa_{1}, j_{1} \Delta \kappa_{1}, i_{2} \Delta \kappa_{2}, j_{2} \Delta \kappa_{2})\\
    & e^{\iota(\Phi_{i_{1}i_{2}}^{(1|2)} + \Phi_{j_{1}j_{2}}^{(1|2)} +\beta(i_{1} \Delta \kappa_{1}, j_{1} \Delta \kappa_{1}, i_{2} \Delta \kappa_{2}, j_{2} \Delta \kappa_{2}))}] \Big] e^{\iota(n_{1}\Delta \kappa_{1} m_{1} \Delta x_{1} - n_{2}\Delta \kappa_{2} m_{2} \Delta x_{2})}
\label{eqn:fft_2d_2_parts}
\end{aligned}
\end{equation}
where $\Phi^{(1|2)}$ represents the appropriate phase angles, $\Phi^{(1)}$ or $\Phi^{(2)}$, selected as: $\Phi^{(1)}$ for $i_{2}$ or $j_{2} > 0$, and $\Phi^{(2)}$ for $i_{2}$ or $j_{2} < 0$. This can be simplified to
 \begin{equation}
\begin{aligned}
    A(m_{1} \Delta x_{1}, m_{2} \Delta x_{2}) &=2  \sum_{n_{2}=0}^{N_{2}}\sum_{n_{1}=0}^{N_1} \Big[ B_{n_{1}n_{2}}e^{\iota(n_{1}\Delta \kappa_{1} m_{1} \Delta x_{1} + n_{2} \Delta \kappa_{2} m_{2} \Delta x_{2})}\\
    & + \overline{B}_{n_{1}n_{2}} e^{\iota(n_{1} \Delta \kappa_{1} m_{1} \Delta x_{1} - n_{2} \Delta \kappa_{2} m_{2} \Delta x_{2})} \Big] \\
\end{aligned}
\end{equation}
where
\begin{equation}
\begin{aligned}
	B_{n_{1}n_{2}} &= \sum_{n_{2}=0}^{N_{2}}\sum_{n_{1}=0}^{N_1}\sqrt{S(n_{1} \Delta \kappa_{1}, n_{2} \Delta \kappa_{2})\Delta\kappa_{1}\Delta\kappa_{2}}\\
    & \Big[ \sqrt{(1 - \sum_{i_{1} + j_{1} = n_{1}}^{i_{1} \geq j_{1} \geq 0}\sum_{i_{2} + j_{2} = n_{2}}^{n_{2} \geq |i_{2}| \geq |j_{2}| \geq 0}b_{p}^{2}(i_{1} \Delta \kappa_{1}, j_{1} \Delta \kappa_{1}, i_{2} \Delta \kappa_{2}, j_{2} \Delta \kappa_{2})}e^{\iota\Phi_{n_{1}n_{2}}^{(1)}}\\
    & + [\sum_{i_{1} + j_{1} = n_{1}}^{i_{1} \geq j_{1} \geq 0}\sum_{i_{2} + j_{2} = n_{2}}^{n_{2} \geq |i_{2}| \geq |j_{2}| \geq 0} b_{p}(i_{1} \Delta \kappa_{1}, j_{1} \Delta \kappa_{1}, i_{2} \Delta \kappa_{2}, j_{2} \Delta \kappa_{2}) e^{\iota(\Phi_{i_{1}i_{2}}^{(1|2)} + \Phi_{j_{1}j_{2}}^{(1|2)} +\beta(i_{1} \Delta \kappa_{1}, j_{1} \Delta \kappa_{1}, i_{2} \Delta \kappa_{2}, j_{2} \Delta \kappa_{2}))}] \Big]\\
\end{aligned}
\end{equation}
and
\begin{equation}
\begin{aligned}
	\overline{B}_{n_{1}n_{2}} &= \sum_{n_{2}=0}^{N_{2}}\sum_{n_{1}=0}^{N_1}\sqrt{S(n_{1} \Delta \kappa_{1}, -n_{2} \Delta \kappa_{2})\Delta\kappa_{1}\Delta\kappa_{2}}\\
    & \Big[ \sqrt{(1 - \sum_{i_{1} + j_{1} = n_{1}}^{i_{1} \geq j_{1} \geq 0}\sum_{i_{2} + j_{2} = -n_{2}}^{|n_{2}| \geq |i_{2}| \geq |j_{2}| \geq 0}b_{p}^{2}(i_{1} \Delta \kappa_{1}, j_{1} \Delta \kappa_{1}, i_{2} \Delta \kappa_{2}, j_{2} \Delta \kappa_{2})}e^{\iota\Phi_{n_{1}n_{2}}^{(2)}}\\
    & + [\sum_{i_{1} + j_{1} = n_{1}}^{i_{1} \geq j_{1} \geq 0}\sum_{i_{2} + j_{2} = -n_{2}}^{|n_{2}| \geq |i_{2}| \geq |j_{2}| \geq 0} b_{p}(i_{1} \Delta \kappa_{1}, j_{1} \Delta \kappa_{1}, i_{2} \Delta \kappa_{2}, j_{2} \Delta \kappa_{2}) e^{\iota(\Phi_{i_{1}i_{2}}^{(1|2)} + \Phi_{j_{1}j_{2}}^{(1|2)} + \beta(i_{1} \Delta \kappa_{1}, j_{1} \Delta \kappa_{1}, i_{2} \Delta \kappa_{2}, j_{2} \Delta \kappa_{2}))}] \Big]\\
\end{aligned}
\end{equation}
which, following Eq.\ \eqref{eqn:FFT_standard}, equates to a sequence of FFTs. Again, let $\bm{A}$ denote the matrix form of the random field and let $\bm{B}$ and $\overline{\bm{B}}$ denote the matrix forms of $B_{n_1n_2}$ and $\overline{B}_{n_1n_2}$. We can express the simulation compactly in terms of FFTs as 
\begin{equation}
    \bm{A} = 2 \left[\Re\{N^{2} (\text{IFFT}_{\kappa_2} \circ \text{IFFT}_{\kappa_1}(\bm{B})) + N (\text{FFT}_{\kappa_2} \circ \text{IFFT}_{\kappa_1}(\overline{\bm{B}}))\}\right]
\end{equation}
where the subscript $\kappa_1$ or $\kappa_2$ specifies the axis of the matrix over which the FFT/IFFT operates.

In the case of quadrant random fields, this further simplifies as 
\begin{equation}
\begin{aligned}
    A(m_{1} \Delta x_{1}, m_{2} \Delta x_{2}) &=  \sum_{n_{2}=0}^{N_{2}}\sum_{n_{1}=0}^{N_1} 2 B_{n_{1}n_{2}} \Big[e^{\iota(n_{1}\Delta \kappa_{1} m_{1} \Delta x_{1})}(e^{\iota(n_{2}\Delta \kappa_{2} m_{2} \Delta x_{2})} + e^{\iota(-n_{2}\Delta \kappa_{2}m_{2}\Delta x_{2})}) \Big] \\
\end{aligned}
\end{equation}
which can again be expressed compactly as
\begin{equation}
    \bm{A} = 2 \left[\Re\{N^{2}(\text{IFFT}_{\kappa_2} \circ \text{IFFT}_{\kappa_1}(\bm{B}))+N (\text{FFT}_{\kappa_2} \circ \text{IFFT}_{\kappa_1}(\bm{B}))\}\right]
\end{equation}

Hence, for quadrant fields, it is necessary to perform only a single IFFT over $\kappa_1$ and to perform both an FFT and an IFFT along $\kappa_2$ on the resultant.

\subsection{FFT implementation for the simulation of d-dimensional random fields}

The simulation formula for the simulation of $d$-dimensional random fields is a direct extension of the 2-dimensional case and is given as follows:
\begin{equation}
\begin{aligned}
    & A(m_{1} \Delta x_{1}, m_{2} \Delta x_{2}, \dots, m_{d} \Delta x_{d}) =\\
    &  \sum_{n_{d}=-N_{d}}^{N_{d}} \dots \sum_{n_{2}=-N_{2}}^{N_{2}}\sum_{n_{1}=0}^{N_1}[\sqrt{2}\sqrt{2S_{p}(n_{1} \Delta \kappa_{1}, n_{2} \Delta \kappa_{2} \dots n_{d} \Delta \kappa_{d})\Delta\kappa_{1}\Delta\kappa_{2} \dots \Delta\kappa_{d}}\\
    & \cos(n_{1}m_{1}\Delta\kappa_{1}\Delta x_{1} + n_{2}m_{2}\Delta\kappa_{2}\Delta x_{2} \dots n_{d}m_{d}\Delta\kappa_{d}\Delta x_{d} + \Phi_{n_{1}n_{2}\dots n_{d}})\\
    & + [\sum_{i_{1} + j_{1} = n_{1}}^{i_{1} \geq j_{1} \geq 0}\sum_{i_{2} + j_{2} = n_{2}}^{|n_{2}| \geq |i_{2}| \geq |j_{2}| \geq 0} \dots \sum_{i_{d} + j_{d} = n_{d}}^{|n_{d}| \geq |i_{d}| \geq |j_{d}| \geq 0}\sqrt{2}\sqrt{2S(n_{1} \Delta \kappa_{1}, n_{2} \Delta \kappa_{2} \dots n_{d} \Delta \kappa_{d})\Delta\kappa_{1}\Delta\kappa_{2}\dots\Delta\kappa_{d}}\\
    & b_{p}(i_{1} \Delta \kappa_{1}, j_{1} \Delta \kappa_{1}, i_{2} \Delta \kappa_{2}, j_{2} \Delta \kappa_{2} \dots \dots i_{d} \Delta \kappa_{d}, j_{d} \Delta \kappa_{d})\\
    &\cos(n_{1}m_{1}\Delta\kappa_{1}\Delta x_{1} + n_{2}m_{2}\Delta\kappa_{2}\Delta x_{2} \dots n_{d}m_{d}\Delta\kappa_{d}\Delta x_{d} + \Phi_{i_{1}i_{2} \dots i_{d}} + \Phi_{j_{1}j_{2} \dots j_{d}} + \\
    & \beta(i_{1} \Delta \kappa_{1}, j_{1} \Delta \kappa_{1}, i_{2} \Delta \kappa_{2}, j_{2} \Delta \kappa_{2}, \dots, \dots i_{d} \Delta \kappa_{d}, j_{d} \Delta \kappa_{d}))]\\
\end{aligned}
\end{equation}
where
\begin{equation}
\begin{aligned}
    & S_{p}(n_{1} \Delta \kappa_{1}, n_{2} \Delta \kappa_{2}, \dots, n_{d} \Delta \kappa_{d}) = S(n_{1} \Delta \kappa_{1}, n_{2} \Delta \kappa_{2}, \dots, n_{d} \Delta \kappa_{d})\\
    & \Big(1 - \sum_{i_{1} + j_{1} = n_{1}}^{i_{1} \geq j_{1} \geq 0} \sum_{i_{2} + j_{2} = n_{2}}^{|n_{2}| \geq |i_{2}| \geq |j_{2}| \geq 0} \dots \sum_{i_{d} + j_{d} = n_{d}}^{|n_{d}| \geq |i_{d}| \geq |j_{d}| \geq 0} b_{p}^{2}(i_{1} \Delta \kappa_{1}, j_{1} \Delta \kappa_{1}, i_{2} \Delta \kappa_{2}, j_{2} \Delta \kappa_{2}, \dots, i_{d} \Delta \kappa_{d}, j_{d} \Delta \kappa_{d} )\Big)\\
	& b_{p}^{2}(i_{1} \Delta \kappa_{1}, j_{1} \Delta \kappa_{1}, i_{2} \Delta \kappa_{2}, j_{2} \Delta \kappa_{2}, \dots, i_{d} \Delta \kappa_{d}, j_{d} \Delta \kappa_{d}) \\
	& = \frac{|B(i_{1} \Delta \kappa_{1}, j_{1} \Delta \kappa_{1}, i_{2} \Delta \kappa_{2}, j_{2} \Delta \kappa_{2}, \dots, i_{d} \Delta \kappa_{d}, j_{d} \Delta \kappa_{d})|^{2}\Delta\kappa_{1}\Delta\kappa_{2}\dots\Delta\kappa_{d}}{S_{p}(i_{1} \Delta \kappa_{1}, i_{2} \Delta \kappa_{2}, \dots, i_{d} \Delta \kappa_{d})S_{p}(j_{1} \Delta\kappa_{1}, j_{2} \Delta \kappa_{2}, \dots, j_{d} \Delta \kappa_{d})S((i_{1} + j_{1}) \Delta \kappa_{1}, \dots, (i_{d} + j_{d}) \Delta \kappa_{d})}\\
\end{aligned}
\end{equation}

Following similar steps involved in the development in the FFT implementation for 1D-1V random field we have
 \begin{equation}
\begin{aligned}
    & A(m_{1} \Delta x_{1}, m_{2} \Delta x_{2}, \dots, m_{d} \Delta x_{d}) =\\
    & 2  \sum_{n_{d}=0}^{N_{d}} \dots \sum_{n_{2}=0}^{N_{2}}\sum_{n_{1}=0}^{N_1}\sum_{I_{1}=1, I_{2}=\pm1,\dots, I_{d}=\pm1}\sqrt{S(I_{1}n_{1} \Delta \kappa_{1}, I_{2}n_{2} \Delta \kappa_{2}, \dots, I_{d}n_{d} \Delta \kappa_{d})\Delta\kappa_{1}\Delta\kappa_{2}\dots\Delta\kappa_{d}}\\
    & \Big[ \sqrt{(1 - \sum_{i_{1} + j_{1} = I_{1}n_{1}}^{i_{1} \geq j_{1} \geq 0}\sum_{i_{2} + j_{2} = I_{2}n_{2}}^{n_{2} \geq |i_{2}| \geq |j_{2}| \geq 0} \dots \sum_{i_{d} + j_{d} = I_{d}n_{d}}^{n_{d} \geq |i_{d}| \geq |j_{d}| \geq 0} b_{p}^{2}(i_{1} \Delta \kappa_{1},j_{1} \Delta \kappa_{1}, i_{2} \Delta \kappa_{2}, j_{2} \Delta \kappa_{2}, \dots, \dots i_{d} \Delta \kappa_{d}, j_{d} \Delta \kappa_{d})}\\
    & e^{\iota\Phi_{n_{1}n_{2} \dots n_{d}}^{I_{1}I_{2} \dots I_{d}}} + [\sum_{i_{1} + j_{1} = I_{1}n_{1}}^{i_{1} \geq j_{1} \geq 0}\sum_{i_{2} + j_{2} = I_{2}n_{2}}^{n_{2} \geq |i_{2}| \geq |j_{2}| \geq 0} \dots \sum_{i_{d} + j_{d} = I_{d}n_{d}}^{n_{d} \geq |i_{d}| \geq |j_{d}| \geq 0} b_{p}(i_{1} \Delta \kappa_{1},j_{1} \Delta \kappa_{1}, i_{2} \Delta \kappa_{2}, j_{2} \Delta \kappa_{2}, \dots, \dots i_{d} \Delta \kappa_{d}, j_{d} \Delta \kappa_{d})\\
    & e^{\iota(\Phi_{i_{1}i_{2} \dots i_{d}}^{I_{1}I_{2} \dots I_{d}} + \Phi_{j_{1}j_{2} \dots j_{d}}^{I_{1}I_{2} \dots I_{d}} +\beta(i_{1} \Delta \kappa_{1},j_{1} \Delta \kappa_{1}, i_{2} \Delta \kappa_{2}, j_{2} \Delta \kappa_{2}, \dots, \dots i_{d} \Delta \kappa_{d}, j_{d} \Delta \kappa_{d}))}] \Big]\\
    & e^{\iota(I_{1}n_{1} m_{1} \Delta \kappa_{1} \Delta x_{1} + I_{2}n_{2} m_{2} \Delta \kappa_{2} \Delta x_{2} + \dots  + I_{d}n_{d} m_{d} \Delta \kappa_{d} \Delta x_{d})}\\
\end{aligned}
\end{equation}
This can be simplified to a form amenable to the FFT implementation as
\begin{equation}
\begin{aligned}
    A(m_{1} \Delta x_{1}, m_{2} \Delta x_{2}, \dots, m_{d} \Delta x_{d}) = & 2  \sum_{n_{d}=0}^{N_{d}} \dots \sum_{n_{2}=0}^{N_{2}} \sum_{n_{1}=0}^{N_1} \sum_{I_{1}=1, I_{2}=\pm1,\dots, I_{d}=\pm1} \\
    & \Big[ B_{n_{1}n_{2} \dots n_{d}}^{{I_{1}I_{2} \dots I_{d}}}e^{\iota(I_{1}n_{1}m_{1}\Delta \kappa_{1}  \Delta x_{1} + I_{2}n_{2}m_{2}\Delta \kappa_{2}  \Delta x_{2} + \dots  + I_{d}n_{d}m_{d}\Delta \kappa_{d}  \Delta x_{d})} \Big] \\
    \label{eqn:nD_FFT}
\end{aligned}
\end{equation}
Again, expressing this in terms of FFT and IFFT operations the following results:
\begin{equation}
\begin{aligned}
    \bm{A} = 2\left[ \sum_{I_{1}=1, I_{2}=\pm1,\dots, I_{d}=\pm1} \Re \{N^{J}\text{FFT}_{\kappa_d}^{I_d} \circ \text{FFT}_{\kappa_{d-1}}^{I_{d-1}} \circ \dots \circ \text{FFT}_{\kappa_1}^{I_1}(\bm{B}^{I_1I_2\dots I_d})\} \right]
\end{aligned}
\end{equation}
where 
\begin{equation}
\begin{aligned}
    J = \sum_{j=1}^d \hat{I}_j, \quad & \hat{I}_j=1 \text{ if } I_j=1,\\ &\hat{I}_j=0 \text{ otherwise}
\end{aligned}
\end{equation}
$\text{FFT}^{I_j}$ equals IFFT if $I_j=1$ and FFT if $I_j=-1$, and $\bm{B}^{I_1I_2\dots I_d}$ are the $d^{th}$-order tensors having components $B_{n_{1}n_{2} \dots n_{d}}^{{I_{1}I_{2} \dots I_{d}}}$  in Eq.\ \eqref{eqn:nD_FFT}. For example, the 3-dimensional implementation takes the following form:
\begin{equation}
\begin{aligned}
    \bm{A} =& \Re \{N^3 \text{IFFT}_{\kappa_3} \circ \text{IFFT}_{\kappa_2} \circ \text{IFFT}_{\kappa_1}(\bm{B}^{111})\\ 
    & + N^2 \text{FFT}_{\kappa_3} \circ \text{IFFT}_{\kappa_2} \circ \text{IFFT}_{\kappa_1}(\bm{B}^{-111})\\
    & + N^2 \text{IFFT}_{\kappa_3} \circ \text{FFT}_{\kappa_2} \circ \text{IFFT}_{\kappa_1}(\bm{B}^{1-11})\\
    & + N \text{FFT}_{\kappa_3} \circ \text{FFT}_{\kappa_2} \circ \text{IFFT}_{\kappa_1}(\bm{B}^{-1-11})\}
\end{aligned}
\end{equation}

In the case of quadrant random fields, the FFT implementation can be further simplified to 
\begin{equation}
\begin{aligned}
    A(m_{1} \Delta x_{1}, m_{2} \Delta x_{2}, \dots, m_{d} \Delta x_{d}) &=  2 \sum_{n_{d}=0}^{N_{d}} \dots \sum_{n_{2}=0}^{N_{2}}\sum_{n_{1}=0}^{N_1}\\
    & B_{n_{1}n_{2} \dots n_{d}} \Big[e^{\iota(n_{1}m_{1}\Delta \kappa_{1} \Delta x_{1})}(e^{\iota(n_{2}m_{2} \Delta \kappa_{2} \Delta x_{2})} + e^{-\iota(n_{2}m_{2} \Delta \kappa_{2} \Delta x_{2})}) \\
    & \dots (e^{\iota(n_{d}m_{d} \Delta \kappa_{d} \Delta x_{d})} + e^{-\iota(n_{d}m_{d} \Delta \kappa_{d} \Delta x_{d})}) \Big] \\
\end{aligned}
\label{eqn:nd_quadrant}
\end{equation}
In terms of FFT and IFFT operators, it takes the following form:
\begin{equation}
\begin{aligned}
    \bm{A} = 2\left[ \sum_{I_{1}=1, I_{2}=\pm1,\dots, I_{d}=\pm1} \Re \{N^{J}\text{FFT}_{\kappa_d}^{I_d} \circ \text{FFT}_{\kappa_{d-1}}^{I_{d-1}} \circ \dots \circ \text{FFT}_{\kappa_1}^{I_1}(\bm{B}\} \right]
\end{aligned}
\end{equation}
where $\bm{B}$ is the $d^{th}$-order tensor having terms $B_{n_{1}n_{2} \dots n_{d}}$ in Eq.\ \eqref{eqn:nd_quadrant}.




\subsection{Notes on the use of the FFT technique}

It is well known that the application of the FFT technique requires that certain conditions be satisfied. One such condition relates the spatial and wave number discretizations as follows:
\begin{equation}
\begin{aligned}
	& \Delta \kappa_{1} \Delta x_{1} = \frac{2\pi}{N_{1}}\\
	& \Delta \kappa_{2} \Delta x_{2} = \frac{2\pi}{N_{2}}\\
	& \vdots\\
	& \Delta \kappa_{d} \Delta x_{d} = \frac{2\pi}{N_{d}}\\
\end{aligned}
\label{eqn:FFT_conditions}
\end{equation}
This means is equivalent to a condition that the spatial domain over which the samples are generated is always equal to one period.

The general procedure for simulating $d$-dimensional third-order random fields with the FFT implementation is as follows:
\begin{enumerate}
\item Assign the appropriate the wave number discretisation over the $d$ dimensions of the power spectrum and the bispectrum. The associated spatial increments follow from Eq.\ \eqref{eqn:FFT_conditions}.
\item Ensure that the spatial increments satisfy the conditions in Eq.\ \eqref{eqn:aliasing} to avoid aliasing.
\item Generate the necessary $2^{d-1}$, $d^{th}$-order random phase tensors $\boldsymbol{\Phi}^{I_1I_2\dots I_d}$ for general fields or a single $d^{th}$-order random phase tensor for quadrant random fields.
\item Compute the $2^{d-1}$, $d^{th}$-order spectral tensors $\bm{B}^{I_1I_2\dots I_d}$ for general fields or a single $d^{th}$-order spectral tensor for quandrant random fields.
\item Apply FFT and IFFT appropriately along the different axes of the d-dimensional spectral tensor(s) $\bm{B}$ according to the equations above.
\end{enumerate}

The major advantage of the FFT implementation is the computational expense. Each FFT has well-known complexity of the order $O(M \log N)$, whereas each summation of cosines has complexity of the order $O(MN)$. Because the summations in the original formulation are nested over each dimension, the complexity grows exponentially with dimension as $O((MN)^d)$. However, as we can see from the above expressions, the FFT implementation requires only $2^{d-1}d$ FFTs and therefore has complexity of order $O(d2^{d-1}M \log N) \ll O((MN)^d)$. For quadrant random fields, this is reduced even further having order $O(2dM\log N)$ and therefore only scales linearly with dimension. 

The result is a drastic reduction in the computational time, without which the simulation of multidimensional higher-order random fields becomes infeasible. To illustrate the savings, Table \ref{table:fft_1d_comparision} shows a comparison of the computation time for the non-FFT and the FFT implementations for a 1-dimensional random field for different number of sample functions generated. 
\begin{table}[!ht]
\centering
\begin{tabular}{l l l}
\hline
 & \multicolumn{2}{c}{\textbf{Time (sec)}}\\
\hline
\textbf{No. of Samples} & \textbf{Non-FFT} & \textbf{FFT}\\
\hline
128 & 14.842 & 0.0893\\
512 & 26.891 & 0.0957\\
1024 & 48.383 & 0.1399\\
2048 & 82.525 & 0.3750\\
4096 & 456.100 & 1.9270\\
\hline
\end{tabular}
\caption{Comparison of the computation time for simulation of 1D third-order random fields using the standard and FFT implementations.}
\label{table:fft_1d_comparision}
\end{table}
On average the FFT calculations are three orders of magnitude faster.

While Table \ref{table:fft_1d_comparision} illustrates the huge savings for one-dimensional fields, it is particular interest to observe how these computation times scale with the dimension. Table \ref{table:dimensionality_comparision} shows computation times for 2- and 3-dimensional random fields using the FFT implementation remain modest.
\begin{table}[!ht]
\centering
\begin{tabular}{l l l}
\hline
 & \multicolumn{2}{c}{\textbf{Time (sec)}}\\
\hline
\textbf{No. of Samples} & \textbf{2D} & \textbf{3D}\\
\hline
1 & 0.224 & 20.651\\
16 & 0.225 & 21.839\\
128 & 0.274 & 25.600\\
512 & 0.375 & 37.89\\
\hline
\end{tabular}
\caption{Computational time for the simulation of 2D and 3D third-order random fields using the FFT implementation. Standard implementation is not shown because the calculations are impractical on a desktop computer.}
\label{table:dimensionality_comparision}
\end{table}
Note, however, that we do not compare with the summation of cosines here because these calculations become intractable for dimensions greater than one.

\section{Numerical examples}
In this section, we present examples of the simulation of skewed 2- and 3-dimensional random fields from prescribed power spectra and bispectra. These examples, although purely mathematical in nature and not corresponding to any physically meaningful random field, have been developed to call attention to specific features of the proposed methodology.

\subsection{Comparison of 2-dimensional $2^{nd}$- and $3^{rd}$-order random fields}

\begin{figure}[!ht]
  \centering\includegraphics[width=0.6\linewidth]{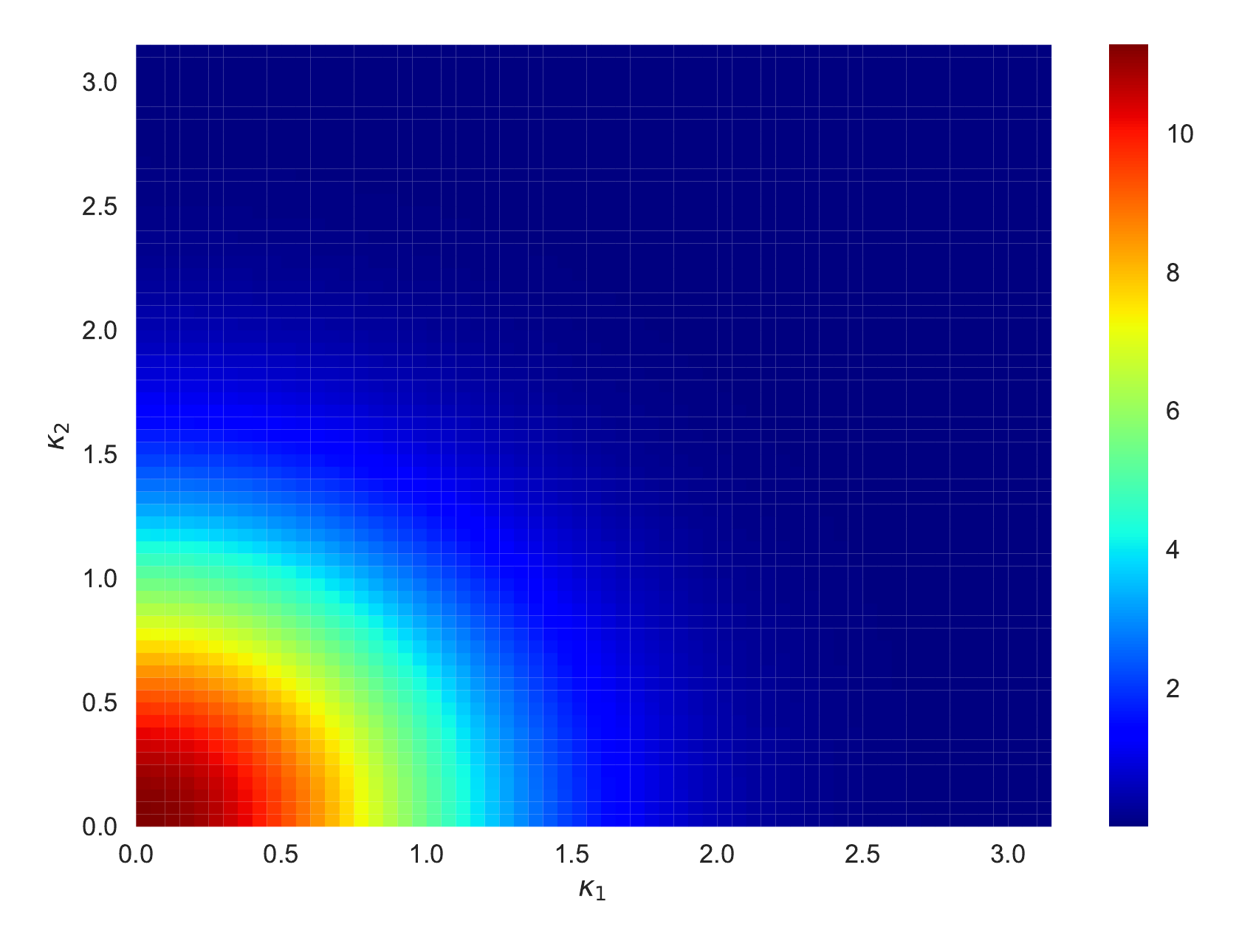}
  \caption{2-dimensional power spectrum}
  \label{fig:example_1_power_spectrum}
\end{figure}

The first example compares the simulation of a 2-dimensional random field by the $2^{nd}$-order
SRM and the $3^{rd}$-order SRM. The prescribed power spectrum is given by
\begin{equation}
    S(\kappa_{1}, \kappa_{2}) = \frac{20}{\sqrt{\pi}}\exp{-\frac{1}{2}(\kappa_{1}^{2} + \kappa_{2}^{2})} \ \text{for} \ \kappa_{1}, \kappa_{2} \geq 0
    \label{eqn:power_spectrum}
\end{equation}
and is plotted in Fig \ref{fig:example_1_power_spectrum}, yielding a random field with zero mean and variance 75. Notice that the power spectrum is radially symmetric.  
The prescribed bispectrum for the 3rd-order random field is given by
\begin{equation}
\begin{aligned}
\Re B(\kappa_{11}, \kappa_{12}, \kappa_{21}, \kappa_{22}) &= \Im  B(\kappa_{11}, \kappa_{12}, \kappa_{21}, \kappa_{22}) \\
&=\frac{58}{\pi}\exp{-(\kappa_{11}^{2} + \kappa_{12}^{2} + \kappa_{21}^{2} + \kappa_{22}^{2})}
& \text{for} \ \kappa_{11},\kappa_{12},\kappa_{21},\kappa_{22}\geq0
\end{aligned}
\end{equation}

Visualisation of the 2-dimensional bispectrum, which is a $4^{th}$-order tensor, is not trivial and is not presented here. Of particular note here is that the bispectrum is symmetric across all dimensions, i.e.\ it has the same rate of decay along each axis. This implies that the coupling of the waves is the same in both dimensions.

\begin{figure}[!ht]
\centering
\begin{subfigure}{.49\textwidth}
  \centering
  \includegraphics[width=\linewidth]{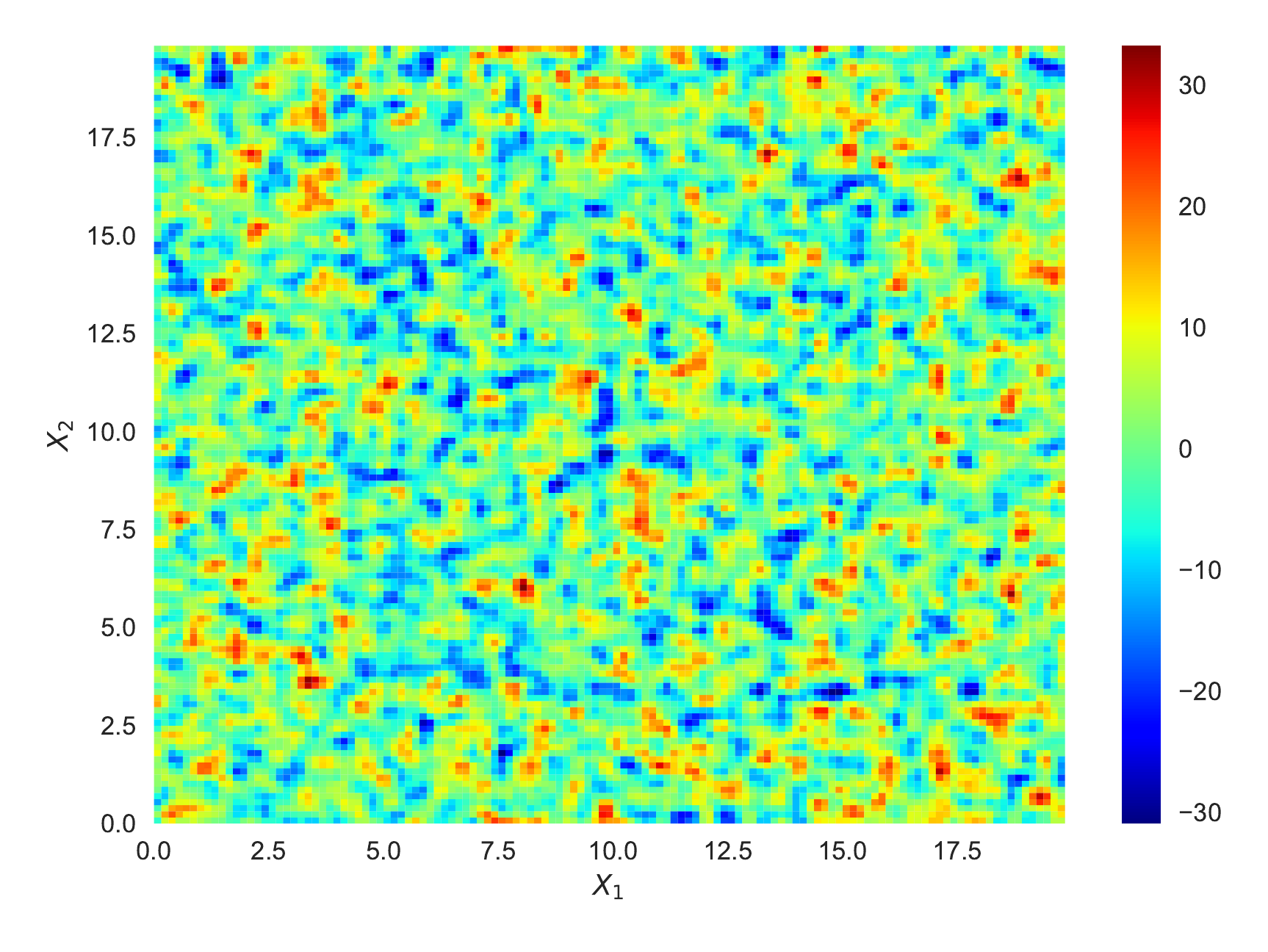}
  \caption{$2^{nd}$-order}
  \label{fig:example_1_SRM_sample}
\end{subfigure}
\begin{subfigure}{.49\textwidth}
  \centering
  \includegraphics[width=\linewidth]{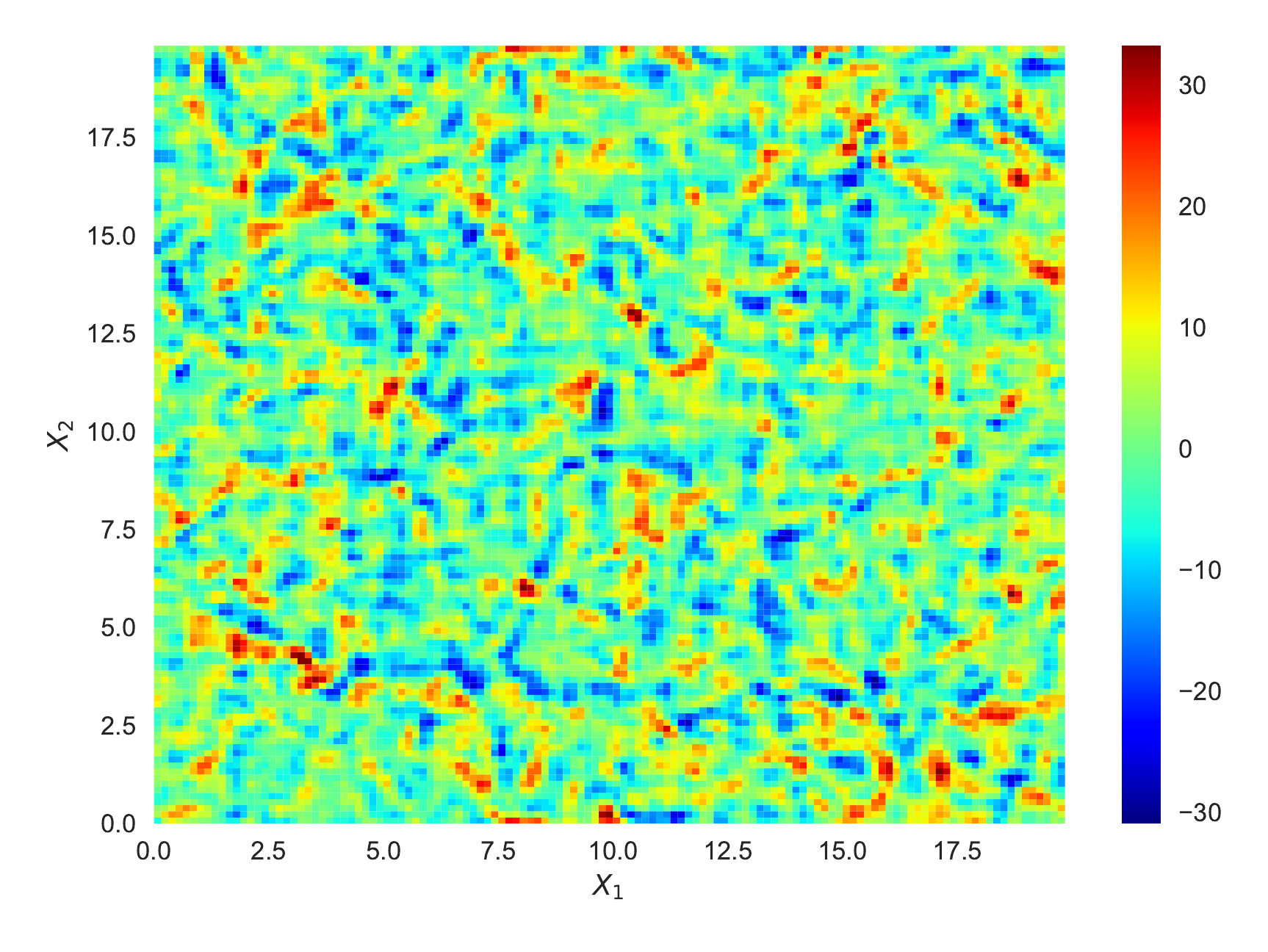}
  \caption{$3^{rd}$-order}
  \label{fig:example_1_BSRM_sample}
\end{subfigure}%
\caption{Example 1 -- $2^{nd}$- and $3^{rd}$-order, 2-dimensional random fields simulated by SRM}
\label{fig:example_1_SR}
\end{figure}

\begin{figure}[!ht]
  \centering\includegraphics[width=0.6\linewidth]{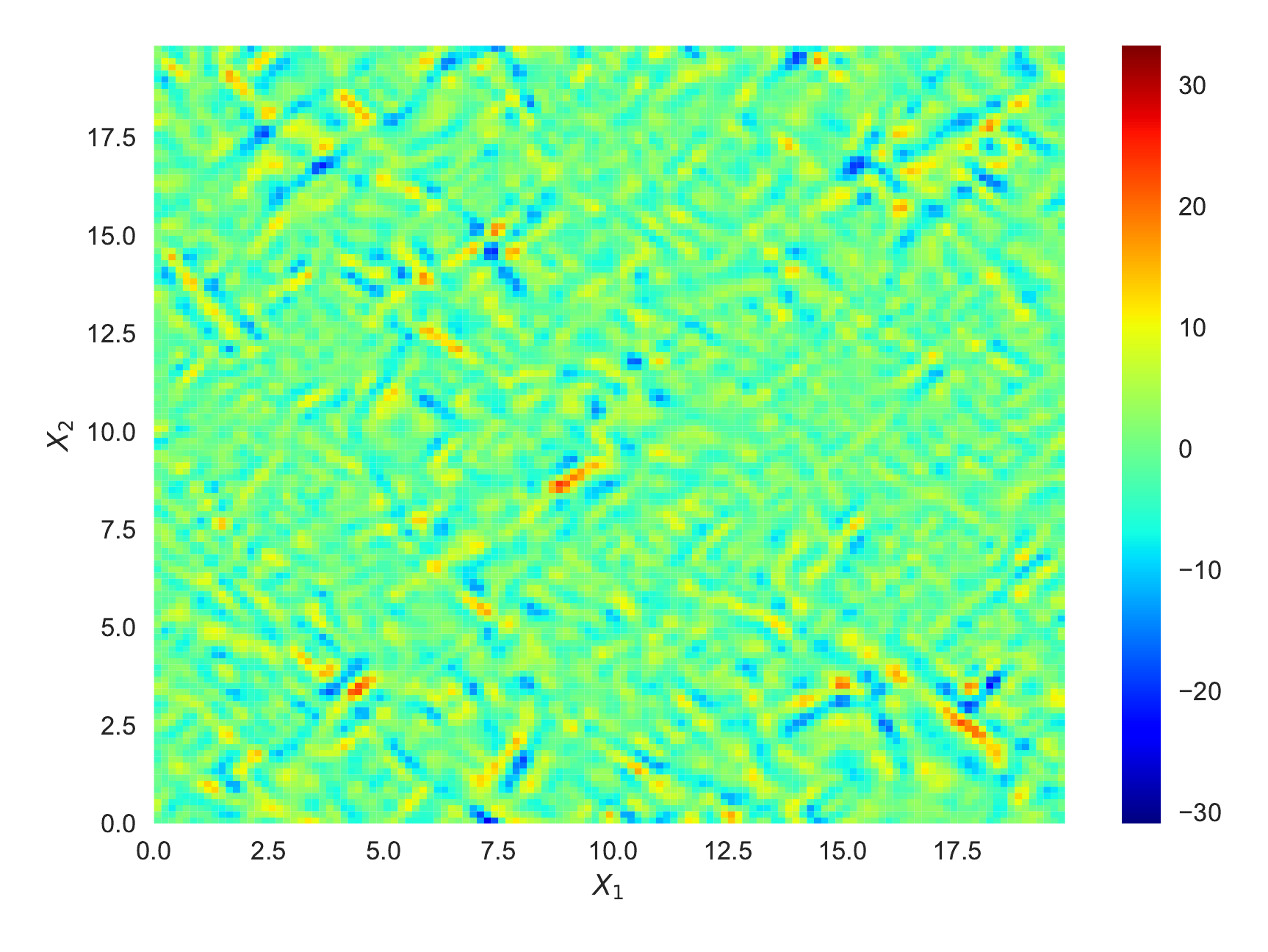}
  \caption{Example 1 -- Difference between sample realizations from the $3^{rd}$-order and $2^{nd}$-order SRM simulations having identical phase angles.}
  \label{fig:example_1_comparision}
\end{figure}

1000 samples of the $2^{nd}$- and $3^{rd}$-order random fields are simulated using the SRM with the following parameters.
\begin{equation}
\begin{aligned}
	&\Delta x_{1} = \Delta x_{2} = 0.78125\\
	&\Delta \kappa_{1} = \Delta \kappa_{2} = 0.0628\\
	& N_{1} = N_{2} = 64\\
	& M_{1} = M_{2} = 128\\
\end{aligned}
\end{equation}
Plots of representative $2^{nd}$- and $3^{rd}$-order sample realizations having identical phase angles are presented in Figure \ref{fig:example_1_SR}.
On initial inspection, both sample realizations look similar. However a closer inspection of the samples and their statistical properties reveals interesting characteristics. The difference between the sample realizations of the $2^{nd}$- and $3^{rd}$-order random fields is shown in Figure \ref{fig:example_1_comparision}. The plot clearly shows that there are significant differences between the two sample realizations. These differences arise from asymmetry introduced by the proposed methodology. Also note that the differences are oriented along a $\arctan(1)=45^{\circ}$ and $\arctan(-1)= -45^{\circ}$ angle relative to the $x_1$ and $x_2$ axes. This arises because the form of the bispectrum is identical in both the $x_1$ and $x_2$ directions. Consequently, the length-scale associated with third-order correlations in the $x_1$ and $x_2$ axes are the same -- resulting in the $45^{\circ}$ and $-45^{\circ}$ ``bands'' of skewness. 

Statistical properties, estimated from the 1000 sample realizations, are presented in Table \ref{table:example_1}, illustrating the ability of the proposed methodology to match the theoretical properties upto to the third-order.
\begin{table}[!ht]
\centering
\begin{tabular}{l l l l}
\hline
\textbf{Moments} & \textbf{Target} & \textbf{$3^{rd}$-order} & \textbf{$2^{nd}$-order}\\
\hline
Mean & 0.00 & 0.0173 & 0.0173 \\
Variance & 74.4874 & 74.4764 & 74.4734\\
Skewness & 0.2022 & 0.2040 & 0.0008\\
\hline
\end{tabular}
\caption{Example 1 -- Target and estimated moments of random fields generated by the 2nd and 3rd order SRM.}
\label{table:example_1}
\end{table}
The original SRM, on the other hand, matches the properties of the process only up to second-order.

\subsection{2-dimensional random fields with different bispectra}

In the second example, we modify the bispectrum such that wave interactions occur on different lenth-scales in the $\kappa_1$ and $\kappa_2$ directions and illustrate how the asymmetric features of random field differ with these relative length-scales. We generate two sets of random fields with the same power spectrum given above in Eq.\ \eqref{eqn:power_spectrum} and shown in Figure \ref{fig:example_1_power_spectrum}. However, we consider two different bispectra as follows below
\begin{equation}
\begin{aligned}
    \Re B_{1}\kappa_{11}, \kappa_{12}, \kappa_{21}, \kappa_{22}) = &\Im  B_{1}(\kappa_{11}, \kappa_{12}, \kappa_{21}, \kappa_{22}) \\
    &= \frac{140}{\pi}\exp{-(\kappa_{11}^{2} + 10\kappa_{12}^{2} + \kappa_{21}^{2} + 10\kappa_{22}^{2})}
\end{aligned}
\end{equation}
\begin{equation}
\begin{aligned}
    \Re B_{2}(\kappa_{11}, \kappa_{12}, \kappa_{21}, \kappa_{22}) = &\Im  B_{2}(\kappa_{11}, \kappa_{12}, \kappa_{21}, \kappa_{22}) \\
    &= \frac{140}{\pi}\exp{-(10\kappa_{11}^{2} + \kappa_{12}^{2} + 10\kappa_{21}^{2} + \kappa_{22}^{2})}
\end{aligned}
\end{equation}
Again, visualisation of the 2-dimensional bispectra is not included.



\begin{figure}[!ht]
\centering
\begin{subfigure}{.5\textwidth}
  \centering
  \includegraphics[width=\linewidth]{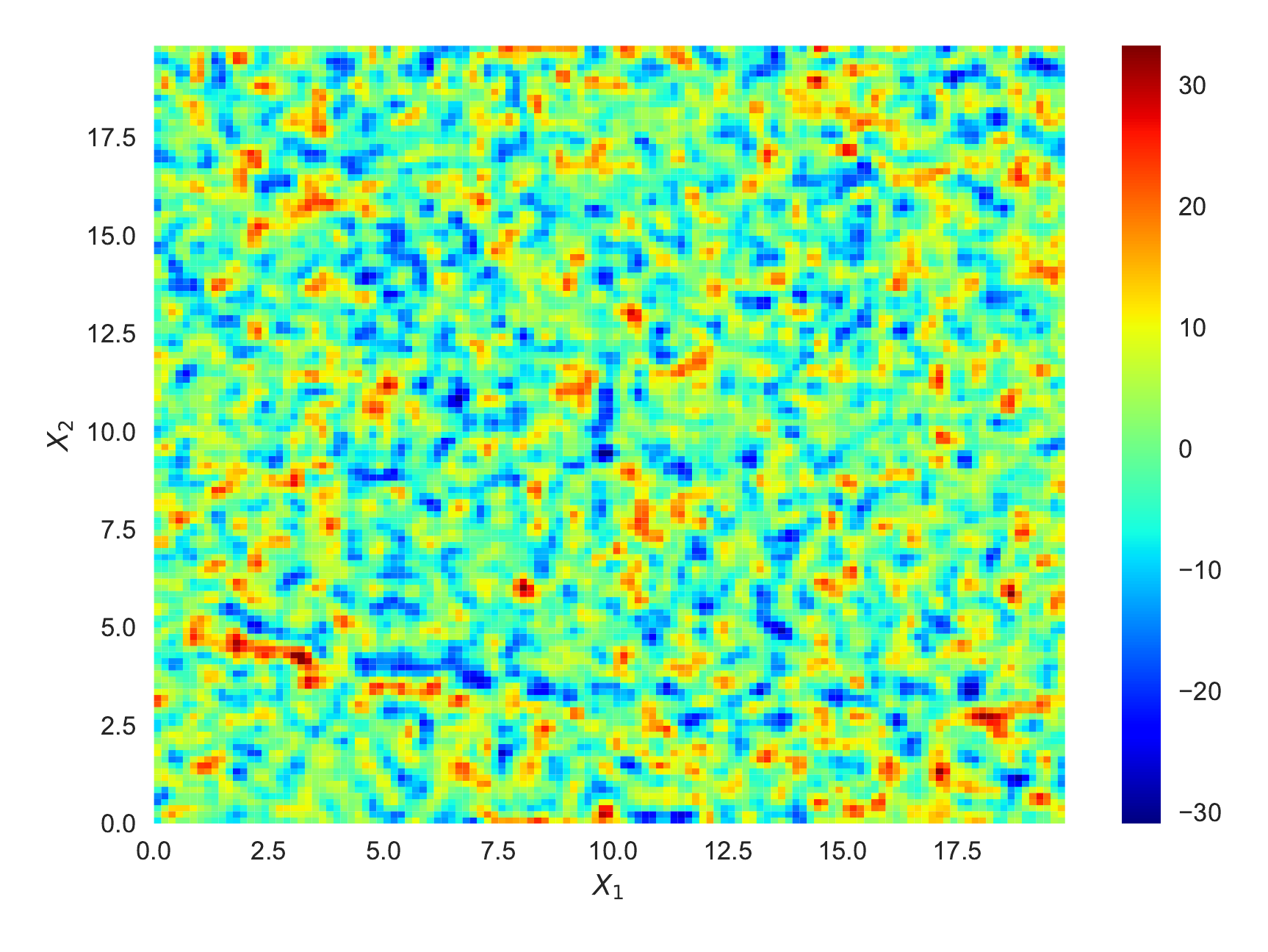}
  \caption{Random field with bispectrum $B_1$.}
  \label{fig:example_2_BSRM_1_sample}
\end{subfigure}%
\begin{subfigure}{.5\textwidth}
  \centering
  \includegraphics[width=\linewidth]{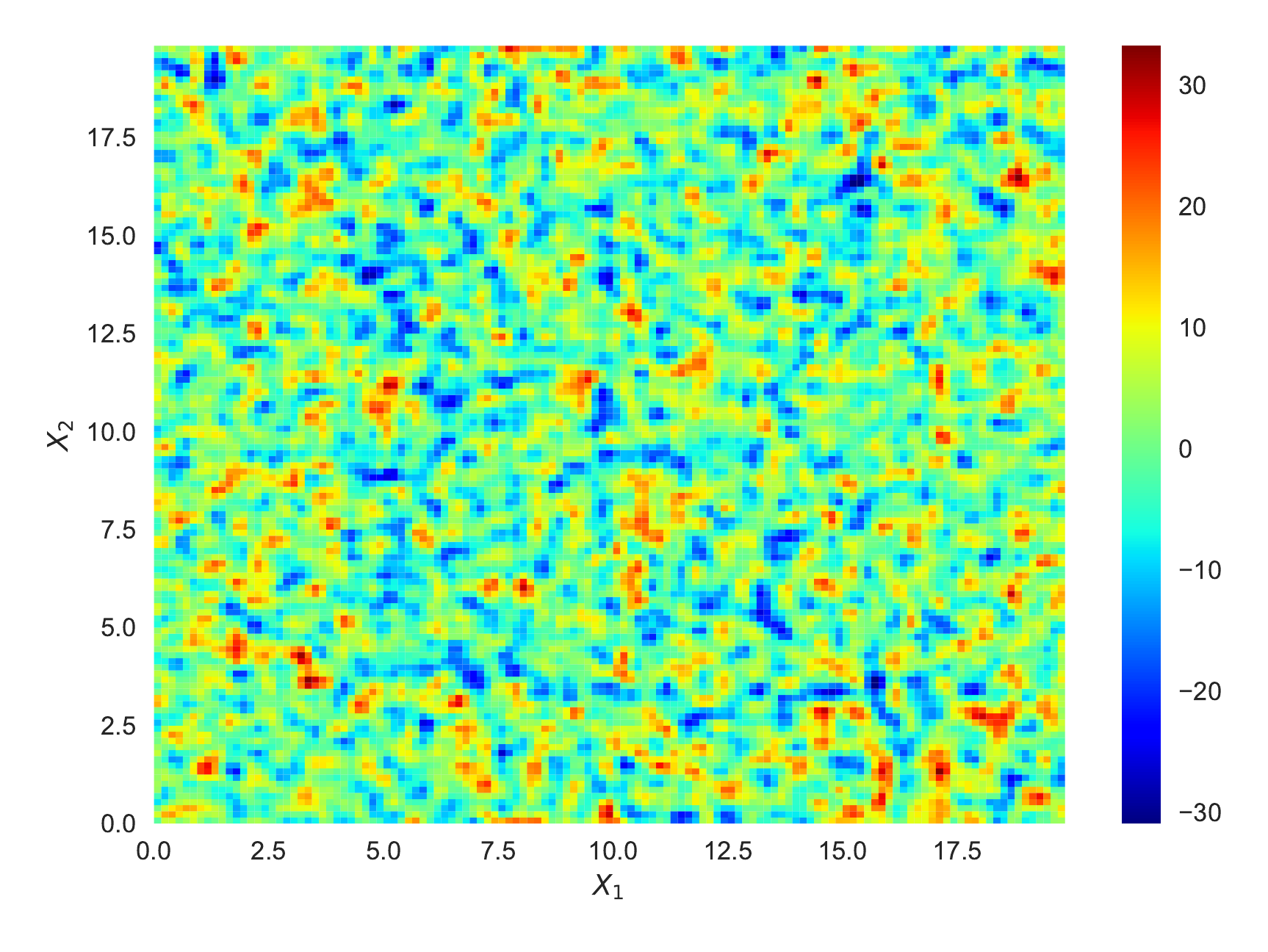}
  \caption{Random field with bispectrum $B_2$.}
  \label{fig:example_2_BSRM_2_sample}
\end{subfigure}
\caption{2-dimensional random fields generated from the two bispectra using the $3^{rd}$-order SRM.}
\label{fig:example_2_SR}
\end{figure}

The first bispectrum shows accelerated decay along the $x_{2}$ axis whereas the second bispectrum has accelerated decay along the $x_{1}$ axis. The samples are again simulated using the FFT implementation of the $3^{rd}$-order SRM. Plots of two sample realisations, again having the same discretization and random phase angles as in example 1 for direct comparison with the $2^{nd}$-order, are presented in Figure \ref{fig:example_2_SR}.


Again, to the naked eye, the sample realisations look similar to the second-order (Figure \ref{fig:example_1_SRM_sample}). But taking difference between the sample realisations of $2^{nd}$- and $3^{rd}$-order fields, we now see that the asymmetric features are elongated along a particular axis. In the case of $B_1$, the asymmetric features lie most prominently along the $x_2$ axis where the decay in bispectrum is more rapid. In other words, since $B_1$ has faster rate of decay along the $\kappa_{2}$ axis, the features are elongated along the $x_{1}$ axis. Interestingly, the asymmetric features occur at an angle approximately $\arctan(\sqrt{10})\approx 73^{\circ}$ and $\arctan(-\sqrt{10})\approx -73^{\circ}$ from the $x_2$ axis (or $\arctan(\sqrt{0.1})\approx 18^{\circ}$ and $\arctan(-\sqrt{0.1})\approx -18^{\circ}$ from the $x_1$ axis). The inverse is true for $B_2$.

\begin{figure}[!ht]
\centering
\begin{subfigure}{.5\textwidth}
  \centering
  \includegraphics[width=\linewidth]{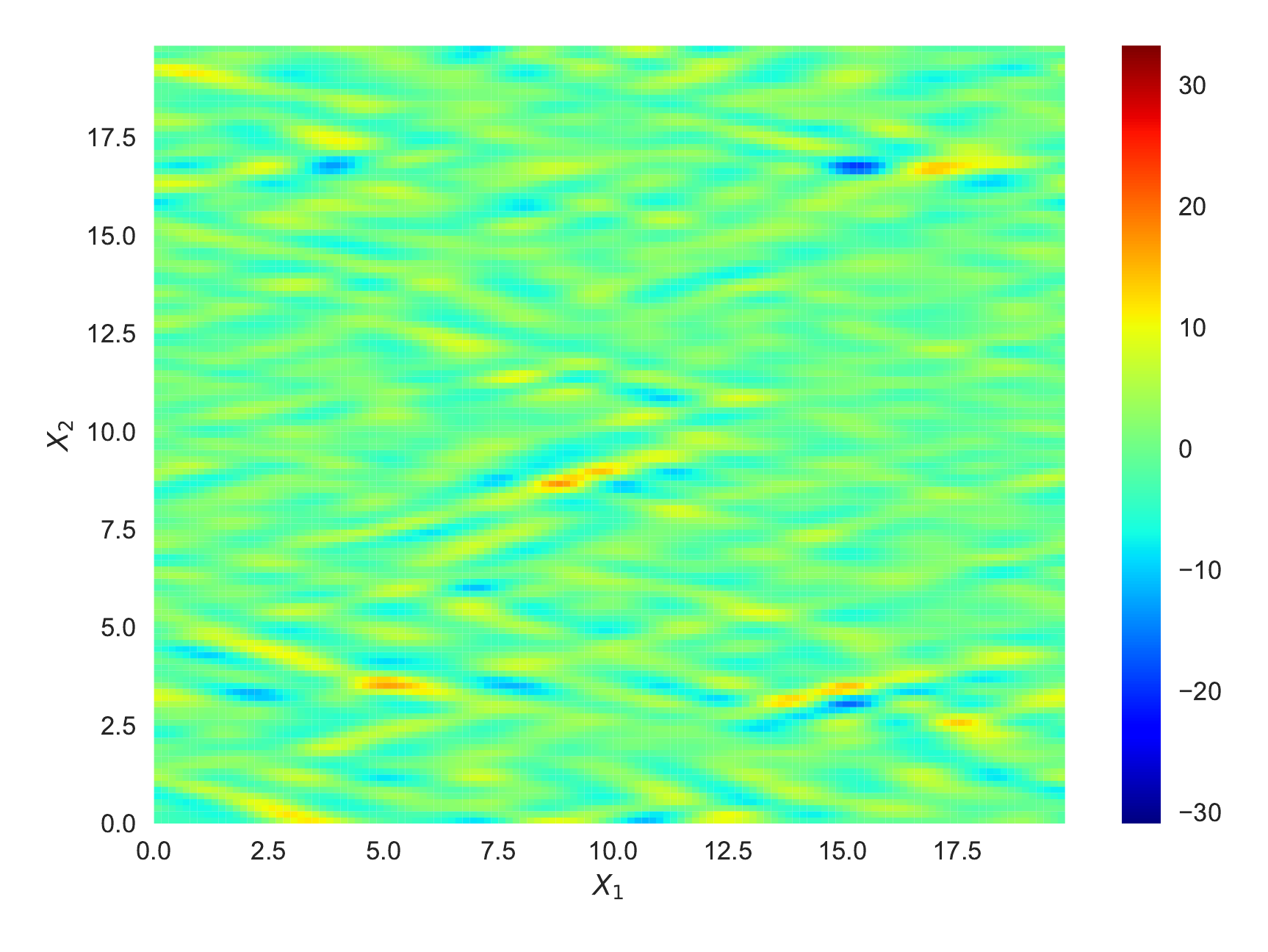}
  \caption{Random field with bispectrum $B_1$.}
  \label{fig:example_2_BSRM_1_SRM}
\end{subfigure}%
\begin{subfigure}{.5\textwidth}
  \centering
  \includegraphics[width=\linewidth]{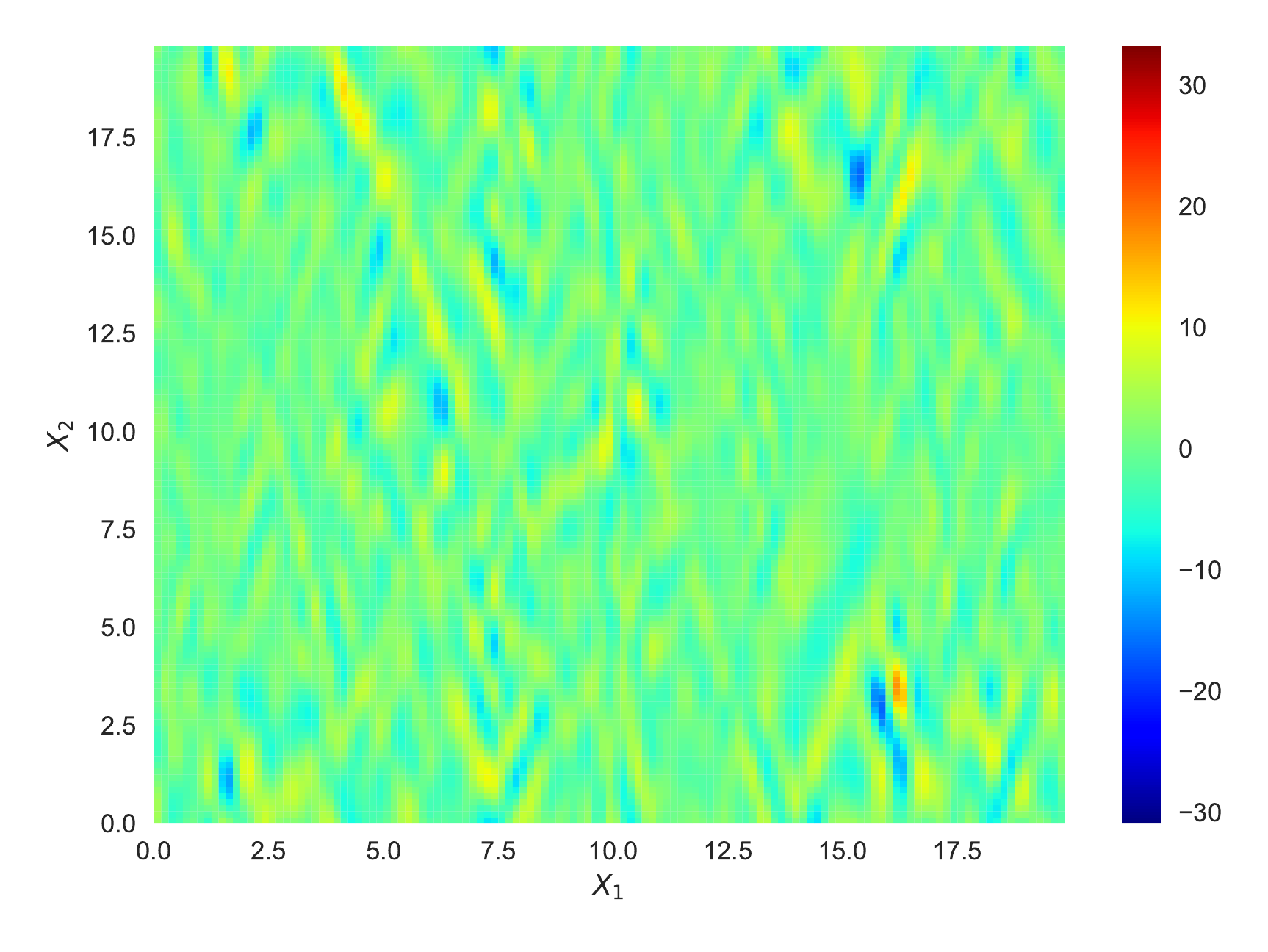}
  \caption{Random field with bispectrum $B_1$.}
  \label{fig:example_2_BSRM_2_SRM}
\end{subfigure}
\caption{Difference between samples generated by BSRM and SRM simulations for both the bispectra}
\label{fig:example_2_comparision}
\end{figure}


Again, we generated 1000 samples of the $2^{nd}$- and $3^{rd}$-order random fields with the following parameters

\begin{equation}
\begin{aligned}
	&\Delta x_{1} = \Delta x_{2} = 0.78125\\
	&\Delta \kappa_{1} = \Delta \kappa_{2} = 0.0628\\
	& N_{1} = N_{2} = 64\\
	& M_{1} = M_{2} = 128\\
\end{aligned}
\end{equation}

and the statistical properties of the sample realizations are presented in Table \ref{table:example_2}.
\begin{table}[!ht]
\centering
\begin{tabular}{l l l l l}
\hline
\textbf{Moments} & \textbf{Target} & \textbf{$3^{rd}$-order, $B_1$} & \textbf{$3^{rd}$-order, $B_2$} & \textbf{$2^{nd}$-order}\\
\hline
Mean & 0.00 & 0.0173 & 0.0173 & 0.0173 \\
Variance & 74.4874 & 74.5158 & 74.4968 & 74.4963\\
Skewness & 0.04559 & 0.04629 & 0.04678 & 0.0008\\
\hline
\end{tabular}
\caption{Target and estimated moments of random fields generated by the 2nd and 3rd order SRM}
\label{table:example_2}
\end{table}
Again, all of the random fields possess the correct mean and variance. However, only the $3^{rd}$-order SRM samples possess the correct skewness. Moreover, they possess the full bispectra but this cannot be visualized.


\subsection{Comparison of 3-dimensional $2^{nd}$- and $3^{rd}$-order random fields}

\begin{figure}[!ht]
\centering
  \includegraphics[width=0.6\linewidth]{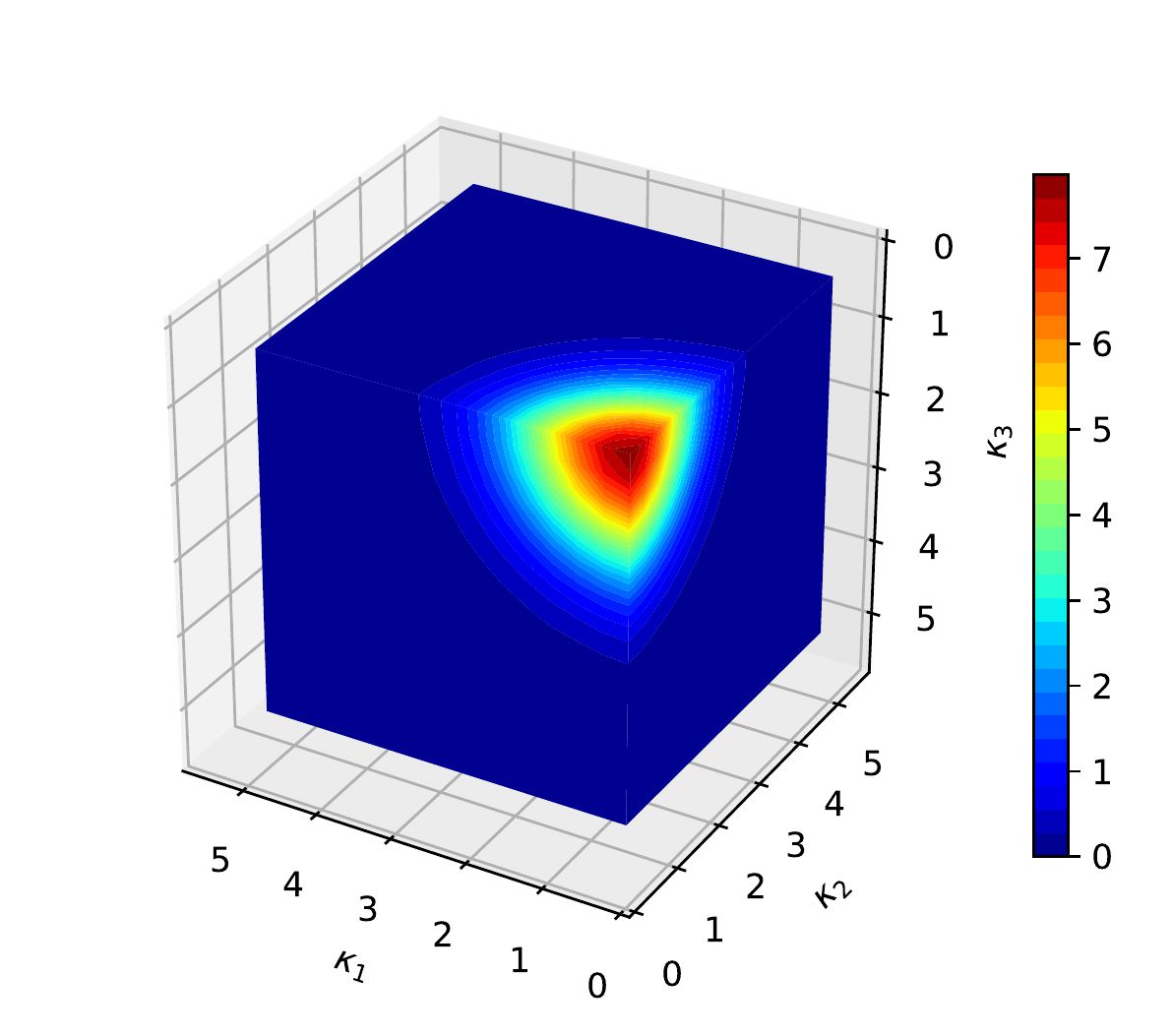}
  \caption{3-dimensional power spectrum}
  \label{fig:example_3_ps}
\end{figure}

In this example, we compare simulations of 3-dimensional random fields having a prescribed power spectrum ($2^{nd}$-order) and power spectrum and bispectrum ($3^{rd}$-order). Both random fields have a power spectrum given by:
\begin{equation}
    S(\omega_{1}, \omega_{2}, \omega_{3}) = \frac{20}{\sqrt{2\pi}}\exp{-\frac{1}{2}(\omega_{1}^{2} + \omega_{2}^{2} + \omega_{3}^{2})}
\end{equation}
and plotted in Figure \ref{fig:example_3_ps}. 

The third-order random field has bispectrum given by
\begin{equation}
\begin{aligned}
    &\Re B(\omega_{11}, \omega_{12}, \omega_{21}, \omega_{22}, \omega_{31}, \omega_{32}) = \Im  B(\omega_{11}, \omega_{12}, \omega_{21}, \omega_{22}, \omega_{31}, \omega_{32}) \\
    &= \frac{22}{2\pi}\exp{-(\omega_{11}^{2} + \omega_{12}^{2} + \omega_{21}^{2} + \omega_{22}^{2} + \omega_{31}^{2} + \omega_{32}^{2})}
\end{aligned}
\end{equation}
Visualisation of this 3-dimensional bispectrum, which is a $6^{th}$-order tensor is not trivial and is therefore not presented here. 

\begin{figure}[!ht]
\centering
\begin{subfigure}{.47\textwidth}
  \centering
  \includegraphics[width=\linewidth]{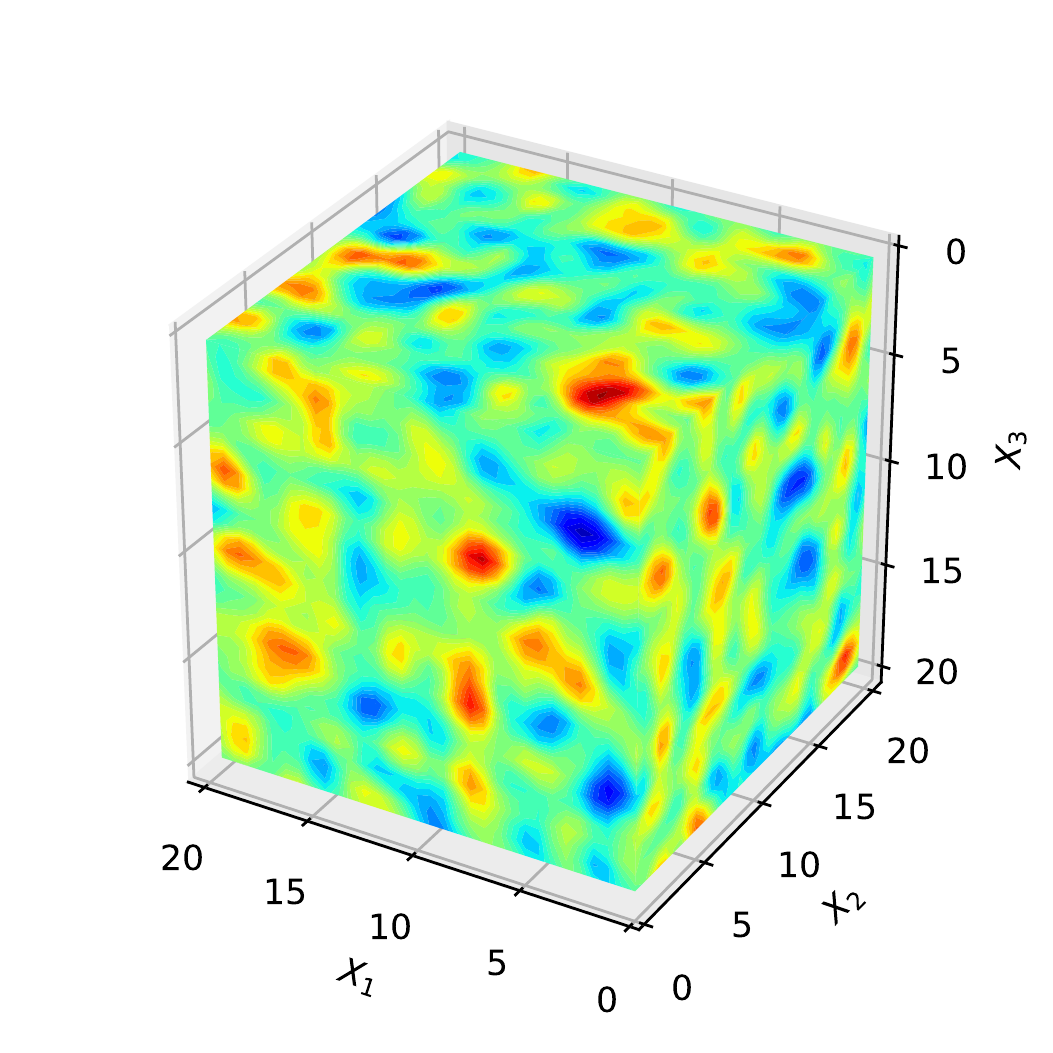}
  \caption{$2^{nd}$-order}
  \label{fig:example_3_SRM_sample}
\end{subfigure}%
\begin{subfigure}{.53\textwidth}
  \centering
  \includegraphics[width=\linewidth]{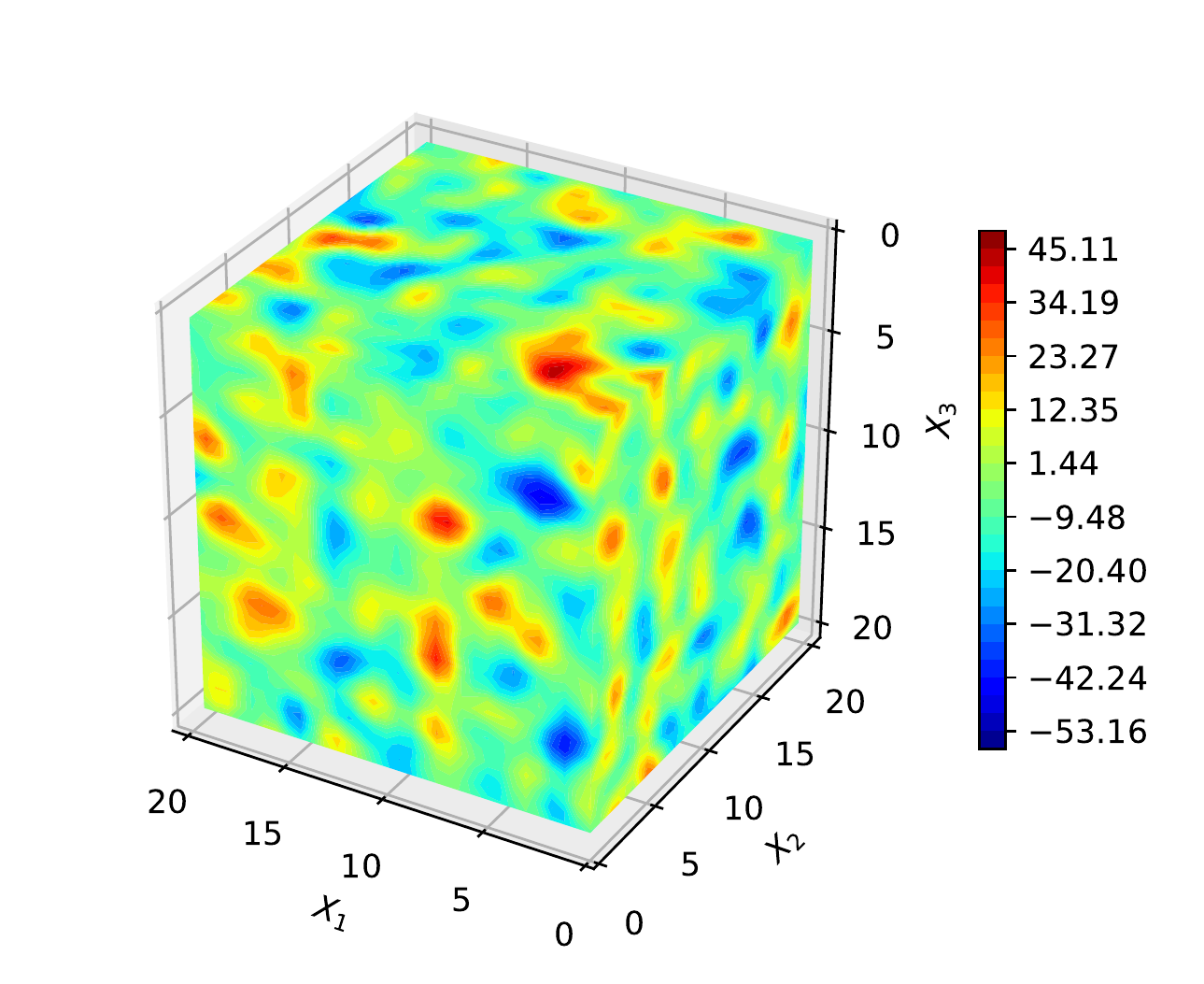}
  \caption{$3^{rd}$-order}
  \label{fig:example_3_BSRM_sample}
\end{subfigure}
\caption{3-dimensional random fields simulated by the $2^{nd}$- and $3^{rd}$-order SRM.}
\label{fig:example_3_SR}
\end{figure}

1000 samples with the following discretization were simulated
\begin{equation}
\begin{aligned}
	&\Delta x_{1} = \Delta x_{2} = \Delta x_{3} = 0.625\\
	&\Delta \kappa_{1} = \Delta \kappa_{2}  = \Delta \kappa_{2} = 0.314\\
	& N_{1} = N_{2} = N_{3} = 16\\
	& M_{1} = M_{2} = M_{3} = 32\\
\end{aligned}
\end{equation}
Plots of representative sample realisationsof the $2^{nd}$- and $3^{rd}$-order random fields, having identical phase angles, are presented in Figure \ref{fig:example_3_SR}.

\begin{figure}[!ht]
  \centering\includegraphics[width=0.6\linewidth]{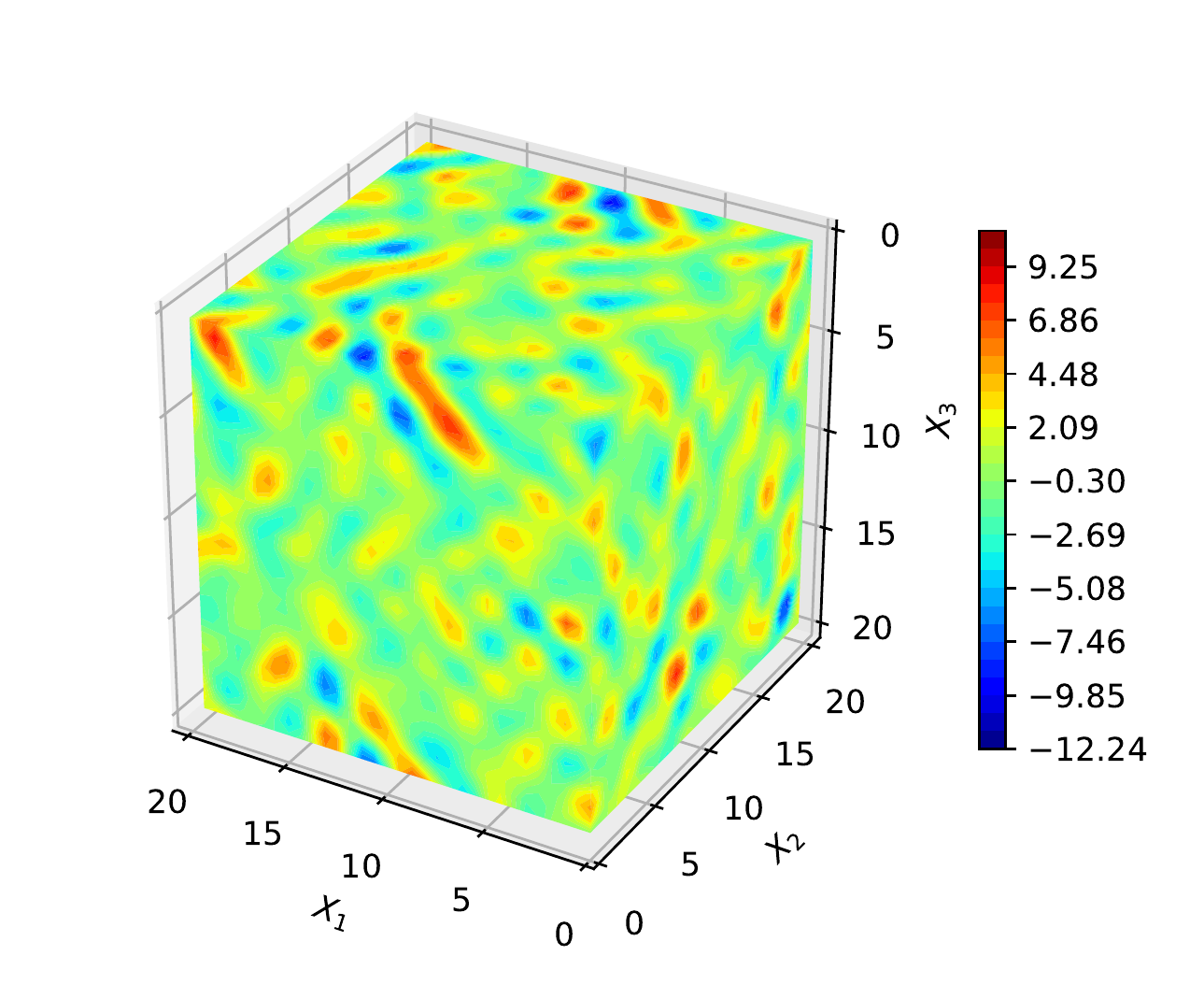}
  \caption{Difference between sample realizations from $2^{nd}$- and $3^{rd}$-order SRM simulations}
  \label{fig:example_3_comparision}
\end{figure}

As in the 2-dimensional case, the sample realisations look similar. The difference between the $2^{nd}$- and $3^{rd}$-order sample realisations is shown in Figure \ref{fig:example_3_comparision}. This difference is the result of the asymmetric non-Gaussianity introduced by the bispectrum. Here, similar to example 1, the asymmetric features in the difference plot are inclined along a $45^{\circ}$ and $-45^{\circ}$ angle along on each plane ($x_1-x_2$, $x_1-x_3$, and $x_2-x_3$) of the sample realisation. The uniformity of the bispectrum across all frequencies gives rise to this.

Sample statistics are given in Table \ref{table:example_3} from the 1000 simulations, which demonstrates the ability of the $3^{rd}-order$ simulations to match the moments up to the skewness. 
\begin{table}[h]
\centering
\begin{tabular}{l l l l}
\hline
\textbf{Moments} & \textbf{Target} & \textbf{$3^{rd}$-order} & \textbf{$2^{nd}$-order}\\
\hline
Mean & 0.00 & 0.0364 & 0.0364 \\
Variance & 179.0812 & 178.9807 & 178.9271\\
Skewness & 0.02107 & 0.02205 & 0.00081\\
\hline
\end{tabular}
\caption{Target and estimated moments of random fields generated by the 2nd and 3rd order SRM}
\label{table:example_3}
\end{table}
The samples also possess the prescribed bispectrum, but it is not feasible to illustrate this.

\subsection{3-dimensional random fields with different bispectra}

In this final example, we further investigate the effects of variations in the bispectrum in 3-dimensional random fields. The random fields simulated here possess the power spectrum from Eq.\ \eqref{fig:example_3_ps} and illustrated in Figure \ref{fig:example_3_ps}. We then generate $3^{rd}$-order random fields with 3 different bispectra given by
\begin{equation}
\begin{aligned}
    &\Re B_{1}(\omega_{11}, \omega_{12}, \omega_{21}, \omega_{22}, \omega_{31}, \omega_{32}) = \Im  B_{1}(\omega_{11}, \omega_{12}, \omega_{21}, \omega_{22}, \omega_{31}, \omega_{32}) \\
    &= \frac{300}{2\pi}\exp{-(10\omega_{11}^{2} + \omega_{12}^{2} + \omega_{21}^{2} + 10\omega_{22}^{2} + \omega_{31}^{2} + \omega_{32}^{2})}
\end{aligned}
\end{equation}
\begin{equation}
\begin{aligned}
    &\Re B_{2}(\omega_{11}, \omega_{12}, \omega_{21}, \omega_{22}, \omega_{31}, \omega_{32}) = \Im  B_{2}(\omega_{11}, \omega_{12}, \omega_{21}, \omega_{22}, \omega_{31}, \omega_{32}) \\
    &= \frac{300}{2\pi}\exp{-(\omega_{11}^{2} + 10\omega_{12}^{2} + \omega_{21}^{2} + \omega_{22}^{2} + 10\omega_{31}^{2} + \omega_{32}^{2})}
\end{aligned}
\end{equation}
\begin{equation}
\begin{aligned}
    &\Re B_{3}(\omega_{11}, \omega_{12}, \omega_{21}, \omega_{22}, \omega_{31}, \omega_{32}) = \Im  B_{3}(\omega_{11}, \omega_{12}, \omega_{21}, \omega_{22}, \omega_{31}, \omega_{32}) \\
    &= \frac{300}{2\pi}\exp{-(\omega_{11}^{2} + \omega_{12}^{2} + 10\omega_{21}^{2} + \omega_{22}^{2} + \omega_{31}^{2} + 10\omega_{32}^{2})}
\end{aligned}
\end{equation}

\begin{figure}[!ht]
\centering
\begin{subfigure}{.3\textwidth}
  \centering
  \includegraphics[width=\linewidth]{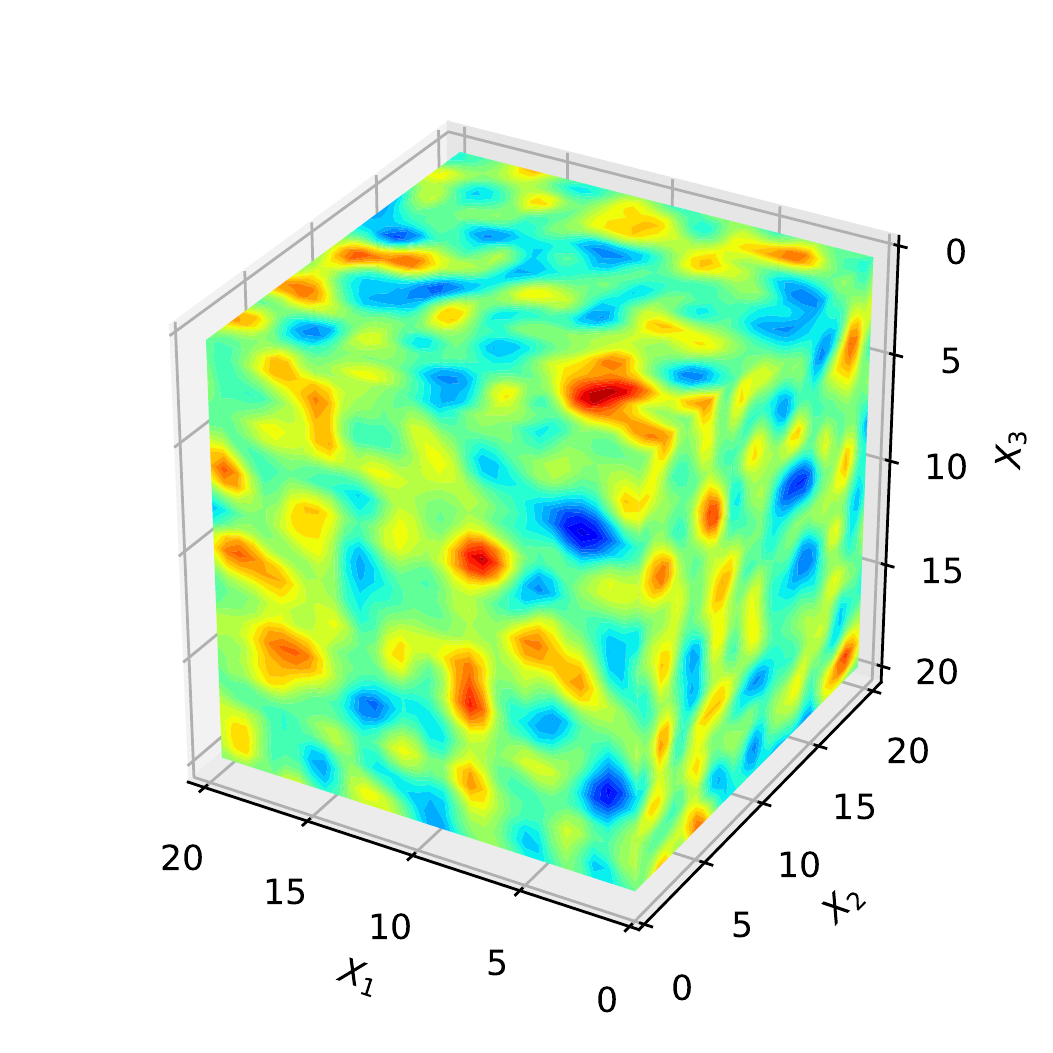}
  \caption{Bispectrum $B_1$}
  \label{fig:example_4_BSRM_1_sample}
\end{subfigure}%
\begin{subfigure}{.3\textwidth}
  \centering
  \includegraphics[width=\linewidth]{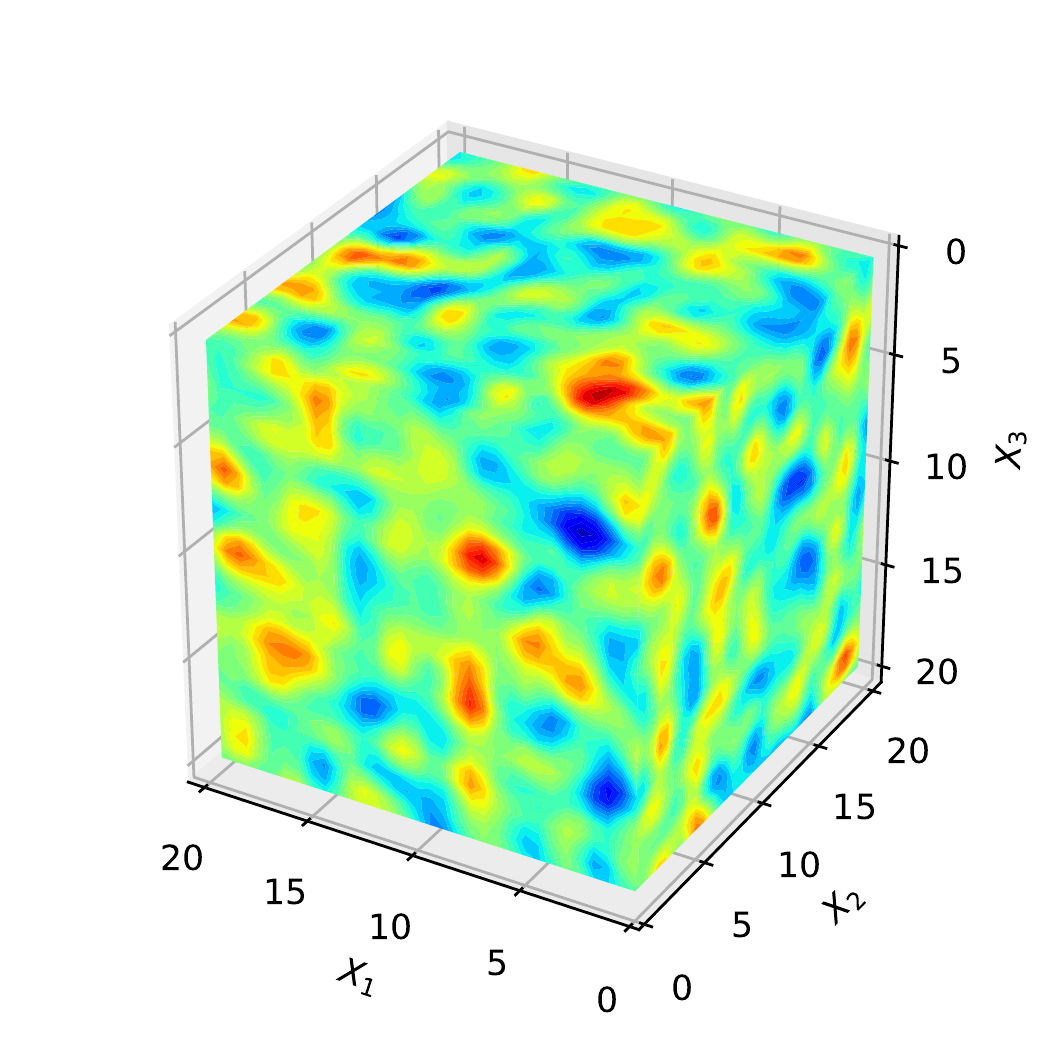}
  \caption{Bispectrum $B_2$}
  \label{fig:example_4_BSRM_2_sample}
\end{subfigure}
\begin{subfigure}{.36\textwidth}
  \centering
  \includegraphics[width=\linewidth]{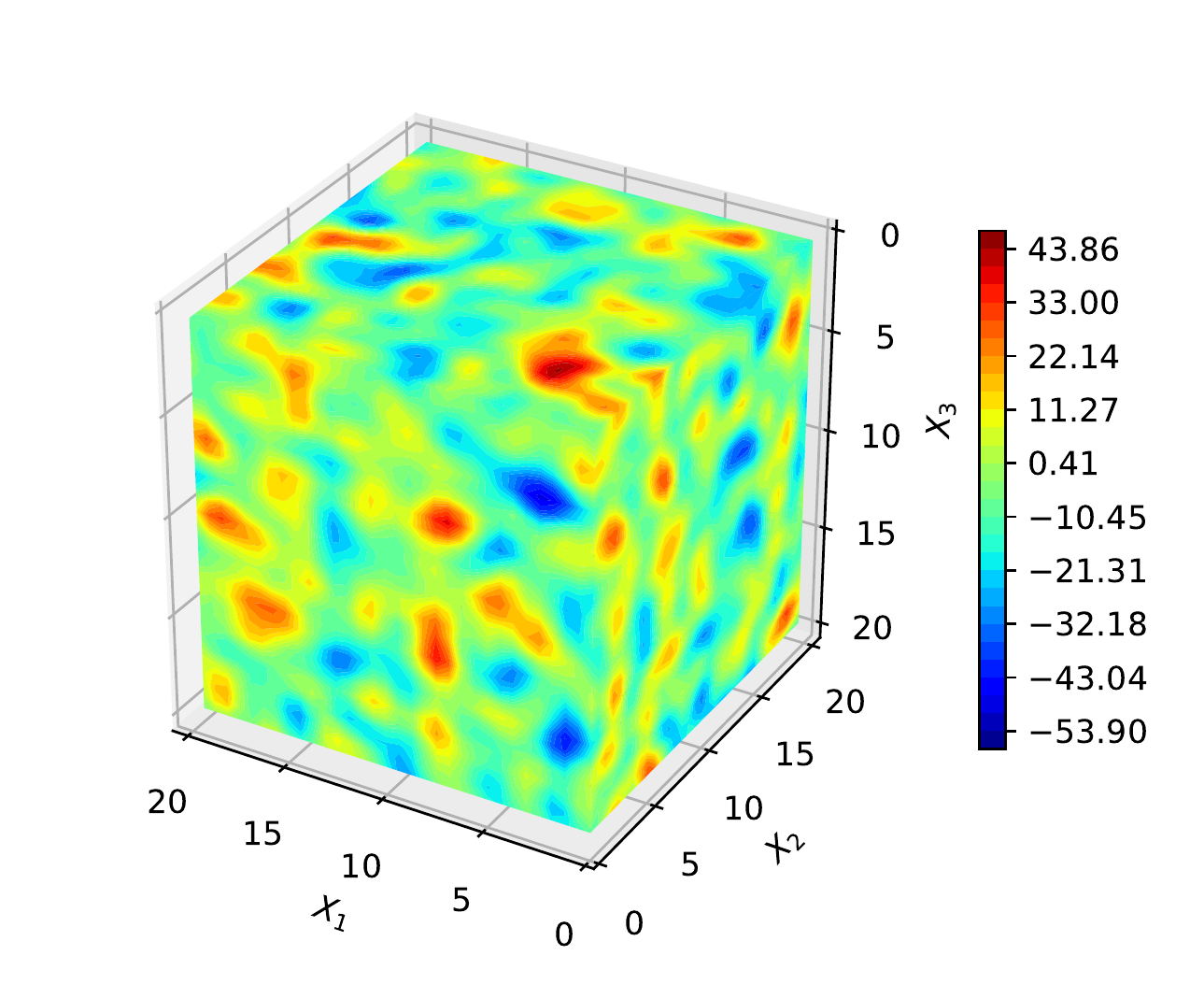}
  \caption{Bispectrum $B_3$}
  \label{fig:example_4_BSRM_3_sample}
\end{subfigure}
\caption{3-dimensional random fields generated using 3 different bispectra.}
\label{fig:example_4_SR}
\end{figure}
\begin{figure}[!ht]
\centering
\begin{subfigure}{.3\textwidth}
  \centering
  \includegraphics[width=\linewidth]{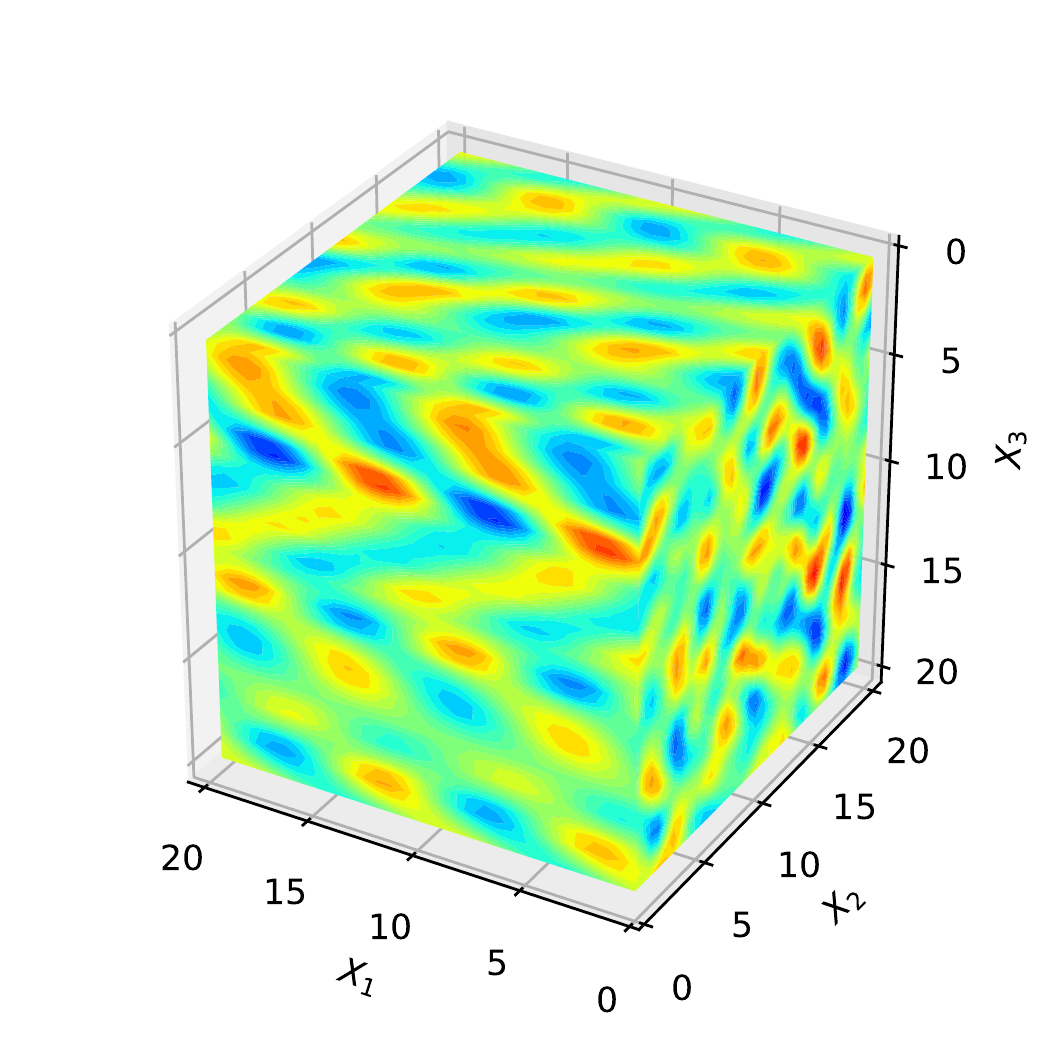}
  \caption{Bispectrum $B_1$}
  \label{fig:example_4_BSRM_1_SRM}
\end{subfigure}%
\begin{subfigure}{.3\textwidth}
  \centering
  \includegraphics[width=\linewidth]{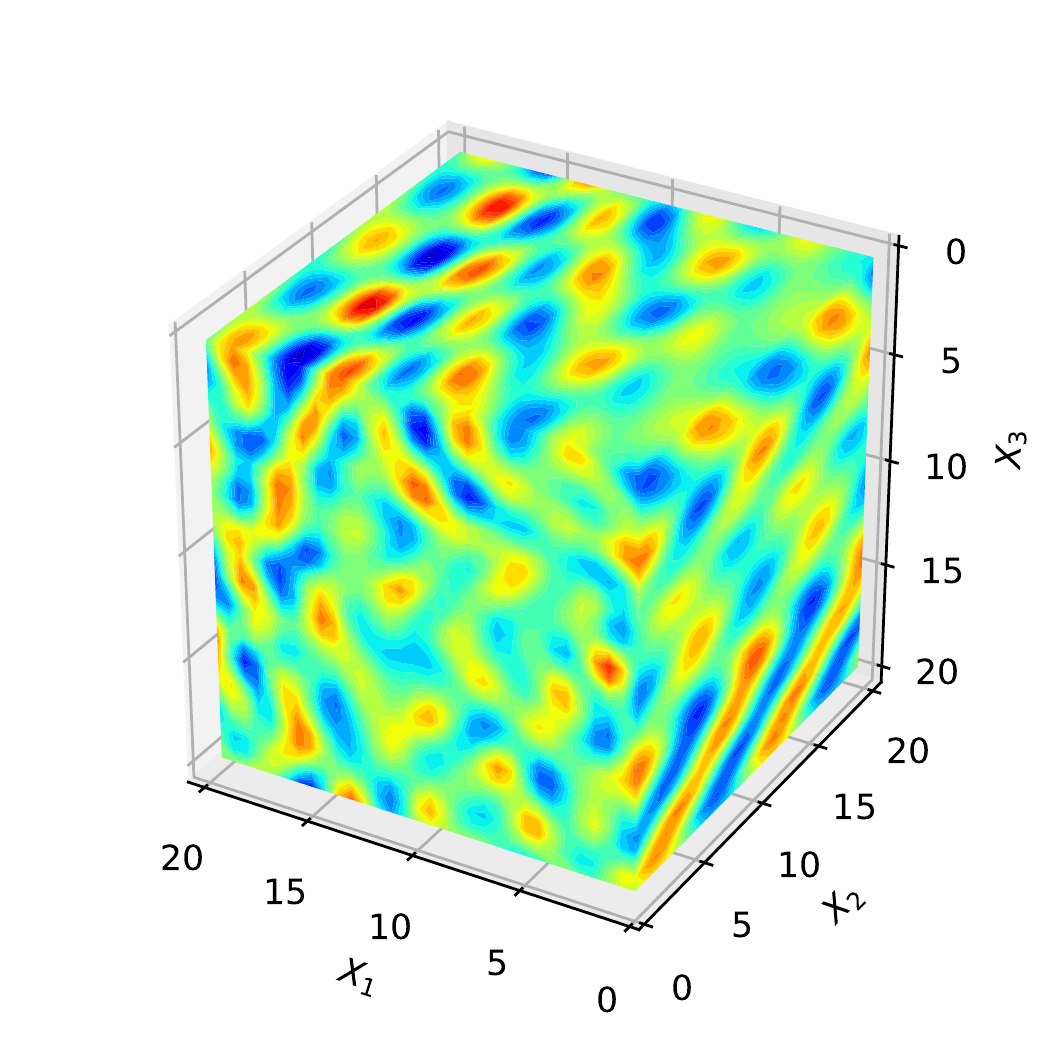}
  \caption{Bispectrum $B_2$}
  \label{fig:example_4_BSRM_2_SRM}
\end{subfigure}
\begin{subfigure}{.36\textwidth}
  \centering
  \includegraphics[width=\linewidth]{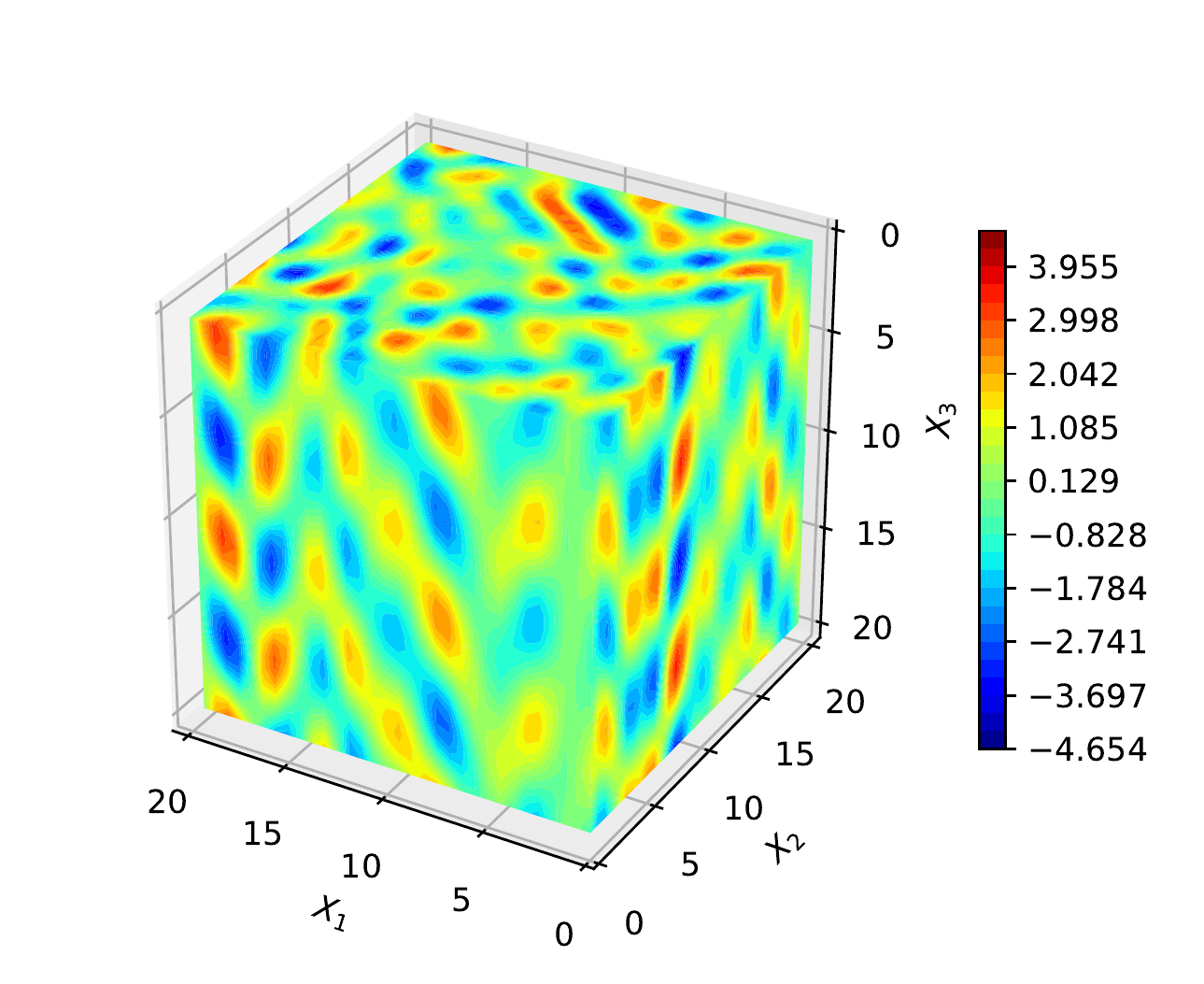}
  \caption{Bispectrum $B_3$}
  \label{fig:example_4_BSRM_3_SRM}
\end{subfigure}
\caption{Difference between 3-dimensional random fields generated using the $2^{nd}$- and $3^{rd}$-order spectral representation methods.}
\label{fig:example_4_comparision}
\end{figure}

Plots of the sample realisations from the $3^{rd}$-order SRM, having identical phase angles as those in the previous example, are presented in Figure \ref{fig:example_4_SR}.
As in previous examples, the random field realizations look very similar. Figure \ref{fig:example_4_comparision} shows the difference between these samples and the $2^{nd}$-order field simulated in Figure \ref{fig:example_3_SRM_sample}.
Here we see that by taking the difference between the samples generated by the $2^{nd}$- and $3^{rd}$-order Spectral Representation Methods, we have asymmetric features elongated in particular axes. Specifically, for realizations with bispectrum $B_1$, the length-scale of asymmetric features is much smaller in the $x_1$-axis as in the $x_2$, $x_3$ axes. Moreover, we see that the asymmetric features lie along different angles on different planes. In the $x_1-x_3$ and $x_1-x_2$ planes, the asymmetric features lie along angles of $\arctan(\sqrt{10})\approx 73^{\circ}$ and $\arctan(-\sqrt{10})\approx -73^{\circ}$ relative to the $x_1$ axis. Meanwhile, on the $x_2-x_3$ plane corresponding to a plane with equal bispectral length-scales, the asymmetric features lie at approximately $\arctan(1)=45^{\circ}$. The same observations can be made for the samples from random fields with bispectra $B_2$ and $B_3$. In all cases, the length-scale of the asymmetric features follows directly from the form of the bispectrum and the angle at which these features lie on a given plane relates to the arctan of the relative length scales of the bispectrum in that plane.


Lastlly, 1000 samples with the following discretization were simulated
\begin{equation}
\begin{aligned}
	&\Delta x_{1} = \Delta x_{2} = \Delta x_{3} = 0.625\\
	&\Delta \kappa_{1} = \Delta \kappa_{2}  = \Delta \kappa_{2} = 0.314\\
	& N_{1} = N_{2} = N_{3} = 16\\
	& M_{1} = M_{2} = M_{3} = 32\\
\end{aligned}
\end{equation}
and the statistics of the resulting random fields were calculated as shown in Table \ref{table:example_4}.
\begin{table}[!ht]
\centering
\begin{tabular}{l l l l l l}
\hline
\textbf{Moments} & \textbf{Target} & \textbf{$3^{rd}$-order, $B_1$} & \textbf{$3^{rd}$-order, $B_2$} & \textbf{$3^{rd}$-order, $B_3$} & \textbf{$2^{nd}$-order}\\
\hline
Mean & 0.00 & 0.0364 & 0.0364 & 0.0364 & 0.0364\\
Variance & 179.0812 & 178.9703 & 178.9787 & 178.9605 & 178.9270\\
Skewness & 0.00580 & 0.00680 & 0.00682 & 0.00661 & 0.0008\\
\hline
\end{tabular}
\caption{Target and estimated moments of random fields generated by the 2nd and 3rd order SRM}
\label{table:example_4}
\end{table}
Again, the third-order samples are shown to possess the appropriate $2^{nd}$- and $3^{rd}$-order statistics. While they also possess the proper bispectra, this cannot be feasibly illustrated.

\section{Conclusions}

In this paper, the $3^{rd}$-order Spectral Representation Method has been extended for the simulation of multi-dimensional random fields. This simulation formula has been derived for 2-dimension, 3-dimensional, and general $d$-dimensional fields. A fast Fourier transform implementation of the $3^{rd}$-order SRM has also presented, which leads to enormous computational gains -- making the generation of 2D and 3D fields feasible for implementation on a desktop computer. Numerical examples of 2D and 3D random fields are provided, which highlight the effectiveness of the proposed methodology and some interesting features associated with the asymmetries of the generated random fields.

\section{Acknowledgements}
This work has been supported by the National Science Foundation under award number 1652044.





\bibliographystyle{model1-num-names}
\bibliography{elsarticle-template-1-num.bib}







\newpage
\appendix

\section{Mean, 2-, and 3-point correlations functions of the simulated random fields}
\label{A:1}

Here, we show that the generated random fields possess the proper mean value, 2-point, and 3-point autocorrelation functions.

\subsection{Mean value}


The expected value of the real random field $A(x_1,x_2)$ expressed in terms of the Cramer spectral representation in Eq.\ \eqref{eqn:Cramer_real} is given below. Applying the orthogonality conditions in \eqref{eqn:1_order_orthogonality} gives the following.
\begin{equation}
\begin{aligned}
    & \mathbb{E}[A(x_{1}, x_{2})] \\
    & = \mathbb{E}[\int_{-\infty}^{\infty}\int_{0}^{\infty}[\cos(\kappa_{1}x_{1} + \kappa_{2}x_{2})du(\kappa_{1}, \kappa_{2}) + \sin(\kappa_{1}x_{1} + \kappa_{2}x_{2})dv(\kappa_{1}, \kappa_{2})]]\\
    & = \int_{-\infty}^{\infty}\int_{0}^{\infty}[\cos(\kappa_{1}x_{1} + \kappa_{2}x_{2})\mathbb{E}[du(\kappa_{1}, \kappa_{2})] + \sin(\kappa_{1}x_{1} + \kappa_{2}x_{2})\mathbb{E}[dv(\kappa_{1}, \kappa_{2})]]\\
    & = 0
\end{aligned}
\end{equation}
Hence, if the selected orthogonal increments possess the first-order orthogonality condition, the random field will have zero mean.

\subsection{2-point correlation}
Let us define the two-point correlation function of $A(x_1,x_2)$ as below. Applying the second-order orthogonality conditions in Eqs.\ \eqref{eqn:2_order_orthogonal_increments} and basic trigonometric identities yields:
\begin{equation}
\begin{aligned}
    & \mathbb{E}[A(x_{1}, x_{2})A(x_{1} + \xi_{1}, x_{2} + \xi_{2})]\\
    & = \mathbb{E}[\int_{-\infty}^{\infty}\int_{0}^{\infty}[\cos(\kappa_{1}x_{1} + \kappa_{2}x_{2})du(\kappa_{1}, \kappa_{2}) + \sin(\kappa_{1}x_{1} + \kappa_{2}x_{2})dv(\kappa_{1}, \kappa_{2})]\\
    & \int_{-\infty}^{\infty}\int_{0}^{\infty}[\cos(\kappa'_{1}(x_{1} + \xi_{1}) + \kappa'_{2}(x_{2} + \xi_{2}))du(\kappa'_{1}, \kappa'_{2})\\
    & + \sin(\kappa'_{1}(x_{1} + \xi_{1}) + \kappa'_{2}(x_{2} + \xi_{2}))dv(\kappa'_{1}, \kappa'_{2})]]\\
    & = \int_{-\infty}^{\infty}\int_{0}^{\infty}\int_{-\infty}^{\infty}\int_{0}^{\infty}\\
    & [\cos(\kappa_{1}x_{1} + \kappa_{2}x_{2})\cos(\kappa'_{1}(x_{1} + \xi_{1}) + \kappa'_{2}(x_{2} + \xi_{2}))\mathbb{E}[du(\kappa_{1}, \kappa_{2})du(\kappa'_{1}, \kappa'_{2})]\\
    & + \cos(\kappa_{1}x_{1} + \kappa_{2}x_{2})\sin(\kappa'_{1}(x_{1} + \xi_{1}) + \kappa'_{2}(x_{2} + \xi_{2}))\mathbb{E}[du(\kappa_{1}, \kappa_{2})dv(\kappa'_{1}, \kappa'_{2})]\\
    & + \sin(\kappa_{1}x_{1} + \kappa_{2}x_{2})\cos(\kappa'_{1}(x_{1} + \xi_{1}) + \kappa'_{2}(x_{2} + \xi_{2}))\mathbb{E}[dv(\kappa_{1}, \kappa_{2})du(\kappa'_{1}, \kappa'_{2})]\\
    & + \sin(\kappa_{1}x_{1} + \kappa_{2}x_{2})\sin(\kappa'_{1}(x_{1} + \xi_{1}) + \kappa'_{2}(x_{2} + \xi_{2}))\mathbb{E}[dv(\kappa_{1}, \kappa_{2})dv(\kappa'_{1}, \kappa'_{2})]]\\
    & = \int_{-\infty}^{\infty}\int_{0}^{\infty}[\cos(\kappa_{1}x_{1} + \kappa_{2}x_{2})\cos(\kappa_{1}(x_{1} + \xi_{1}) + \kappa_{2}(x_{2} + \xi_{2}))S_{1}(\kappa_{1}, \kappa_{1})d\kappa_{1}d\kappa_{2} + \\
    & \sin(\kappa_{1}x_{1} + \kappa_{2}x_{2})\sin(\kappa_{1}(x_{1} + \xi_{1}) + \kappa_{2}(x_{2} + \xi_{2}))S_{1}(\kappa_{1}, \kappa_{1})d\kappa_{1}d\kappa_{2}] \\
    & = \int_{-\infty}^{\infty}\int_{0}^{\infty}\cos(\kappa_{1}\xi_{1} + \kappa_{2}\xi_{2})S_{1}(\kappa_{1}, \kappa_{1})d\kappa_{1}d\kappa_{2} \\
    & = \int_{-\infty}^{\infty}\int_{-\infty}^{\infty}e^{\iota(\kappa_{1}\xi_{1} + \kappa_{2}\xi_{2})}S(\kappa_{1}, \kappa_{1})d\kappa_{1}d\kappa_{2} \\ 
    & = R_{2}(\xi_{1}, \xi_{2})\\
\end{aligned}
\end{equation}

It should be noted here that the integrals $\int_{0}^{\infty}$ correspond to $\kappa_{1}$ and $\kappa'_{1}$ whereas the integrals $\int_{-\infty}^{\infty}$ correspond to $\kappa_{2}$ and $\kappa'_{2}$. Therefore, if the proposed orthogonal increments satisfy the second-order orthogonality conditions, the random field will possess the appropriate 2-point correlation function.

\subsection{3-point correlation}
Finally, the $3$-point autocorrelation function can be derived as follows for the Cramer spectral representation:

\begin{equation}
\begin{aligned}
    & \mathbb{E}[A(x_{1}, x_{2})A(x_{1} + \xi_{11}, x_{2} + \xi_{21})A(x_{1} + \xi_{12}, x_{2} + \xi_{22})]\\
    & = \mathbb{E}[\int_{-\infty}^{\infty}\int_{0}^{\infty}[\cos(\kappa_{1}x_{1} + \kappa_{2}x_{2})du(\kappa_{1}, \kappa_{2}) + \sin(\kappa_{1}x_{1} + \kappa_{2}x_{2})dv(\kappa_{1}, \kappa_{2})]\\
    & \int_{-\infty}^{\infty}\int_{0}^{\infty}[\cos(\kappa'_{1}(x_{1} + \xi_{11}) + \kappa'_{2}(x_{2} + \xi_{21}))du(\kappa'_{1}, \kappa'_{2})+ \sin(\kappa'_{1}(x_{1} + \xi_{11}) + \kappa'_{2}(x_{2} + \xi_{21}))dv(\kappa'_{1}, \kappa'_{2})]\\
    & \int_{-\infty}^{\infty}\int_{0}^{\infty}[\cos(\kappa''_{1}(x_{1} + \xi_{12}) + \kappa''_{2}(x_{2} + \xi_{22}))du(\kappa''_{1}, \kappa''_{2}) + \sin(\kappa''_{1}(x_{1} + \xi_{12}) + \kappa''_{2}(x_{2} + \xi_{22}))dv(\kappa''_{1}, \kappa''_{2})]]\\
    & = \int_{-\infty}^{\infty}\int_{0}^{\infty}\int_{-\infty}^{\infty}\int_{0}^{\infty}\int_{-\infty}^{\infty}\int_{0}^{\infty}\\
    & [\cos(\kappa_{1}x_{1} + \kappa_{2}x_{2})\cos(\kappa'_{1}(x_{1} + \xi_{11}) + \kappa'_{2}(x_{2} + \xi_{21}))\cos(\kappa''_{1}(x_{1} + \xi_{12}) + \kappa''_{2}(x_{2} + \xi_{22}))\\
    & \mathbb{E}[du(\kappa_{1}, \kappa_{2})du(\kappa'_{1}, \kappa'_{2})du(\kappa''_{1}, \kappa''_{2})]\\
    & + \cos(\kappa_{1}x_{1} + \kappa_{2}x_{2})\cos(\kappa'_{1}(x_{1} + \xi_{11}) + \kappa'_{2}(x_{2} + \xi_{21}))\sin(\kappa''_{1}(x_{1} + \xi_{12}) + \kappa''_{2}(x_{2} + \xi_{22}))\\
    & \mathbb{E}[du(\kappa_{1}, \kappa_{2})du(\kappa'_{1}, \kappa'_{2})dv(\kappa''_{1}, \kappa''_{2})]\\
    & + \cos(\kappa_{1}x_{1} + \kappa_{2}x_{2})\sin(\kappa'_{1}(x_{1} + \xi_{11}) + \kappa'_{2}(x_{2} + \xi_{21}))\cos(\kappa''_{1}(x_{1} + \xi_{12}) + \kappa''_{2}(x_{2} + \xi_{22}))\\
    & \mathbb{E}[du(\kappa_{1}, \kappa_{2})dv(\kappa'_{1}, \kappa'_{2})du(\kappa''_{1}, \kappa''_{2})]\\
    & + \cos(\kappa_{1}x_{1} + \kappa_{2}x_{2})\sin(\kappa'_{1}(x_{1} + \xi_{11}) + \kappa'_{2}(x_{2} + \xi_{21}))\sin(\kappa''_{1}(x_{1} + \xi_{12}) + \kappa''_{2}(x_{2} + \xi_{22}))\\
    & \mathbb{E}[du(\kappa_{1}, \kappa_{2})dv(\kappa'_{1}, \kappa'_{2})dv(\kappa''_{1}, \kappa''_{2})]\\
    & + \sin(\kappa_{1}x_{1} + \kappa_{2}x_{2})\cos(\kappa'_{1}(x_{1} + \xi_{11}) + \kappa'_{2}(x_{2} + \xi_{21}))\cos(\kappa''_{1}(x_{1} + \xi_{12}) + \kappa''_{2}(x_{2} + \xi_{22}))\\
    & \mathbb{E}[dv(\kappa_{1}, \kappa_{2})du(\kappa'_{1}, \kappa'_{2})du(\kappa''_{1}, \kappa''_{2})]\\
    & + \sin(\kappa_{1}x_{1} + \kappa_{2}x_{2})\cos(\kappa'_{1}(x_{1} + \xi_{11}) + \kappa'_{2}(x_{2} + \xi_{21}))\sin(\kappa''_{1}(x_{1} + \xi_{12}) + \kappa''_{2}(x_{2} + \xi_{22}))\\
    & \mathbb{E}[dv(\kappa_{1}, \kappa_{2})du(\kappa'_{1}, \kappa'_{2})dv(\kappa''_{1}, \kappa''_{2})]\\
    & + \sin(\kappa_{1}x_{1} + \kappa_{2}x_{2})\sin(\kappa'_{1}(x_{1} + \xi_{11}) + \kappa'_{2}(x_{2} + \xi_{21}))\cos(\kappa''_{1}(x_{1} + \xi_{12}) + \kappa''_{2}(x_{2} + \xi_{22}))\\
    & \mathbb{E}[dv(\kappa_{1}, \kappa_{2})dv(\kappa'_{1}, \kappa'_{2})du(\kappa''_{1}, \kappa''_{2})]\\
    & + \sin(\kappa_{1}x_{1} + \kappa_{2}x_{2})\sin(\kappa'_{1}(x_{1} + \xi_{11}) + \kappa'_{2}(x_{2} + \xi_{21}))\sin(\kappa''_{1}(x_{1} + \xi_{12}) + \kappa''_{2}(x_{2} + \xi_{22}))\\
    & \mathbb{E}[dv(\kappa_{1}, \kappa_{2})dv(\kappa'_{1}, \kappa'_{2})dv(\kappa''_{1}, \kappa''_{2})]]\\
\end{aligned}
\end{equation}

Using the orthogonality conditions in Eq.\ \eqref{eqn:3_order_orthogonal_increments}, we have
\begin{equation}
\begin{aligned}
	& \mathbb{E}[A(x_{1}, x_{2})A(x_{1} + \xi_{11}, x_{2} + \xi_{21})A(x_{1} + \xi_{12}, x_{2} + \xi_{22})] = \\
	& \int_{-\infty}^{\infty}\int_{0}^{\infty}\int_{-\infty}^{\infty}\int_{0}^{\infty}\\
	& \cos(\kappa_{1}x_{1} + \kappa_{2}x_{2})\cos((\kappa'_{1} + \kappa''_{1})x_{1} + (\kappa'_{2} + \kappa''_{2})x_{2} + \kappa'_{1}\xi_{11} + \kappa''_{1}\xi_{21} + \kappa'_{2}\xi_{12} + \kappa''_{2}\xi_{22})]d\Re B_{1}(\kappa'_{1}, \kappa'_{2}, \kappa'_{11}, \kappa'_{12})\\
	& - \cos(\kappa_{1}x_{1} + \kappa_{2}x_{2})\sin((\kappa'_{1} + \kappa''_{1})x_{1} + (\kappa'_{2} + \kappa''_{2})x_{2} + \kappa'_{1}\xi_{11} + \kappa''_{1}\xi_{21} + \kappa'_{2}\xi_{12} + \kappa''_{2}\xi_{22})]d\Im B_{1}(\kappa'_{1}, \kappa'_{2}, \kappa'_{11}, \kappa'_{12})\\
	& + \sin(\kappa_{1}x_{1} + \kappa_{2}x_{2})\cos((\kappa'_{1} + \kappa''_{1})x_{1} + (\kappa'_{2} + \kappa''_{2})x_{2} + \kappa'_{1}\xi_{11} + \kappa''_{1}\xi_{21} + \kappa'_{2}\xi_{12} + \kappa''_{2}\xi_{22})]d\Im B_{1}(\kappa'_{1}, \kappa'_{2}, \kappa'_{11}, \kappa'_{12})\\
	& + \sin(\kappa_{1}x_{1} + \kappa_{2}x_{2})\sin((\kappa'_{1} + \kappa''_{1})x_{1} + (\kappa'_{2} + \kappa''_{2})x_{2} + \kappa'_{1}\xi_{11} + \kappa''_{1}\xi_{21} + \kappa'_{2}\xi_{12} + \kappa''_{2}\xi_{22})]d\Re B_{1}(\kappa'_{1}, \kappa'_{2}, \kappa'_{11}, \kappa'_{12})\\
\end{aligned}
\end{equation}
This can be further expressed as:
\begin{equation}
\begin{aligned}
	& \mathbb{E}[A(x_{1}, x_{2})A(x_{1} + \xi_{11}, x_{2} + \xi_{21})A(x_{1} + \xi_{12}, x_{2} + \xi_{22})] = \\
	& \int_{-\infty}^{\infty}\int_{0}^{\infty}\int_{-\infty}^{\infty}\int_{0}^{\infty}\\
	& \cos(\kappa'_{1}\xi_{11} + \kappa''_{1}\xi_{21} + \kappa'_{2}\xi_{12} + \kappa''_{2}\xi_{22})]\Re B_{1}(\kappa'_{1}, \kappa'_{2}, \kappa'_{11}, \kappa'_{12})\\
	& - \sin(\kappa'_{1}\xi_{11} + \kappa''_{1}\xi_{21} + \kappa'_{2}\xi_{12} + \kappa''_{2}\xi_{22})]\Im B_{1}(\kappa'_{1}, \kappa'_{2}, \kappa'_{11}, \kappa'_{12})\\
    & = \int_{-\infty}^{\infty}\int_{-\infty}^{\infty}\int_{-\infty}^{\infty}\int_{-\infty}^{\infty}e^{\iota(\kappa_{11}\xi_{11} + \kappa_{21}\xi_{21} + \kappa_{12}\xi_{12} + \kappa_{22}\xi_{22})}B(\kappa_{11}, \kappa_{21}, \kappa_{12}, \kappa_{22})d\kappa_{11}d\kappa_{21}d\kappa_{12}d\kappa_{22}\\
    & = R_{3}(\xi_{11}, \xi_{21}, \xi_{12}, \xi_{22})\\
\end{aligned}
\end{equation}
Here the integrals $\int_{0}^{\infty}$ correspond to $\kappa'_{1}$ and $\kappa''_{1}$ whereas the integrals $\int_{-\infty}^{\infty}$ correspond to $\kappa'_{2}$ and $\kappa''_{2}$. Again, if the proposed orthogonal increments satisfy the third-order orthogonality conditions, then the random field will possess the specified 3-point correlation function.

\section{Orthogonality proofs}
\label{A:2}

\begin{equation}
\begin{aligned}
    &\mathbb{E}[du(\kappa_{1}, \kappa_{2})] = \mathbb{E}[\sqrt{2}A_{pn_{1}n_{2}}\cos\Phi_{n_{1}n_{2}}\\
    & + \sum_{i_{1} + j_{1} = n_{1}}^{i_{1} \geq j_{1} \geq 0}\sum_{i_{2} + j_{2} = n_{2}}^{|n_{2}| \geq |i_{2}| \geq |j_{2}| \geq 0}\sqrt{2}A_{n_{1}n_{2}}b_{p}(\kappa_{1i_{1}}, \kappa_{1j_{1}}, \kappa_{2i_{2}}, \kappa_{2j_{2}})\cos(\Phi_{i_{1}i_{2}} + \Phi_{j_{1}j_{2}} + \beta(\kappa_{1i_{1}}, \kappa_{1j_{1}}, \kappa_{2i_{2}}, \kappa_{2j_{2}}))]\\
    & = \mathbb{E}[\sqrt{2}A_{pn_{1}n_{2}}\cos\Phi_{n_{1}n_{2}}]\\
    & + \mathbb{E}[\sum_{i_{1} + j_{1} = n_{1}}^{i_{1} \geq j_{1} \geq 0}\sum_{i_{2} + j_{2} = n_{2}}^{|n_{2}| \geq |i_{2}| \geq |j_{2}| \geq 0}\sqrt{2}A_{n_{1}n_{2}}b_{p}(\kappa_{1i_{1}}, \kappa_{1j_{1}}, \kappa_{2i_{2}}, \kappa_{2j_{2}})\cos(\Phi_{i_{1}i_{2}} + \Phi_{j_{1}j_{2}} + \beta(\kappa_{1i_{1}}, \kappa_{1j_{1}}, \kappa_{2i_{2}}, \kappa_{2j_{2}}))]\\
    & = 0 + 0\\
\end{aligned}
\end{equation}

\begin{equation}
\begin{aligned}
    &\mathbb{E}[du^{2}(\kappa_{1}, \kappa_{2})] = \mathbb{E}[(\sqrt{2}A_{pn_{1}n_{2}}\cos\Phi_{n_{1}n_{2}}\\
    & + \sum_{i_{1} + j_{1} = n_{1}}^{i_{1} \geq j_{1} \geq 0}\sum_{i_{2} + j_{2} = n_{2}}^{|n_{2}| \geq |i_{2}| \geq |j_{2}| \geq 0}\sqrt{2}A_{n_{1}n_{2}}b_{p}(\kappa_{1i_{1}}, \kappa_{1j_{1}}, \kappa_{2i_{2}}, \kappa_{2j_{2}})\cos(\Phi_{i_{1}i_{2}} + \Phi_{j_{1}j_{2}} + \beta(\kappa_{1i_{1}}, \kappa_{1j_{1}}, \kappa_{2i_{2}}, \kappa_{2j_{2}})))^{2}]\\
    & = \mathbb{E}[(\sqrt{2}A_{pn_{1}n_{2}}\cos\Phi_{n_{1}n_{2}})^{2}]\\
    & + \mathbb{E}[(\sum_{i_{1} + j_{1} = n_{1}}^{i_{1} \geq j_{1} \geq 0}\sum_{i_{2} + j_{2} = n_{2}}^{|n_{2}| \geq |i_{2}| \geq |j_{2}| \geq 0}\sqrt{2}A_{n_{1}n_{2}}b_{p}(\kappa_{1i_{1}}, \kappa_{1j_{1}}, \kappa_{2i_{2}}, \kappa_{2j_{2}})\cos(\Phi_{i_{1}i_{2}} + \Phi_{j_{1}j_{2}} + \beta(\kappa_{1i_{1}}, \kappa_{1j_{1}}, \kappa_{2i_{2}}, \kappa_{2j_{2}})))^{2}]\\
    & + \mathbb{E}[4A^{2}_{pn_{1}n_{2}}\cos^{2}\Phi_{n_{1}n_{2}}\\
    & \sum_{i_{1} + j_{1} = n_{1}}^{i_{1} \geq j_{1} \geq 0}\sum_{i_{2} + j_{2} = n_{2}}^{
|i_{2}| \geq |j_{2}| \geq 0}A^{2}_{n_{1}n_{2}}b^{2}_{p}(\kappa_{1i_{1}}, \kappa_{1j_{1}}, \kappa_{2i_{2}}, \kappa_{2j_{2}})\cos^{2}(\Phi_{i_{1}i_{2}} + \Phi_{j_{1}j_{2}} + \beta(\kappa_{1i_{1}}, \kappa_{1j_{1}}, \kappa_{2i_{2}}, \kappa_{2j_{2}}))]\\
    & = 2\Big(S_{1p}(\kappa_{1n_{1}}, \kappa_{2n_{2}})\Delta\kappa_{1}\Delta\kappa_{2} + \sum_{i_{1} + j_{1} = n_{1}}^{i_{1} \geq j_{1} \geq 0}\sum_{i_{2} + j_{2} = n_{2}}^{|i_{2}| \geq |j_{2}| \geq 0}S_{1}(\kappa_{1n_{1}}, \kappa_{2n_{2}})b_{p}^{2}(\kappa_{1i_{1}}, \kappa_{1j_{1}}, \kappa_{2i_{2}}, \kappa_{2j_{2}})\Delta\kappa_{1}\Delta\kappa_{2}\Big)\\
    & = S(\kappa_{1n_{1}}, \kappa_{2n_{2}})\Delta\kappa_{1}\Delta\kappa_{2}\\
\end{aligned}
\end{equation}

\begin{equation}
\begin{aligned}
    &\mathbb{E}[du(\kappa_{1}, \kappa_{2})du(\kappa'_{1}, \kappa'_{2})du(\kappa''_{1}, \kappa''_{2})] = \mathbb{E}[(\sqrt{2}A_{pn_{1}n_{2}}\cos\Phi_{n_{1}n_{2}}\\
    & + \sum_{i_{1} + j_{1} = n_{1}}^{i_{1} \geq j_{1} \geq 0}\sum_{i_{2} + j_{2} = n_{2}}^{|i_{2}| \geq |j_{2}| \geq 0}\sqrt{2}A_{n_{1}n_{2}}b_{p}(\kappa_{1i_{1}}, \kappa_{1j_{1}}, \kappa_{2i_{2}}, \kappa_{2j_{2}})\cos(\Phi_{i_{1}i_{2}} + \Phi_{j_{1}j_{2}} + \beta(\kappa_{1i_{1}}, \kappa_{1j_{1}}, \kappa_{2i_{2}}, \kappa_{2j_{2}})))\\
    & (\sqrt{2}A_{pn'_{1}n'_{2}}\cos\Phi_{n'_{1}n'_{2}}\\
    & + \sum_{i'_{1} + j'_{1} = n_{1}}^{i'_{1} \geq j'_{1} \geq 0}\sum_{i'_{2} + j'_{2} = n'_{2}}^{|i'_{2}| \geq |j'_{2}| \geq 0}\sqrt{2}A_{n'_{1}n'_{2}}b_{p}(\kappa_{1i'_{1}}, \kappa_{1j'_{1}}, \kappa_{2i'_{2}}, \kappa_{2j'_{2}})\cos(\Phi_{i'_{1}i'_{2}} + \Phi_{j'_{1}j'_{2}} + \beta(\kappa_{1i'_{1}}, \kappa_{1j'_{1}}, \kappa_{2i'_{2}}, \kappa_{2j'_{2}})))\\
    & (\sqrt{2}A_{pn''_{1}n''_{2}}\cos\Phi_{n''_{1}n''_{2}}\\
    & + \sum_{i''_{1} + j''_{1} = n''_{1}}^{i''_{1} \geq j''_{1} \geq 0}\sum_{i''_{2} + j''_{2} = n''_{2}}^{|i''_{2}| \geq |j''_{2}| \geq 0}\sqrt{2}A_{n''_{1}n''_{2}}b_{p}(\kappa_{1i''_{1}}, \kappa_{1j''_{1}}, \kappa_{2i''_{2}}, \kappa_{2j''_{2}})\cos(\Phi_{i''_{1}i''_{2}} + \Phi_{j''_{1}j''_{2}} + \beta(\kappa_{1i''_{1}}, \kappa_{1j''_{1}}, \kappa_{2i''_{2}}, \kappa_{2j''_{2}})))\\
    & = \mathbb{E}[(\sqrt{2}A_{pi_{1}i_{2}}\cos\Phi_{i_{1}i_{2}})(\sqrt{2}A_{pj_{1}j_{2}}\cos\Phi_{j_{1}j_{2}})\\
    & (\sqrt{2}A_{n_{1}n_{2}}b_{p}(\kappa_{1i_{1}}, \kappa_{1j_{1}}, \kappa_{2i_{2}}, \kappa_{2j_{2}})\cos(\Phi_{i_{1}i_{2}} + \Phi_{j_{1}j_{2}} + \beta(\kappa_{1i_{1}}, \kappa_{1j_{1}}, \kappa_{2i_{2}}, \kappa_{2j_{2}})))]\\
    & = B(\kappa_{1i_{1}}, \kappa_{1j_{1}}, \kappa_{2i_{2}}, \kappa_{2j_{2}})
\end{aligned}
\end{equation}

\end{document}